\documentclass[aos,seceqn,citesort,dvips]{arximspdf}

\doi{10.1214/09-AOS756}
\volume{38}
\issue{3}
\pubyear{2010}
\firstpage{1478}
\lastpage{1545}

\makeatletter

\newproclaim{example}{Example}[section]
\newproclaim{rem}[example]{Remark}

\newtheorem{theo}[example]{Theorem}
\newtheorem{lem}[example]{Lemma}
\newtheorem{prop}[example]{Proposition}

\newtheorem{Hypothesis}{Hypothesis}

\makeatother

\begin{document}
\begin{frontmatter}

\title{Limit theorems for moving averages of discretized processes plus noise}
\runtitle{Limit theorems for moving averages}

\begin{aug}
\author[A]{\fnms{Jean} \snm{Jacod}\ead[label=e1]{jean.jacod@upmc.fr}},
\author[B]{\fnms{Mark} \snm{Podolskij}\corref{}\ead[label=e2]{mark.podolskij@math.ethz.ch}} and
\author[C]{\fnms{Mathias} \snm{Vetter}\thanksref{t1}\ead[label=e3]{mathias.vetter@rub.de}}
\runauthor{J. Jacod, M. Podolskij and M. Vetter}
\affiliation{UPMC (Universit\'e Paris-6), ETH Z\"urich and
Ruhr-Universit\"at Bochum}
\address[A]{J. Jacod\\
Institut de Math\'ematiques de Jussieu\\
UPMC (Universit\'e Paris-6)\\
175 rue du Chevaleret\\
75 013 Paris\\
France\\
\printead{e1}}
\address[B]{M. Podolskij\\
ETH Z\"urich\\
Department of Mathematics\\
8092 Z\"urich\\
Switzerland\\
\printead{e2}}
\address[C]{M. Vetter\\
Ruhr-Universit\"at Bochum\\
Fakult\"at f\"ur Mathematik\\
Universit\"atsstr. 150\\
44780 Bochum\\
Germany\\
\printead{e3}}
\end{aug}

\thankstext{t1}{Supported by Deutsche Forschungsgemeinschaft through
SFB 823.}

\received{\smonth{12} \syear{2008}}
\revised{\smonth{6} \syear{2009}}

%
\begin{abstract}
This paper presents some limit theorems for certain functionals of
moving averages of semimartingales plus noise which are observed at
high frequency. Our method generalizes the pre-averaging approach (see
[\textit{Bernoulli} \textbf{15} (2009) 634--658, \textit{Stochastic
Process. Appl.} \textbf{119} (2009) 2249--2276]) and provides consistent
estimates for various characteristics of general semimartingales.
Furthermore, we prove the associated multidimensional (stable) central
limit theorems. As expected, we find central limit theorems with a
convergence rate $n^{-1/4}$, if $n$ is the number of observations.
\end{abstract}

%

\end{frontmatter}

\section{Introduction} \label{Real-Intro}

The last years have witnessed a considerable development of the
statistics of processes observed at very high frequency due to the
recent availability of such data. This is particularly the case for
market prices of stocks, currencies and other financial
instruments. Correlatively, the technology for the analysis of such
data has grown rapidly. The emblematic problem is the question of
how to estimate daily volatility for financial prices (in stochastic
process terms, the quadratic variation of log prices).

However, those high-frequency data are almost always corrupted by
some noise. This may be recording or measurement errors, a situation
which can be modeled by an additive white noise. For financial data
we also have a different sort of ``noise'' due to the fact that
prices are recorded as multiples of the basic currency unit so that
some rounding is necessarily performed, and the level of rounding is
far from being negligible for very high frequency data in comparison
to the intrinsic variability of the underlying process. For these
reasons, it is commonly acknowledged that the underlying process of
interest, such as the price semimartingale, is latent rather than
observed.

A large amount of work has already been devoted to the subject,
especially for additive white noise, but also for some other types of
noise like rounding effects. A comprehensive discussion of the noise
models and the effect of noise on the inference for the underlying
process may be found in \cite{LM}. Various statistical procedures for
getting rid of the noise have been proposed (see, e.g.,
\cite{AMZ1,BR,BS1,ZMA,Z} and, more closely related to the present work,
\cite{BHLS,PV,PV2,JLMPV}).

Most of the aforementioned papers are concerned with the estimation of
the integrated volatility, that is, the quadratic variation, for a
continuous semimartingale. Only Podolskij and Vetter \cite{PV,PV2} deal
with estimation of various volatility functionals and robustness to
jumps in the discontinuous semimartingale setting with i.i.d. noise. So
there is a lack of more general results, allowing, for example, one to
estimate other powers of the volatility (like the ``quarticity'') or
the sum of some powers of the jumps, for a general It\^o
semimartingale. These quantities have proved extremely useful for a
number of estimation or testing problems in the context of
high-frequency data, but they have been studied when the process is
observed without noise. Recall that the typical statistical problems in
the noise-free framework are (i) estimation of the quadratic variation
(see \cite {BS2,J2}), (ii) tests for the presence of jumps (see
\cite{AJ1,BS1}), (iii) tests for the presence of the continuous
component (see \cite{CM,AJ3}) or (iv) estimation of the ``activity
index'' of the jump part (see \cite{AJ2,TT}).

The aim of this paper is to provide probabilistic tools to solve (some
of) the aforementioned statistical problems in the presence of noise.
Thus this is a rather probabilistic paper, but the interest and
motivation of the forthcoming results lie essentially in potential
applications; therefore, after the main results we give hints toward
how to apply the results for concrete statistical questions, but not a
full account of these applications (see, e.g., Remarks \ref {R111a} and
\ref{R222} or Theorem \ref{TVQ1}).

Let us be more specific. We consider an It\^o semimartingale $X$
which is corrupted by noise. The observed process $Z=(Z_t)_{t\geq0}$
is given as
\[
Z_t = X_t + \chi_t,\qquad t\geq0,
\]
where $(\chi_t)_{t\geq0}$ are errors which are, conditionally on the
process $X$, centered and independent. The process $Z$ is assumed to be
observed at equidistant time points $i\Delta_n$, $i=0,1, \ldots,
[t/\Delta_n]$, with $\Delta_n\rightarrow0$ as $n\rightarrow\infty$.
This structure of noise allows for an additive white noise but also for
noise involving rounding effects since $\chi_t$ may depend on $X_t$, or
even on the whole past of $X$ before time $t$. It rules out, though,
some other interesting types of noise, like an additive colored noise.
Note, however, that the $\chi_t$ are not necessarily independent (the
independence is only ``conditional on $X$'').

In the no-noise case (i.e., $\chi\equiv0$) an extensive theory has been
developed in various papers which allows for estimating quantities like
${\sum_{s\leq t}}|\Delta X_s|^p$ where $\Delta X_s$ denotes the jump size
of $X$ at time $s$, or $\int_0^t|\sigma_s|^p\,ds$ where $\sigma$ is the
volatility. See, for instance, \cite{BGJPS} or \cite{J2} among others.
Typically, these quantities are estimated by sums of powers of the
successive increments of $X$, that is, they are limits of such sums.
When noise is present, these estimators are inadequate because they
converge toward some characteristics of the noise rather than toward
the characteristics of the process $X$ in which we are interested.
There are currently three main approaches to overcome this difficulty,
mainly for the estimation of the quadratic variation in the continuous
case: the subsampling method \cite{Z}, the realized kernel method
\cite{BHLS} and the pre-averaging method \cite{PV,JLMPV} (see also
\cite{GJ} for a comprehensive theory in the parametric setting). All
these approaches achieve the optimal rate of $\Delta_n^{1/4}$. In this
paper we use the pre-averaging method to derive rather general
estimators.

More precisely, we choose a (smooth enough) weight function $g$ on $[0,1]$
and an appropriate sequence $k_n$
with which we associate the (observed) variables,
\begin{eqnarray*}
\overline{Z}(g)^n_i &=& \sum_{j=1}^{k_n-1}g(j/k_n)\bigl(Z_{(i+j)\Delta_n} -
Z_{(i+j-1)\Delta_n}\bigr),
\\
\widehat{Z}(g)^n_i &=& \sum_{j=1}^{k_n}\bigl(g(j/k_n) -
g\bigl((j-1)/k_n\bigr)\bigr)^2 \bigl(Z_{(i+j)\Delta_n} - Z_{(i+j-1)\Delta_n}\bigr)^2.
\end{eqnarray*}
Our aim is to study the asymptotic behavior of the following
functionals:
\[
V(Z,g,p,r)^n_t=\sum_{i=0}^{[t/\Delta_n]-k_n}
|\overline{Z}(g)^n_i|^p|\widehat{Z}(g)^n_i|^r
\]
for suitable powers $p,r\geq0$. The local smoothing performed by the
quantity $\overline{Z}(g)^n_i$ is somewhat related to the idea
proposed in
\cite{ZMA} for the estimation of a certain conditional variance.
Its role is the reduction of the
influence of the noise process $\chi$ whereas $\widehat{Z}(g)^n_i$ is used
for bias corrections. The asymptotic theory for the functionals
$V(Z,g,p,0)^n_t$ in the absence of jumps is (partially) derived in
\cite{JLMPV} and~\cite{PV2}, but here we extend these
results to the case of general semimartingales.

Quite naturally, the asymptotic behavior of $V(Z,g,p,r)^n_t$ is
different according to whether the process $X$ is continuous or
not. In particular, different scaling is required to obtain
nontrivial limits for $V(Z,g,p,r)^n_t$. More precisely, we show
the following ($\stackrel{\mathbb{P}}{\longrightarrow}$ means
convergence in probability, and
$\stackrel{\mathrm{u.c.p.}}{\longrightarrow}$ means convergence in
probability uniformly over all
finite time intervals):
\begin{longlist}
\item For all\vspace*{-2pt} semimartingales $X$ it holds that
$\frac1{k_n} V(Z,g,p,0)^n_t\stackrel{\mathbb{P}}{\longrightarrow}
\overline{g}(p)\times\break\sum_{s\leq t}|\Delta X_s|^p$ for
$p>2$\vspace*{1pt} and $\frac1{k_n} V(Z,g,2,0)^n_t-\frac1{2k_n}
V(Z,g,0,1)^n_t\stackrel {\mathbb{P}}{\longrightarrow}
\overline{g}(2)[X,X]_t$ where the $\overline{g}(p)$'s are known
constants (which depend on $g$), and $[X,X]$ is the quadratic variation
of $X$.
\item When $X$ is a continuous It\^o semimartingale it holds that
$\Delta_n^{1-p/4}V(Z,\break g,p$,$0)^n_t\stackrel{\mathrm
{u.c.p.}}{\longrightarrow}m_p\int_0^t |\theta
\overline{g}(2)\sigma_s^2+\frac{\overline{g}{}'(2)}{\theta} \alpha
_s^2 |^{p/2}\,ds$ where
$m_p$, $\theta$ are certain constants, $(\sigma_s^2)$~is the volatility
process and $(\alpha_s^2)$ is the local conditional variance of the
noise process $\chi$. Furthermore, a proper linear combination of
$V(Z,g,p,r)^n_t$ for integers $p,r$ with $p+2r=l$ converges in
probability to $\int_0^t |\sigma_s|^{l} \,ds$ when $l$ is an even integer.
\end{longlist}
For each of the aforementioned cases we prove a joint stable
central limit theorem for a given family of weight functions
$(g_i)_{1\leq i\leq d}$ [for the first functional in (i) we
additionally have to assume that $p>3$]. The corresponding
convergence rate is~$\Delta_n^{1/4}$.

We end this introduction by emphasizing that only the one-dimensional
case for $X$ is studied here. The extension to multi-dimensional
semimartingales is possible, and even mathematically rather straightforward,
but extremely cumbersome.

This paper is organized as follows: in Section \ref{sec-Intro} we
introduce the setting and the assumptions. Sections \ref{sec-LLN}
and \ref{sec-CLT} are devoted to stating the results, first
the various convergences in probability and second the associated
central limit theorems. The proofs are gathered in Section
\ref{sec-PP}.

\section{The setting}\label{sec-Intro}

We have a one-dimensional underlying process $X=(X_t)_{t\geq0}$,
and observation times $i\Delta_n$ for all $i=0,1,\ldots,k,\ldots$
with $\Delta_n\to0$. We suppose that $X$ is a semimartingale which
can thus be written as
%
%
\begin{equation}\label{SET1}
X=X_0+B+X^c+\bigl(x1_{\{|x|\leq1\}}\bigr)\star(\mu-\nu)+\bigl(x1_{\{|x|>1\}}\bigr)\star
\mu.
\end{equation}
Here $\mu$ is the jump measure of $X$ with $\nu$ its predictable
compensator; $X^c$ is the continuous (local) martingale part of $X$,
and $B$ is the drift. All these are defined on some filtered
probability space $(\Omega^{(0)},\mathcal F^{(0)},(\mathcal
F_t^{(0)})_{t\geq0},\mathbb {P}^{(0)})$. We use here the usual notation
of stochastic calculus, and for any unexplained (but standard) notation
we refer to \cite{JS}; for example
$\psi\star(\mu-\nu)_t=\int_0^t\int_{\mathbb{R}}\psi
(s,x)(\mu-\nu)(ds,dx)$ is the stochastic integral of the predictable
function $\psi(\omega,t,x)$ with respect to the martingale measure
$\mu-\nu$, when it exists.

The process $X$ is observed with an error; that is, at stage $n$, and
instead of the values $X^n_i=X_{i\Delta_n}$ for $i\geq0$, we observe
$X^n_i+\chi^n_i$ where the $\chi^n_i$'s are ``errors'' which are,
conditionally on the process $X$, centered and independent (this allows
for errors which are depending on $X$ and thus may be unconditionally
dependent). It is convenient to define the noise $\chi_t$ for any time
$t$, although at stage $n$ only the values $\chi_{i\Delta_n}$ are
really used.

Mathematically speaking, this can be formalized as such: for each
$t\geq0$, we have a transition probability $Q_t(\omega^{(0)},dz)$ from
$(\Omega^{(0)},\mathcal F_t^{(0)})$ into $\mathbb{R}$. We endow the
space $\Omega^{(1)}=\mathbb{R}^{[0,\infty)}$ with the product Borel
$\sigma$-field $\mathcal F^{(1)}$ and the ``canonical process''
$(\chi_t\dvtx t\geq0)$ and with the probability
$\mathbb{Q}(\omega^{(0)},d\omega^{(1)})$ which is the product
$\bigotimes_{t\geq0}Q_t(\omega^{(0)},\cdot)$. We introduce the filtered
probability space $(\Omega,\mathcal F,(\mathcal F_t)_{t\geq0},\mathbb
{P})$ and the filtration $(\mathcal G_t)$ as follows:
%
%
\begin{equation}\label{S1}\quad
\left.
\begin{array}{l}
\Omega=\Omega^{(0)}\times\Omega^{(1)},\qquad
\mathcal F=\mathcal F^{(0)}\otimes\mathcal F^{(1)},\\
\mathcal F_t=\mathcal F^{(0)}_t\otimes\sigma\bigl(\chi_s\dvtx s\in
[0,t)\bigr),\qquad
\mathcal G_t=\mathcal F^{(0)}\otimes\sigma\bigl(\chi_s\dvtx
s\in[0,t)\bigr),\\[2pt]
\mathbb{P}\bigl(d\omega^{(0)},d\omega^{(1)}\bigr)=\mathbb{P}^{(0)}\bigl(d\omega
^{(0)}\bigr)\mathbb{Q}
\bigl(\omega^{(0)},d\omega^{(1)}\bigr).
\end{array}
\right\}
\end{equation}
Any variable or process which is defined on either $\Omega^{(0)}$ or
$\Omega^{(1)}$ can be considered in the usual way as a variable or a
process on $\Omega$. Note that $X$ is still a semimartingale with the
same decomposition (\ref{SET1}) on $(\Omega,\mathcal F,(\mathcal
F_t)_{t\geq0},\mathbb{P})$ despite the fact that
the filtration $(\mathcal F_t)$ is not right-continuous. On the other hand,
the ``process'' $\chi$ typically has no measurable property in
time since under $\mathbb{Q}(\omega^{(0)},\cdot)$ it is constituted of independent
variables; as mentioned before, only the values of $\chi$ at the
observation times are relevant, and the extension as a process
indexed by $\mathbb{R}_+$ is for notational convenience only.

At time $t$, instead of $X_t$, we observe the variable
%
%
\begin{equation}\label{SET10}
Z_t=X_t+\chi_t.
\end{equation}
We make the following crucial assumption on the noise, for some $q\geq2$:
\renewcommand{\theHypothesis}{(\textit{N}-$q$)}
\begin{Hypothesis}\label{hypoNq}
There is a sequence of $(\mathcal F^{(0)}_t)$-stopping times $(T_n)$
increasing to $\infty$, such that $\int Q_t(\omega^{(0)},dz) |z|^q\leq
n$ whenever $t<T_n(\omega^{(0)})$. We write for any integer $r\leq q$,
%
%
\begin{equation}\label{SET12}
\beta(r)_t\bigl(\omega^{(0)}\bigr) = \int Q_t\bigl(\omega^{(0)},dz\bigr)
z^r,\qquad
\alpha_t = \sqrt{\beta(2)_t},
\end{equation}
and we also assume that
%
%
\begin{equation}\label{SET11}
\beta(1) \equiv0.
\end{equation}
\end{Hypothesis}

In most applications, the local boundedness of the $q$th moment
of the noise, even for all $q>0$, is not a serious restriction.
Condition (\ref{SET11}), on the other hand, is a quite
serious restriction (see \cite{JLMPV} for a discussion of the implications
of this assumption, and below are some examples).
\begin{example}\label{E2.0} The structure of the noise allows for
an additive white noise [all $Q_t(\omega^{(0)},\cdot)$ are equal to
a fixed probability measure, independent of $(\omega^{(0)},t)$, with
mean $0$]. It also allows for
some sort of rounding which means that the observed process $Z_t$
takes its values in $\alpha\mathbb{Z}$ where $\alpha>0$ is the
rounding level; for
example if the $\xi_t$ are i.i.d. uniform on
$[0,\alpha]$ and independent of $\mathcal F^{(0)}$ (hence of $X$) and
$Z_t=\alpha[(X_t+\xi_t)/\alpha]$ (here $[x]$ denotes the integer
part of the real~$x$), we have $Q_t(\cdot,dx)=
(X_t/\alpha-[X_t/\alpha])\varepsilon_{\alpha[X_t/\alpha
]+1-X_t}(dx)+(1-X_t/\alpha
+[X_t/\alpha])
\varepsilon_{\alpha[X_t/\alpha]-X_t}(dx)$ which satisfies Hypothesis \ref{hypoNq}
for all $q$
(here $\varepsilon$ denotes the Dirac measure). Many other
specifications of rounding errors are possible, obviously.

However, it unfortunately does not allow for ``pure rounding,'' that is,
$Z_t=\alpha[X_t/\alpha]$; although in this case we have the
structure (\ref{S1}), the property (\ref{SET11}) is violated. In this case,
there is no way of estimating the
integrated volatility in a consistent way because
this quantity is not even a function of the path $t\mapsto Z_t$
in the ``completely observed'' case.
\end{example}

We choose a sequence of integers $k_n$ satisfying for some $\theta>0$,
%
%
\begin{equation}\label{SET3}
k_n\sqrt{\Delta_n} = \theta+\mathrm{o}(\Delta_n^{1/4});\qquad
\mbox{we write }
u_n = k_n\Delta_n.
\end{equation}
We will also consider weight functions $g$ on $[0,1]$, satisfying
%
%
\begin{equation}\label{1}
\left.
\begin{array}{l}
g \mbox{ is continuous, piecewise $C^1$ with a piecewise Lipschitz
derivative $g'$,}\\
g(0)=g(1)=0,\qquad \displaystyle\int_0^1g(s)^2\,ds > 0.
\end{array}
\right\}\hspace*{-28pt}
\end{equation}
It is convenient to extend such a $g$ to the whole of $\mathbb{R}$ by setting
$g(s)=0$ if $s\notin[0,1]$. We associate with $g$ the following
numbers [where $p\in(0,\infty)$ and $i\in\mathbb{Z}$]:
%
%
\begin{equation}\label{3}
\left.
\begin{array}{l}
g^n_i = g(i/k_n), \qquad
g'^n_i = g^n_i-g^n_{i-1},\\[2pt]
\displaystyle \overline{g}(p)_n=\sum_{i=1}^{k_n}|g^n_i|^p, \qquad
\overline{g}{}'(p)_n=\sum_{i=1}^{k_n}|g'^n_i|^p.
\end{array}
\right\}
\end{equation}
If $g,h$ are bounded functions with support in $[0,1]$, and $p>0$ and
$t\in\mathbb{R}$, we set
%
%
\begin{equation}\label{30}
\overline{g}(p) = \int|g(s)|^p \,ds,\qquad
\overline{(gh)}(t) = \int g(s)h(s-t) \,ds.
\end{equation}
For example, $\overline{g}{}'(p)$ is associated with $g'$ by the first definition
above, and $\overline{g}(2)=\overline{(gg)}(0)$. Note that, as $n\to
\infty$,
%
%
\begin{equation}\label{SET4}
\overline{g}(p)_n = k_n \overline{g}(p)+\mathrm{O}(1),\qquad
\overline{g}{}'(p)_n = k_n^{1-p} \overline{g}{}'(p)+\mathrm{O}(k_n^{-p}).
\end{equation}

With any process $Y=(Y_t)_{t\geq0}$ we associate the following
random variables:
%
%
\begin{equation}\label{RR2}
\left.
\begin{array}{l}
Y_i^n = Y_{i\Delta_n},\qquad \Delta^n_iY = Y_{i\Delta
_n}-Y_{(i-1)\Delta_n},\\[2pt]
\displaystyle\overline{Y}(g)^n_i = \sum_{j=1}^{k_n-1}g^n_j\Delta^n_{i+j}Y =
-\sum_{j=1}^{k_n}g'^n_jY^n_{i+j-1},\\
\displaystyle\widehat{Y}(g)^n_i = \sum_{j=1}^{k_n}(g'^n_j\Delta^n_{i+j}Y)^2,
\end{array}
\right\}
\end{equation}
and we define the $\sigma$-fields $\mathcal F^n_i=\mathcal F_{i\Delta
_n}$ and
$\mathcal G^n_i=\mathcal G_{i\Delta_n}$.

Now we can define the processes of interest for this paper. Below,
$p$ and $r$ are nonnegative reals, and typically the process
$Y$ will be $X$ or $Z$.
%
%
\begin{equation}\label{SET5}
V(Y,g,p,r)^n_t=\sum_{i=0}^{[t/\Delta_n]-k_n}|\overline
{Y}(g)^n_i|^p |\widehat{Y}(g)^n_i|^r.
\end{equation}
\begin{rem}\label{R2.0} The process $V(Z,g,p,0)^n_t$ is the
realized $p$-variation of moving averages of the observations
$Z_{i\Delta_n}$
over a window of size
$u_n=k_n\Delta_n$ and is designed to wipe out the influence of the noise.
The influence of the noise after using this procedure is of order
of magnitude $1/\sqrt{k_n}$ because the averaging uses $k_n$ observations.
On the other hand when there is no noise
but we still take moving averages, the rate of convergence
of our functionals are typically $\sqrt{u_n}$ because at time $t$
the summands (the number of which is about $t/\Delta_n$)
are strongly dependent; if we want enough independence
to obtain a CLT we basically have to consider nonoverlapping intervals
whose number is about $t/u_n$.

The ``overall'' rate of convergence is of order $\frac1{\sqrt{k_n}}
\vee\sqrt{u_n}$ ; this explains the choice (\ref{SET3}) for $k_n$
which amounts to optimizing the rate. Of course, doing so does
not completely wipe out the noise which then comes as a bias; this is
why we need the complicated processes $V(Z,g,p,r)^n_t$ in order to
remove this bias (see Remark~\ref{R2.1} below).
\end{rem}

Finally we state our assumptions on $X$. One of these
is that $X$ is an \textit{It\^o semimartingale}.
This means that its characteristics are absolutely continuous with
respect to Lebesgue measure, or equivalently that it can be written as
%
%
\begin{eqnarray}\label{CD300}
X_t &=& X_0+\int_0^tb_s\,ds+\int_0^t\sigma_s\,dW_s
\nonumber\\[-8pt]\\[-8pt]
&&{} +\bigl(\delta1_{\{|\delta|\leq1\}}\bigr)\star(\underline{\mu}-\underline
{\nu})_t+\bigl(\delta1_{\{
|\delta|>1\}}\bigr)\star\underline{\mu}_t,\nonumber
\end{eqnarray}
where $W$ is a Brownian motion and $\underline{\mu}$ and $\underline
{\nu}$ are
a Poisson random measure on $\mathbb{R}_+\times E$, and its compensator
$\underline{\nu}(dt,dz)=dt\otimes\lambda(dz)$ [where $(E,\mathcal
E)$ is an auxiliary space
and $\lambda$ a $\sigma$-finite measure]. The required regularity and
boundedness conditions on the coefficients $b,\sigma,\delta$ are gathered
in the following:
\renewcommand{\theHypothesis}{(H)}
\begin{Hypothesis}\label{hypoH}
The process $X$ has the form (\ref{CD300})
[on $(\Omega^{(0)},\mathcal F^{(0)},\break (\mathcal F_t^{(0)})$,$\mathbb
{P}^{(0)})$], and:

\textup{(a)} the process $(b_t)$ is optional and locally bounded;

\textup{(b)} the processes $(\sigma_t)$ is c\`adl\`ag ($=$ right-continuous
with left limits) and adapted;

\textup{(c)} the function $\delta$ is predictable, and there is a bounded function
$\gamma$ in $\mathbb{L}^2(E,\mathcal E,\lambda)$ such that
the process $\sup_{z\in E}(|\delta(\omega^{(0)},t,z)|\wedge
1)/\gamma(z)$ is
locally bounded.
\end{Hypothesis}

In particular, a \textit{continuous It\^o semimartingale} is of the form
%
%
\begin{equation}\label{CIS}
X_t = X_0+\int_0^tb_s\,ds+\int_0^t\sigma_s\,dW_s,
\end{equation}
where the processes $b$ and $\sigma$ are optional [relative to
$(\mathcal F^{(0)}_t)$] and such that the integrals above make sense.
When this is the case, we sometimes need the process $\sigma$ itself
to be an It\^o semimartingale; it can then be written
as in (\ref{CD300}), but another way of expressing this property is as
follows [we are again on the space
$(\Omega^{(0)},\mathcal F^{(0)},(\mathcal F_t^{(0)}),\mathbb{P}^{(0)})$]:
%
%
\begin{equation}\label{CC2}
\sigma_t=\sigma_0+\int_0^t\widetilde{b}_s\,ds+\int_0^t\widetilde
{\sigma}_s\,dW_s+M_t
+\sum_{s\leq t}\Delta\sigma_s 1_{\{|\Delta\sigma_s|>v\}},
\end{equation}
where $M$ is a local martingale orthogonal to $W$ and with bounded
jumps and $\langle M,M\rangle_t=\int_0^ta_s\,ds$, and the compensator of
$\sum_{s\leq t}1_{\{|\Delta\sigma_s|>v\}}$ is $\int_0^ta'_s\,ds$, and where
$\widetilde{b}_t$, $a_t$, $a'_t$ and $\widetilde{\sigma}_t$ are
optional processes; the first
three being locally integrable and the fourth being locally
square-integrable. Then we set the following:
\renewcommand{\theHypothesis}{(K)}
\begin{Hypothesis}\label{hypoK}
We have (\ref{CIS}) and (\ref{CC2}), and
the processes $\widetilde{b}_t$, $a_t$, $a'_t$ are locally bounded
whereas the processes $b_t$ and $\widetilde{\sigma}_t$ are
left-continuous with right limits.
\end{Hypothesis}
\begin{rem}\label{R2.1}
(i) The intuition behind the quantities $\overline{Z}(g)^n_i$ and
$\widehat{Z}(g)^n_i$
can be
explained as follows. Assume for simplicity that $X$ is the continuous
It\^o
semimartingale (\ref{CIS}) and the noise process
$\chi$ is independent of $X$. Now, conditionally on $\mathcal F^n_i$, it
holds that
\[
\Delta_n^{-1/4} \overline{Z}(g)^n_i \stackrel{\mathrm{asy}}{\sim
} \mathcal N \biggl(0, \theta\overline
{g}(2) \sigma^2_{i\Delta_n} +
\frac{\overline{g}{}'(2)}{\theta} \alpha^2_{i\Delta_n} \biggr),
\]
when the processes
$\alpha$ and $\sigma$ are continuous on the interval $(i\Delta_n, (i
+k_n)\Delta_n]$. On the other hand, we have that
\[
\widehat{Z}(g)^n_i \approx\frac{2\overline{g}{}'(2)}{k_n} \alpha
^2_{i\Delta_n},
\]
when the process $\alpha$ is continuous on the interval $(i\Delta_n,
(i +k_n)\Delta_n]$ (this approximation holds even for all
semimartingales $X$). It is now intuitively clear that a certain
combination of the quantities $\overline{Z}(g)^n_i$ and $\widehat
{Z}(g)^n_i$ can be
used to estimate some functions of $\sigma_{i\Delta_n}$ (which is usually
the main object of interest). In particular, a proper linear
combination of $V(Y,g,p-2l,l)^n_t$, $l=0, \ldots, p/2$, for an
even number $p$, converges in probability to $\int_0^t
|\sigma_s|^p \,ds$. This intuition is formalized in Theorems \ref{TLNN1}
and \ref{TL3N1}.

(ii) In the continuous case the quantities $\overline{Z}(g)^n_i$ and
$\widehat{Z}(g)^n_i$
are asymptotically $k_n$-dependent, that is, $\overline{Z}(g)^n_i$
[resp., $\widehat{Z}(g)^n_i$]
is asymptotically (conditionally) independent of $\overline{Z}(g)^n_j$
[resp., $\widehat{Z}(g)^n_j$]
when $|i-j|>k_n$. Thus we will apply a classical block splitting
technique for $m$-dependent
variables to derive the central limit theorem for $V(Y,g,p,r)^n_t$ when
$X$ is continuous
(see Section \ref{BS}).
\end{rem}

\section{Results: The laws of large numbers}\label{sec-LLN}

\subsection{LLN for all semimartingales}\label{ssec:LLNS}

We consider here an LLN which holds for \textit{all} semimartingales,
and we start with the version without noise, that is, $Z=X$. For the
sake of comparison, we recall the following classical result:
%
%
\begin{equation}\label{LLN50}
\sum_{i=1}^{[t/\Delta_n]}|\Delta^n_iX|^p \stackrel{\mathbb
{P}}{\longrightarrow}
\cases{
\displaystyle\sum_{s\leq t}|\Delta X_s|^p, &\quad if $p>2$,\vspace*{2pt}\cr
[X,X]_t, &\quad if $p=2$.}
\end{equation}
Below, and throughout the paper, $g$ always denotes a weight function
satisfying~(\ref{1}).
\begin{theo}\label{TLLN1} For any $t\geq0$ which is not a fixed
time of discontinuity of $X$, we have
%
%
\begin{equation}\label{LLN1}
\frac1{k_n} V(X,g,p,0)^n_t \stackrel{\mathbb{P}}{\longrightarrow}
\cases{
\displaystyle\overline{g}(p) \sum_{s\leq t}|\Delta X_s|^p, &\quad if $p>2$,\vspace*{2pt}\cr
\overline{g}(2) [X,X]_t,&\quad if $p=2$,}
\end{equation}
as soon as $k_n\to\infty$ and $u_n=k_n\Delta_n\to0$ [that is, we do not
need (\ref{SET3}) here].
\end{theo}

This convergence also holds for any $t$ such that $t/\Delta_n$ is an
integer for all $n$, if this happens, but it \textit{never} holds in the
Skorokhod sense, except of course when $X$ is continuous. Taking
in (\ref{SET5}) test functions of the form $f(x)=|x|^p$ is
essential here: the convergence of $\sum_{i=0}^{[t/\Delta
_n]-k_n}f(\overline{X}(g)^n_i)$ for
more general $f$ is so far an open question.

Next we have the version with noise, again for an arbitrary
semimartingale $X$. In the previous
theorem nothing is said about $V(X,g,p,r)^n_t$ when $r\geq1$
which are of little interest. However, when noise
is present, we need those processes to remove an intrinsic bias, and so we
provide their behavior, or at least some (rough) estimates on them.
\begin{theo}\label{TLLN2} \textup{(a)} For any $t\geq0$ which is not a fixed
time of discontinuity of $X$ we have
%
%
\begin{equation}\label{LLN111}\qquad
p>2 \mbox{ and (N-$p$) holds}\quad \Rightarrow\quad
\frac1{k_n} V(Z,g,p,0)^n_t \stackrel{\mathbb{P}}{\longrightarrow
} \overline{g}(p)\sum_{s\leq t}|\Delta X_s|^p.
\end{equation}
Moreover, if $r>0$ and $p+2r>2$ and if [N-$(p+2r)$] holds, then
%
%
\begin{equation}\label{LLN119}
\mbox{the sequence
$ \bigl(k_n^{r-({p+4r})/({p+2r})} V(Z,g,p,r)^n_t \bigr)$ is tight.}
\end{equation}

\textup{(b)} Under (N-$2$) we have for all $t$ as above,
%
%
\begin{equation}\label{LLN112}
\frac1{k_n} V(Z,g,2,0)^n_t-\frac1{2k_n} V(Z,g,0,1)^n_t \stackrel
{\mathbb{P}}{\longrightarrow}
\overline{g}(2) [X,X]_t.
\end{equation}
\end{theo}

It is worth emphasizing that the behaviors of $V(Z,g,p,0)^n$
and $V(X,g,p,0)^n$ are basically the same when $p>2$, at least
for the convergence in probability because the jumps dominate in these
processes both the ``continuous martingale part'' and the noise, and, in
particular, by using the pre-averaging procedure, we wipe out the
noise completely in this case. On the opposite, when $p=2$ the two processes
$V(Z,g,2,0)^n$ and $V(X,g,2,0)^n$ behave differently, even at the
level of convergence in probability.

\subsection{\texorpdfstring{LLN for continuous It\^o
semimartingales---1}{LLN for continuous Ito
semimartingales---1}}
\label{ssec:LLNCIS1}

When $X$ is continuous, Theorem \ref{TLLN2} gives a vanishing limit
when $p>2$, so it is natural in this case to look for a normalization
which provides a nontrivial limit. Exactly as when there is no noise
(see \cite{J2}) this is possible only when
$X$ is a continuous It\^o semimartingale of the form (\ref{CIS}).
\begin{theo}\label{TLNN1} Assume Hypothesis \ref{hypoNq} for some $q>2$ and that $X$
is given by (\ref{CIS}). Assume also
that $b$ is locally bounded and that $\sigma$ and $\alpha$ are
c\`adl\`ag. Then if $0<p\leq q/2$ we have
%
%
\begin{equation}\label{LNN5}
\Delta_n^{1-p/4} V(Z,g,p,0)^n_t \stackrel{\mathit{
u.c.p.}}{\longrightarrow} m_p\int_0^t \biggl|\theta
\overline{g}(2)\sigma_s^2+\frac{\overline{g}{}'(2)}{\theta} \alpha
_s^2 \biggr|^{p/2} \,ds,
\end{equation}
where $m_p$ denotes the $p$th absolute moment of $\mathcal N(0,1)$.
\end{theo}

This result should be compared to the fact
that, under the same assumptions on $X$,
the processes $\Delta_n^{1-p/2}\sum_{i=1}^{[t/\Delta_n]}|\Delta
^n_iX|^p$ converge to
the limiting process $m_p\int_0^t|\sigma_s|^{p} \,ds$.

This theorem is not really satisfactory; unlike Theorem \ref
{TLLN2}(a), the limit depends on the noise,
through $\alpha_s$, and further, we
do not know how to prove a CLT associated to it because of the
intrinsic bias due to the noise (see Remark \ref{R2.1}). However, at
least when $p$ is an even integer (the most interesting case in practice),
we have a useful substitute.
That is, by an application of the binomial formula and the estimation of
the terms that involve the process $\alpha_s$, we obtain (up to a constant
factor) the process $\int_0^t|\sigma_s|^{p} \,ds$ in the limit.

For any even integer $p\geq2$ we introduce the numbers
$\rho_{p,l}$ for $l=0,\ldots,p/2$ which are the solutions of the
following triangular system of linear equations ($C^p_q=\frac{q!}
{p!(q-p)!}$ denote the binomial coefficients):
%
%
\begin{equation}\label{L3N2}
\left.
\begin{array}{l}
\rho_{p,0} = 1,\\[2pt]
\displaystyle\sum_{l=0}^j2^l m_{2j-2l} C_{p-2l}^{2j-2l} \rho_{p,l} = 0,\qquad
j=1,2,\ldots,p/2.
\end{array}
\right\}
\end{equation}
These could, of course, be explicitly computed, and, for example, we have
%
%
\begin{equation}\label{L3N3}
\rho_{p,1} = -\tfrac12 C_p^2,\qquad
\rho_{p,2} = \tfrac34 C_p^4,\qquad \rho_{p,3} = -\tfrac{15}8 C_p^6.
\end{equation}
Then for any process $Y$ and for $p\geq2$ an even integer we set
%
%
\begin{equation}\label{L3N4}
\overline{V}(Y,g,p)^n_t = \sum_{l=0}^{p/2}\rho_{p,l} V(Y,g,p-2l,l)^n_t.
\end{equation}
\begin{theo}\label{TL3N1} \textup{(a)} Let $X$ be an arbitrary semimartingale,
and assume (N-$p$) for some even integer
$p\geq2$. Then for all $t\geq0$ we have
%
%
\begin{equation}\label{L3N5}
\frac1{k_n} \overline{V}(Z,g,p)^n_t \stackrel{\mathbb
{P}}{\longrightarrow} \cases{
\displaystyle\overline{g}(p)\sum_{s\leq t}|\Delta X_s|^p, &\quad if $p\geq
4$,\vspace*{2pt}\cr
\overline{g}(2) [X,X]_t, &\quad if $p=2$.}
\end{equation}

\textup{(b)} Let $X$ satisfy (\ref{CIS}), and assume (N-$2p$) for some even
integer $p\geq2$. Assume also that $b$ is locally bounded and that
$\sigma$
and $\alpha$ are c\`adl\`ag. Then we have
%
%
\begin{equation}\label{L3N6}
\Delta_n^{1-p/4} \overline{V}(Z,g,p)^n_t \stackrel{\mathit
{u.c.p.}}{\longrightarrow}
m_p (\theta\overline{g}(2))^{p/2}\int_0^t|\sigma_s|^p \,ds.
\end{equation}
\end{theo}

The first part of (\ref{L3N5}) is an obvious consequence of (a) of Theorem
\ref{TLLN2} whereas the second part of (\ref{L3N5}) is nothing other
than (\ref{LLN112}) because $\rho_{2,1}=-1/2$.

\subsection{\texorpdfstring{LLN for continuous It\^o
semimartingales---2}{LLN for continuous Ito semimartingales---2}}
\label{ssec:LLNCIS2}

Statistical applications require ``estimators''
for the conditional variance which will appear in the CLTs
associated with some of the previous LLNs. In other words, we
need to provide some other laws of large numbers, which a
priori seem artificial but are motivated by potential
applications.

To this end we need auxiliary processes to be used
also for the CLTs below. Let $W^1$ and $W^2$ be two independent Brownian
motions on another auxiliary filtered probability space
$(\Omega', \mathcal F',(\mathcal F'_t)_{t\geq0},\mathbb{P}')$. With
any function $g$ satisfying (\ref{1}), and
extended as before on $\mathbb{R}$ by setting it to be $0$ outside
$[0,1]$, we
define the following Wiener integral processes:
%
%
\begin{equation}\label{KL1}
L(g)_t = \int g(s-t)\,dW^1_s,\qquad L'(g)_t = \int g'(s-t)\,dW^2_s.
\end{equation}
If $h$ is another function satisfying (\ref{1}), we define $L(h)$ and
$L'(h)$ likewise, \textit{with the same $W^1$ and $W^2$}. The
four-dimensional process $U:=(L(g),L'(g),L(h)$, $L'(h))$ is continuous in
time, centered, Gaussian and stationary. Clearly $(L(g)$, $L(h))$ is
independent of $(L'(g),L'(h))$, and the variables $U_t$ and $U_{t+s}$
are independent if $s\geq1$.

The process $L(g)$ comes in naturally as the limit of $\overline{W}(g)^n_i$
[that is,
$\overline{X}(g)$ when $X=W$]; indeed, we will see
that $L(g)_t$ is the limit in law of $\frac1{\sqrt{u_n}} \overline
{W}(g)^n_{[\Delta_n/t]}$,
and we need the whole process $L(g)_t$ to account for the dependency of the
variables $\overline{W}(g)^n_i$ when $i$ varies. In the same way,
$k_n\widehat{\chi}(g)^n_
{[\Delta_n/t]}$ converges in law to $2(\alpha_t)^2L'(g)_t$ (see the ``key
Lemma'' \ref{LNN10} below).

Some further notation is needed. We set
%
%
\begin{equation}\label{L3N61}
\left.
\begin{array}{l}
m_p(g;\eta,\zeta) = \mathbb{E}'\bigl(\bigl(\eta L(g)_0+\zeta L'(g)_0\bigr)^p\bigr),\\
m_{p,q}(g,h;\eta,\zeta) \\
\qquad = \displaystyle\int_0^2\mathbb{E}' \bigl(\bigl(\eta
L(g)_1+\zeta L'(g)_1\bigr)^p
\bigl(\eta L(h)_t+\zeta L'(h)_t\bigr)^q \bigr) \,dt.
\end{array}
\right\}
\end{equation}
These could of course be expressed by the mean of expectations with
respect to the joint law of $U$ above and, considered as
functions of $(\eta,\zeta)$, they are $C^\infty$. In particular, since
$L(g)_0$ and $L'(g)_0$ are independent centered Gaussian variables
with respective variances $\overline{g}(2)$ and $\overline{g}{}'(2)$,
when $p$ in an
integer we have
%
%
\begin{equation}\label{L3N611}
m_p(g;\eta,\zeta) = \cases{
\displaystyle\sum_{v=0}^{p/2}C_p^{2v} (\eta^2\overline{g}(2))^v (\zeta
^2\overline{g}{}'(2))^{p/2-v}
m_{2v} m_{p-2v},\vspace*{2pt}\cr
\hspace*{32.8pt}\mbox{if $p$ is even},\cr
0,\qquad \mbox{if $p$ is odd.}}
\end{equation}
Next, recalling (\ref{L3N2}), we set for $p\geq2$ an even integer:
%
%
\begin{equation}\label{L3N612}
\left.
\begin{array}{l}
\displaystyle\mu_p(g;\eta,\zeta) = \sum_{r=0}^{p/2}\rho_{p,r} (2\zeta
^2\overline
{g}'(2))^r
m_{p-2r}(g;\eta,\zeta),\\
\displaystyle\mu_{2p}(g,h;\eta,\zeta)=\sum_{r,r'=0}^{p/2}\rho_{p,r} \rho
_{p,r'}
(2\zeta^2\overline{g}{}'(2))^r(2\zeta^2\overline{h}{}'(2))^{r'}\\
\hspace*{100pt}{}\times
m_{p-2r,p-2r'}(g,h;\eta,\zeta),\\
\displaystyle\overline{\mu}_{2p}(g,h;\eta,\zeta) = \mu_{2p}(g,h;\eta,\zeta)-
2\mu_p(g;\eta,\zeta) \mu_p(h;\eta,\zeta).
\end{array}
\right\}
\end{equation}

The following lemma will be useful in the sequel:
\begin{lem}\label{LL3N11} We have
%
%
\begin{equation}\label{L3N613}
\mu_p(g;\eta,\zeta) = m_p \eta^p \overline{g}(2)^{p/2}.
\end{equation}
Moreover if $g_i$ is a finite family of functions
satisfying (\ref{1}), for any $(\eta,\zeta)$ the matrix with entries
$\overline{\mu}_{2p}(g_i,g_j;\eta,\zeta)$ is symmetric nonnegative.
\end{lem}

Finally, we associate with any process $Y$ and
any even integer $p$ the functionals
%
%
\begin{eqnarray}\label{L3N67}\hspace*{22pt}
&&M(Y,g,h;p)^n_t\nonumber\\
&&\qquad=\sum_{r,r'=0}^{p/2}\rho_{p,r}\rho_{p,r'}
\sum_{i=0}^{[t/\Delta_n]-3k_n}(\widehat{Y}(g)^n_i)^r (\widehat
{Y}(h)^n_i)^{r'}\nonumber\\[-8pt]\\[-8pt]
&&\qquad\quad\hspace*{105.3pt}{}\times
\Biggl(|\overline{Y}(g)^n_{i+k_n}|^{p-2r}
\frac1{k_n}\sum_{j=1}^{2k_n}|\overline{Y}(h)^n_{i+j}|^{p-2r'}\nonumber\\
&&\hspace*{180.4pt}{}
-2|\overline{Y}(g)^n_i|^{p-2r}|\overline{Y}(h)^n_{i+k_n}|^{p-2r'}
\Biggr).\nonumber
\end{eqnarray}
Then our last LLN is as follows:
\begin{theo}\label{TL3N7} Let $X$ satisfy (\ref{CIS}), and let $p\geq2$
be an even integer. Assume (N-$2p$), that $b$ is locally bounded and
that $\sigma$ and $\alpha$ are c\`adl\`ag. Then if $p\leq q/2$ and if
$g$ and $h$
are two functions satisfying (\ref{1}), we have
%
%
\begin{equation}\label{L3N68}
\Delta_n^{1-p/2} M(Z,g,h;p)^n_t \stackrel{\mathit
{u.c.p.}}{\longrightarrow} \theta^{-p/2}
\int_0^t\overline{\mu}_{2p}(g,h;\theta\sigma_s,\alpha_s) \,ds.
\end{equation}
\end{theo}

The reader will observe that the limit in (\ref{L3N68}) is symmetrical
in $g$ and $h$, although $M(Y,g,h;p)^n_t$ is not. The motivation for this
result is that it provides consistent estimators for the conditional variance
to be encountered in the CLT below (see Remark \ref{R111a}). Indeed,
as the summands of the $V(Y,g,p,r)^n_t$ are asymptotically
$k_n$-dependent in the continuous case (see Remark \ref{R2.1}),
the statistic $\Delta_n^{1-p/2}M(Z,g,h;p)^n_t$ is, up to a multiplicative
constant, an empirical analogue of the asymptotic conditional
covariance between $\overline{V}(Z,g,p)^n_t$ and $\overline{V}(Z,h,p)^n_t$.

\eject
\section{Results: The central limit theorems}\label{sec-CLT}

\subsection{\texorpdfstring{CLT for continuous It\^o semimartingales}{CLT for continuous Ito semimartingales}}\label{ssec:CLTCIS}

As mentioned before, we do not know whether a CLT associated with the
convergence (\ref{LNN5}) exists. But there is one associated with
(\ref{L3N6}) when $p$ is an even integer. Below we give a
joint CLT for several weight functions $g$ at the same time.
We use the notation
%
%
\begin{equation}\label{CC1}\qquad
\widetilde{V}(g,p)^n_t = \frac1{\Delta_n^{1/4}} \biggl(\Delta_n^{1-p/4}
\overline{V}(Z,g,p)_t^n-m_p (\theta \overline{g}(2))^{p/2}\int
_0^t|\sigma
_s|^p \,ds \biggr).
\end{equation}

In view of Lemma \ref{LL3N11}, the square-root matrix $\psi$ referred to
below exists, and by a standard selection theorem one can find
a measurable version for it. For the stable convergence in law used below,
we refer, for example, to \cite{JS}.
\begin{theo}\label{TCC1} Assume Hypothesis \ref{hypoK} and (N-$4p$), where $p$ is an even
integer, and also that the processes $\alpha$ and $\beta(3)$ are c\`
adl\`ag.
If\vspace*{1pt} $(g_i)_{1\leq i\leq d}$ is a family of functions satisfying
(\ref{1}), for each $t\geq0$ the variables $(\widetilde
{V}(g_i,p)^n_t)_{1\leq
i\leq d}$ converge stably in law to the $d$-dimensional variable,
%
%
\begin{equation}\label{CC5}
\Biggl(\theta^{1/2-p/2}\sum_{j=1}^d\int_0^t\psi_{ij}(\theta\sigma
_s,\alpha
_s) \,dB^j_s
\Biggr)_{1\leq i\leq d},
\end{equation}
where $B$ is a $d$-dimensional Brownian
motion independent of $\mathcal F$ (and defined on an extension of the space),
and $\psi$ is a measurable $d\times d$ matrix-valued function such
that $(\psi\psi^\star)(\eta,\zeta)$ is the matrix with entries
$\overline{\mu}_{2p}(g_i,g_j;\eta,\zeta)$, as defined by (\ref{L3N612}).
\end{theo}

Up to the multiplicative constant $\theta^{1-p/2}$, the $\mathcal
F$-conditional
covariance of the $j$th and $k$th components of (\ref{CC5})
is exactly the right-hand side of (\ref{L3N68}) for $g=g_j$ and $h=g_k$.
\begin{rem}\label{R111a}
An application of Theorem \ref{TL3N7}
and the
properties of stable convergence give now a
a \textit{feasible} version of Theorem \ref{TCC1}. We obtain, for
example, that
the quantity
\[
\frac{\widetilde{V}(g,p)^n_t}{\sqrt{\theta^{1-p/2} \Delta
_n^{1-p/2} M(Z,g,g;p)^n_t}}
\]
converges stably in law (for any fixed $t$) to a variable $U\sim
\mathcal N(0,1)$
independent of $\mathcal F$. The latter can be used to construct confidence
regions for the quantity $\int_0^t|\sigma_s|^p \,ds$ for even $p$'s.
\end{rem}
\begin{rem}\label{R111} Theorem \ref{TCC1} can be extended
to the convergence along finite families of times,
but we do not know whether a \textit{functional}
convergence holds, although it is quite likely.
\end{rem}

\subsection{\texorpdfstring{CLT for discontinuous It\^o semimartingales}{CLT for discontinuous Ito semimartingales}}\label{ssec:CLTDIS}

Now we turn to the case when $X$ jumps. There is a CLT for Theorem
\ref{TLLN2}, at least when $p=2$ and $p>3$, exactly as in \cite{J2}
for the processes of type (\ref{LLN50}). The CLT for Theorem
\ref{TL3N1}, when $p$ is an even integer, takes the same
form. In this subsection we are interested in the case $p>3$, whereas
the case $p=2$ is dealt with in the next subsection.

In view of statistical applications (see Remark \ref{R222} below), and as
in the previous
subsection, we need to consider a family $(g_i)_{1\leq i\leq d}$
of weight functions. We use the notation
%
%
\begin{equation}\label{CD1}
\widetilde{V}^\star(g,p)^n_t = \frac1{\Delta_n^{1/4}} \biggl(\frac
1{k_n} V(Z,g,p,0)_t^n-
\overline{g}(p)\sum_{s\leq t}|\Delta X_s|^p \biggr)
\end{equation}
and, further, when $p\geq4$ is an even integer,
%
%
\begin{equation}\label{CD100}
\overline{V}{}^\star(g,p)^n_t = \frac1{\Delta_n^{1/4}} \biggl(\frac1{k_n}
\overline{V}(Z,g,p)_t^n-\overline{g}(p)\sum_{s\leq t}|\Delta
X_s|^p \biggr).
\end{equation}

These are the processes whose asymptotic behavior is studied, but
to describe the limit we need some rather cumbersome notation which
involves the $d$ weight functions, $g_j$ satisfying (\ref{1}).
For any real $x$
and any $p>0$ we write $\{x\}^p=|x|^p$ sign$(x)$, and we introduce
four $d\times d$ symmetric matrices $\Psi_{p-}$, $\Psi_{p+}$,
$\overline{\Psi}_{p-}$
and $\overline{\Psi}_{p+}$ with entries:
%
%
\begin{equation}\label{CD4}
\left.
\begin{array}{l}
\displaystyle\Psi^{ij}_{p-}=\int_0^1 \biggl(\int_t^1\{g_i(s)\}^{p-1}g_i(s-t)\,ds
\biggr)\\
\hspace*{46.57pt}{}\times
\displaystyle\biggl(\int_t^1\{g_j(s)\}^{p-1}g_j(s-t)\,ds \biggr) \,dt,\\[6pt]
\displaystyle\Psi^{ij}_{p+}=\int_0^1 \biggl(\int_0^{1-t}\{g_i(s)\}
^{p-1}g_i(s+t)\,ds \biggr)\\
\hspace*{46.57pt}{}\times
\displaystyle\biggl(\int_0^{1-t}\{g_j(s)\}^{p-1}g_j(s+t)\,ds \biggr) \,dt,\\[6pt]
\displaystyle\overline{\Psi}^{ij}_{p-}=\int_0^1 \biggl(\int_t^1\{g_i(s)\}
^{p-1}g'_i(s-t)\,ds \biggr)\\
\hspace*{47.51pt}{}\times
\displaystyle\biggl(\int_t^1\{g_j(s)\}^{p-1}g'_j(s-t)\,ds \biggr) \,dt,\\[6pt]
\displaystyle\overline{\Psi}^{ij}_{p+}=\int_0^1 \biggl(\int_0^{1-t}\{g_i(s)\}
^{p-1}g'_i(s+t)\,ds \biggr)\\
\hspace*{47.51pt}{}\times
\displaystyle\biggl(\int_0^{1-t}\{g_j(s)\}^{p-1}g'_j(s+t)\,ds \biggr) \,dt.
\end{array}
\right\}
\end{equation}
These matrices are semi-definite positive, and we can thus consider
four independent sequences of i.i.d. $d$-dimensional variables
$(U_{m-})_{m\geq1}$, $(U_{m+})_{m\geq1}$, $(\overline
{U}_{m-})_{m\geq1}$ and
$(\overline{U}_{m+})_{m\geq1}$, defined on an extension of the space,
independent of~$\mathcal F$, and such that for each $m$ the
$d$-dimensional variables $U_{m-}$, $U_{m+}$, $\overline{U}_{m-}$ and
$\overline{U}_{m+}$ are centered Gaussian vectors with respective covariances
$\Psi_{p-}$, $\Psi_{p+}$, $\overline{\Psi}_{p-}$ and $\overline
{\Psi}_{p+}$. Note that
these variables also depend on $p$ and on the family $(g_j)$,
although it does not show in the notation.

Now let $(T_m)_{m\geq1}$ be a sequence of stopping times with pairwise
disjoint graphs, such that $\Delta X_t\neq0$ implies that $t=T_m$ for
some $m$. As is well known (see \cite{J2}), the following
$d$-dimensional processes are well defined when $p>3$ and $\alpha$ is
c\`adl\`ag, and are $\mathcal F$-conditional martingales:
%
%
\begin{eqnarray}\label{CD6}\quad
U(p)_t &=& p\sum_{m\geq1} \{\Delta X_{T_m} \}^{p-1} \biggl(\sqrt
{\theta}
\sigma_{T_m-}U_{m-}+\frac{\alpha_{T_m-}}{\sqrt{\theta}} \overline
{U}_{m-}\nonumber\\[-8pt]\\[-8pt]
&&\hspace*{85.7pt}{} +\sqrt{\theta} \sigma_{T_m}U_{m+}+\frac{\alpha_{T_m}}{\sqrt
{\theta}} \overline{U}
_{m+} \biggr)
1_{\{T_m\leq t\}}.\nonumber
\end{eqnarray}
Moreover, although these processes obviously depend on the choice of the
times~$T_m$, their $\mathcal F$-conditional laws do not; so if the stable
convergence in law below holds for a
particular ``version'' of $U(p)_t$, it also holds for
all other versions.
\begin{theo}\label{TCD1} Assume Hypothesis \ref{hypoH} and let $p>3$. Assume also
(N-$2p$) and that the process $\alpha$ is c\`adl\`ag. If
$(g_i)_{1\leq i\leq d}$ is a family of functions satisfying~(\ref{1}),
for each $t\geq0$ the variables $(\widetilde{V}
^*(g_i,p)^n_t)_{1\leq
i\leq d}$ converge stably in law to the $d$-dimensional variable
$U(p)_t$.

The same holds for the sequence $(\overline{V}{}^*(g_i,p)^n_t)_{1\leq
i\leq d}$
if further $p$ is an even integer.
\end{theo}
\begin{rem}\label{R222}
In the spirit of \cite{AJ1}, we can use this
result to test for the presence of jumps in the presence of noise. We choose
two distinct one-dimensional weight functions $g$ and $h$. It follows from
Theorem \ref{TL3N1} that, taking, for example, $p=4$,
\[
\frac{\overline{V}(Z,g,4)^n_t}{\overline{V}(Z,h,4)^n_t} \stackrel
{\mathbb{P}}{\longrightarrow} \cases{
\overline{g}(2)^2/\overline{h}(2)^2, &\quad on the set where $X$ is
continuous on $[0,t]$,\cr
\overline{g}(4)/\overline{h}(4), &\quad on the set where $X$ has jumps
on $[0,t]$.}
\]
We can choose $g$ and $h$ such that the two limits above are different.
Then Theorems \ref{TCC1} and \ref{TCD1} provide central limit theorems
for the statistics\break $\overline{V}(Z,g,4)^n_t/\overline{V}(Z,h,4)^n_t$
in both occurrences,
allowing for feasible testing of the two hypotheses. For instance, when
$X$ is continuous,
we deduce that the sequence
\[
\Delta_n^{-1/4} \biggl(\frac{\overline{V}(Z,g,4)^n_t}{\overline
{V}(Z,h,4)^n_t} - \frac
{\overline{g}(2)^2}{\overline{h}(2)^2} \biggr),
\]
converges stably in law toward a mixed normal random variable with
$\mathcal F
^{(0)}$-conditional variance,
\begin{eqnarray*}
&&\theta^{-3} \biggl(3(\theta\overline{h}(2))^{2}\int_0^t\sigma_s^4 \,ds
\biggr)^{-2} \biggl(1, - \frac{\overline{g}(2)^2}{\overline
{h}(2)^2} \biggr)\\
&&\qquad{}\times
(\overline{\mu}_{8}(g_i,g_j;\eta,\zeta))_{1\leq i\leq2, 1\leq
j\leq2}
\biggl(1, - \frac{\overline{g}(2)^2}{\overline{h}(2)^2}
\biggr)^\star,
\end{eqnarray*}
where the \mbox{($2 \times2$)}-matrix $(\overline{\mu}_{8}(g_i,g_j;\eta
,\zeta))_{1\leq
i\leq2, 1\leq j\leq2}$ is defined by (\ref{L3N612}) and $g_1=g$, $g_2=h$.
Since we are able to consistently estimate the above quantity by virtue
of Theorems \ref{TL3N1} and \ref{TL3N7}, we can immediately
obtain a feasible test for the null hypothesis of no jumps.
\end{rem}

\subsection{CLT for the quadratic variation}\label{ssec:CLTQV}

Finally we give a CLT for the quadratic variation
associated with (\ref{LLN112}) when $p=2$ or, equivalently,
with (\ref{L3N5}) which is exactly the same in this case. In contrast
to the preceding results the function $g$ is kept fixed;
thus we will only show a one-dimensional result.
So the processes of interest are simply
%
%
\begin{equation}\label{VQ1}
\overline{V}{}^n_t = \frac1{\Delta_n^{1/4}} \biggl(\frac
1{k_n} \overline{V}(Z,g,2)_t^n-
\overline{g}(2) [X,X]_t \biggr).
\end{equation}

In order to describe the
limit, we introduce an extension of the space on which are defined a
Brownian motion $B$ and variables $U_{m-},\overline
{U}_{m-},U_{m+},\overline{U}_{m+}$
indexed by $m\geq1$; each of these being independent from the others
and independent of $\mathcal F$, and such that the variables $U_{m-}$, $U_{m+}$,
$\overline{U}_{m-}$, $\overline{U}_{m+}$ are centered Gaussian
variables with respective
variances $\Psi_{2-}^{11}$, $\Psi_{2+}^{11}$, $\overline{\Psi
}_{2-}^{11}$ and
$\overline{\Psi}_{2+}^{11}$, as defined in (\ref{CD4}).

As in the previous section, $(T_m)_{m\geq1}$ is a sequence of stopping
times with pairwise disjoint graphs, such that $\Delta X_t\neq0$ implies
that $t=T_m$ for some $m$. Then we associate with these data the
process $U(2)$ as defined by (\ref{CD6}). The result goes as follows:
\begin{theo}\label{TVQ1} Assume Hypothesis \ref{hypoH}. Assume also
(N-$4$) and that the process $\alpha$ is c\`adl\`ag. Then for each
$t$ the variables $\overline{V}{}^n_t$ converge stably
in law to the variable
%
%
\begin{equation}\label{VQ2}
\overline{U}_t = \theta^{-1/2}\int_0^t\sqrt{\overline{\mu
}_4(g,g;\theta\sigma_s,\alpha
_s)} \,dB_s+U(2)_t,
\end{equation}
where $\overline{\mu}_4(g,g;\eta,\zeta)$ is defined by (\ref
{L3N612}) which
here takes the form
%
%
\begin{eqnarray}\label{VQ3}
\overline{\mu}_4(g,g;\eta,\zeta) &=& 4\int_0^1 \biggl(\eta^2\int
_s^1g(u)g(u-s)\,du\nonumber\\[-8pt]\\[-8pt]
&&\hspace*{26.6pt}{} + \zeta^2\int_s^1g'(u)g'(u-s)\,du \biggr)^2
\,ds.\nonumber
\end{eqnarray}
When further $X$ is continuous, the processes $\overline{V}{}^n$
converge stably
(in the functional sense) to the process (\ref{VQ2}) with $U(2)=0$ in
this case.
\end{theo}

When $X$ is continuous, we exactly recover Theorem
\ref{TCC1} when $d=1$ and $g_1=g$, for $p=2$. Note that we do not need
Hypothesis \ref{hypoK} here because of the special feature of the case $p=2$.
When $X$ has jumps, however, the functional convergence does not
hold.
\begin{example}\label{E4.0}
Notice that the limiting variable
$\overline{U}_t$ is mixed normal with $\mathcal F^{(0)}$-conditional variance,
\begin{eqnarray*}
&&\theta^{-1}\int_0^t\overline{\mu}_4(g,g;\theta\sigma_s,\alpha
_s) \,ds\\
&&\qquad{} + 4\sum_{m\geq
1}|\Delta X_{T_m}|^{2} \biggl(\theta
\sigma_{T_m-}^2 \Psi_{2-}+\frac{\alpha_{T_m-}^2}{\theta} \overline
{\Psi}_{2-}\\
&&\hspace*{106.2pt}{}
+\theta\sigma_{T_m}^2 \Psi_{2+}+\frac{\alpha_{T_m}^2}{\theta
} \overline{\Psi}
_{2+} \biggr)
1_{\{T_m\leq t\}}.
\end{eqnarray*}
For the sake of demonstration let us consider the weight function
$g(x)= \min(x,1-x) 1_{\{0\leq x\leq1\}}$.
In this case we obtain
\[
\Psi_{2+} = \Psi_{2-} = \tfrac{151}{80640},\qquad \overline{\Psi
}_{2+} =
\overline{\Psi}_{2-} = \tfrac{1}{96}
\]
and
\[
\overline{\mu}_4(g,g;\eta,\zeta) = 4 \bigl( \tfrac{151}{80640}
\eta^4 + \tfrac
{1}{48} \eta^2\zeta^2 + \tfrac{1}{6} \zeta^4 \bigr).
\]
\end{example}

\section{The proofs}\label{sec-PP}

It is difficult to describe the scheme of the
proofs in a few words, since they are quite technical. However, we can
state the basic ideas:
\begin{itemize}
\item For the case $p>2$ of Theorem \ref{TLLN1} and Theorems \ref{TLLN2}
and \ref{TCD1}, the ``big'' jumps play the leading role, and so the
results are proved first when all jumps are bigger than some
$\varepsilon>0$
(hence there are finitely many of them); we thus examine what happens
around each jump, and show that the rest is negligible.
\item For the continuous case, we use the approximations
%
%
\begin{equation}\label{AU1}
\overline{Z}(g)^n_i \approx\sigma_{i\Delta_n} \overline
{W}(g)^n_i+\overline{\chi}
(g)^n_i,\qquad
\widehat{Z}(g)^n_i \approx2\overline{g}{}'(2) \alpha^2_{i\Delta_n}/k_n.
\end{equation}
Since the approximating quantities in (\ref{AU1}) are asymptotically
$k_n$-dependent
we apply the block splitting technique to prove Theorem \ref{TCC1}. Precisely,
we split the sum over $i$ in the definition of $\overline
{V}(Z,g,p)^n_t$ into
big blocks of
size $mk_n$ which are separated by small blocks of size $k_n$. The big
blocks become asymptotically conditionally independent, and the
small blocks become negligible as $m\to\infty$. In a second step we prove
a CLT for big blocks, for any fixed $m$.
\item For the quadratic variation (case $p=2$ of Theorem \ref{TLLN1} and
Theorem \ref{TVQ1}) the proof is a sort of mixture of the two approaches.
\end{itemize}

In the whole proof $K$ denotes a constant which may change
from line to line. It may depend on the characteristics
of the process $X$ and the law of the noise $\chi$ on $\theta$ and
the two sequences, $(k_n)_{n\geq1}$ and $(\Delta_n)_{n\geq1}$, but neither
on $n$ itself, nor on the index $i$
of the increments $\Delta^n_iX$ or $\Delta^n_iZ$. If it
depends on an additional parameter~$q$; we write it $K_q$.

For the proof of all the results we can use a localization procedure,
described in detail in \cite{J2}, for instance, and which allows us to
systematically replace the Hypotheses \ref{hypoNq}, \ref{hypoH} or \ref{hypoK}, according
to the case, by the following strengthened versions:
\renewcommand{\theHypothesis}{(SN-$q$)}
\begin{Hypothesis}\label{hypoSNq}
We have Hypothesis \ref{hypoNq}, and further
$\int Q_t(\omega^{(0)},\break dz) |z|^q\leq K$.
\end{Hypothesis}
\renewcommand{\theHypothesis}{(SH)}
\begin{Hypothesis}\label{hypoSH}
We have Hypothesis \ref{hypoH}, and the processes $b_t$,
$\sigma_t$, $\sup_{z\in E}|\delta(t,z)|/\gamma(z)$ and $X$
are bounded.
\end{Hypothesis}
\renewcommand{\theHypothesis}{(SK)}
\begin{Hypothesis}\label{hypoSK}
We have Hypothesis \ref{hypoK}, and the processes $b_t$,
$\sigma_t$, $\widetilde{b}_t$, $a_t$, $a'_t$, $\widetilde{\sigma
}_t$ and $X$ are bounded.
\end{Hypothesis}

Observe that under Hypothesis \ref{hypoSK}, and upon taking $v$ large enough
in (\ref{CC2}) (changing $v$ changes the coefficients
$\widetilde{b}_t$ and $a_t$ without altering their boundedness), we can
also suppose that the last term in (\ref{CC2}) vanishes
identically; that is,
%
%
\begin{equation}\label{CC002}
\sigma_t=\sigma_0+\int_0^t\widetilde{b}_s\,ds+\int_0^t\widetilde
{\sigma}_s\,dW_s+M_t.
\end{equation}

Recall that $|g'^n_j|\leq K/k_n$. Then the fact
that conditionally on $\mathcal F^{(0)}$ the $\chi_t$'s
are independent and centered, plus H\"older's inequality, gives us
that under Hypothesis \ref{hypoSNq} we have [the $\sigma$-fields $\mathcal F^n_i$ and
$\mathcal G^n_i$
have been defined after (\ref{RR2})]
%
%
\begin{equation}\label{LLN141}
\left.
\begin{array}{l}
p\leq q \quad\Rightarrow\quad
\mathbb{E}(|\overline{\chi}(g)^n_i|^p\mid\mathcal G^n_i) \leq
K_pk_n^{-p/2},\\
2r\leq q \quad\Rightarrow\quad
\mathbb{E}(|\widehat{\chi}(g)^n_i|^r\mid\mathcal G^n_i) \leq K_rk_n^{-r}.
\end{array}
\right\}
\end{equation}
We will also often use the following property, valid for all
semimartingales $Y$:
%
%
\begin{eqnarray}\label{LLN14}
\overline{Y}(g)^n_i=\int_{i\Delta_n}^{i\Delta_n+u_n}g_n(s-i\Delta
_n)\,dY_s\nonumber\\[-8pt]\\[-8pt]
\eqntext{\mbox{where }
g_n(s)=\displaystyle\sum_{j=1}^{k_n-1}g^n_j
1_{((j-1)\Delta_n,j\Delta_n]}(s).}
\end{eqnarray}

\subsection[Proof of Theorem 3.1]{Proof of Theorem \protect\ref{TLLN1}}\label{ssec:T1}
We start with an arbitrary semimartingale $X$, written as
(\ref{SET1}). We more or less
follow the scheme of the proof of Theorem 2.2 of \cite{J2},
and we use the simplifying notation
$V(Y,p)^n=V(Y,g,p,0)^n$ and $\overline{Y}{}^n_i=\overline{Y}(g)^n_i$. The
basic idea follows: for $\varepsilon\in(0,1]$, we set
%
%
\begin{equation}\label{LLN2}
\left.
\begin{array}{l}
X(\varepsilon) = \bigl(x1_{\{|x|>\varepsilon\}}\bigr)\star\mu,\qquad
M(\varepsilon) = \bigl(x1_{\{|x|\leq\varepsilon\}}\bigr)\star(\mu-\nu),\\
A(\varepsilon) = \langle M(\varepsilon),M(\varepsilon)\rangle
,\qquad
B(\varepsilon) = B-\bigl(x1_{\{\varepsilon<|x|\leq1\}}\bigr)\star\nu,\\
A'(\varepsilon) = \bigl(x^21_{\{|x|\leq\varepsilon\}}\bigr)\star\nu,\qquad
B'(\varepsilon) = \mbox{variation process of } B(\varepsilon),
\end{array}
\right\}
\end{equation}
so that we have
%
%
\begin{equation}\label{LLN3}
X = X_0+B(\varepsilon)+X^c+M(\varepsilon)+X(\varepsilon).
\end{equation}
Then we basically show that $\frac1{k_n} V(B(\varepsilon),p)^n$ and
$\frac1{k_n} V(M(\varepsilon),p)^n$ are ``negligible'' when
$n\to\infty$ and $\varepsilon\to0$, as well as $\frac1{k_n} V(X^c,p)^n$
when $p>2$ whereas $\frac1{k_n} V(X(\varepsilon),p)^n$ converges to
$\overline{g}(p)\sum_{s\leq t}|\Delta X_s|^p1_{\{|\Delta
X_s|>\varepsilon\}}$ and
$\frac1{k_n} V(X^c,2)^n$ converges to $\overline{g}(2)C$ where
$C=\langle X^c,X^c\rangle$.

\textit{Step} 1. Let $B'$ be the variation process of $B$. The process
$B'+C+(x^2\wedge1)\star\nu$ is
predictable, increasing finite-valued and hence locally bounded. By
an obvious localization procedure it is enough to prove the result
under the assumption that, for some constant $K$,
%
%
\begin{equation}\label{LLN4}
B'_\infty+C_\infty+(x^2\wedge1)\star\nu_\infty\leq K.
\end{equation}

We also denote by $T_n(\varepsilon)$ the successive jump times of
$X(\varepsilon)$
with the convention $T_0(\varepsilon)=0$ (which of course is not a jump
time). If $0<\varepsilon<\eta\leq1$, we have
%
%
\begin{equation}\label{LLN8}
\left.
\begin{array}{l}
A(\varepsilon)\leq A'(\varepsilon),\qquad
\Delta B'(\varepsilon)\leq\varepsilon,\qquad |\Delta M(\varepsilon
)|\leq2\varepsilon,\\
\displaystyle B'(\varepsilon)\leq B'+\frac1{\varepsilon} A'(\eta)+\frac1{\eta
} (x^2\wedge
1)\star\nu.
\end{array}
\right\}
\end{equation}
We set $\theta(Y,u,t)=\sup_{s\leq r\leq s+u,r\leq
t}|Y_r-Y_s|$. Observe that
$\overline{Y}{}^n_i=-\sum_{j=1}^{k_n}(g((j+1)/k_n)-g(j/k_n))
(Y_{(i+j)\Delta_n}-Y_{i\Delta_n})$. Hence, since the derivative $g'$ is
bounded, we obtain
%
%
\begin{equation}\label{LLN10} i\leq
[t/\Delta_n]-k_n+1 \quad\Rightarrow\quad|\overline{Y}{}^n_i| \leq
K\theta(Y,u_n,t).
\end{equation}

\textit{Step} 2. Here we study $B(\varepsilon)$. (\ref{LLN10}) and
$\theta(B(\varepsilon),u,t)\leq\theta(B'(\varepsilon),u,t)$ yield
for $p>1$
\[
V(B(\varepsilon),p)^n_t
\leq Kk_n B'(\varepsilon)_t \theta(B'(\varepsilon),u_n,t)^{p-1}.
\]
Since $\Delta B'(\varepsilon)\leq\varepsilon$ we have $\limsup
_{n\to\infty}
\theta(B'(\varepsilon),u_n,t)\leq\varepsilon$, so by (\ref{LLN4})
and (\ref{LLN8})
we have
$\limsup_n\frac1{k_n} V(B'(\varepsilon),p)^n_t\leq
K\varepsilon^{p-1} (\frac1{\eta}+\frac1{\varepsilon} A'(\eta
)_t )$ for all
$0<\varepsilon<\eta\leq1$. Since $A'(\eta)_t\to0$ as $\eta\to0$, we
deduce (choose first $\eta$ small, then $\varepsilon$ smaller) that for
$p\geq2$,
%
%
\begin{equation}\label{LLN9}
\lim_{\varepsilon\to0} \limsup_n \frac1{k_n} V(B(\varepsilon),p)^n_t = 0.
\end{equation}

\textit{Step} 3. In this step, we consider a square-integrable martingale
$Y$ such that $D=\langle Y,Y\rangle$ is bounded. By (\ref{LLN14}),
\[
\mathbb{E}((\overline{Y}{}^n_i)^2) =
\mathbb{E} \biggl(\int_{i\Delta_n}^{i\Delta_n+u_n}g_n(s-i\Delta_n)^2\,dD_s
\biggr) \leq K\mathbb{E}(D_{i\Delta_n+u_n}-D_{i\Delta_n}).
\]
On the other hand, $\mathbb{E}(\overline{Y}{}^n_i\overline
{Y}^n_{i+j})=0$ whenever $j\geq
k_n$. Therefore,
%
%
\begin{equation}\label{LLN11}
\mathbb{E} ((V(Y,2)^n_t)^2 ) \leq k_n\sum_{i=0}^{[t/\Delta
_n]-k_n}\mathbb{E}((\overline{Y}{}^n_i)^2)
\leq Kk_n^2\mathbb{E}(D_t).
\end{equation}

We first apply this with $Y=M(\varepsilon)$, hence $D=A(\varepsilon
)$. In view of
(\ref{LLN11}) and since $A'(\varepsilon)_t\to0$ as $\varepsilon\to
0$ and $A'(\varepsilon
)_t\leq
K$, we deduce
\[
\lim_{\varepsilon\to0} \sup_n \mathbb{E} \biggl( \biggl(\frac1{k_n}
V(M(\varepsilon
),2)^n_t \biggr)^2 \biggr)
= 0.
\]
Since by (\ref{LLN10}) we have $V(M(\varepsilon),p)^n_t\leq K
V(M(\varepsilon),2)^n_t \theta(M(\varepsilon),u_n,t)^{p-2}$ when
$p>2$, and since
$\limsup_n\theta(M(\varepsilon),u_n,t)\leq2\varepsilon$, we get
for $p\geq2$,
%
%
\begin{equation}\label{LLN12}
p\geq2, \eta>0 \quad\Rightarrow\quad
\lim_{\varepsilon\to0} \limsup_{n\to\infty}
\mathbb{P} \biggl(\frac1{k_n} V(M(\varepsilon),p)^n_t>\eta\biggr) = 0.
\end{equation}

Next, (\ref{LLN11}) with $Y=X^c$ yields that the sequence
$\frac1{k_n}V(X^c,2)^n_t$ is bounded in $\mathbb{L}^2$. Using exactly the
same argument as above, where now $\theta(X^c,u_n,t)\to0$, yields
%
%
\begin{equation}\label{LLN13}
p>2, \eta>0 \quad\Rightarrow\quad
\lim_{\varepsilon\to0} \limsup_{n\to\infty}
\mathbb{P} \biggl(\frac1{k_n} V(X^c,p)^n_t>\eta\biggr) = 0.
\end{equation}

\textit{Step} 4. In this step we study $V(X(\varepsilon),p)^n_t$.
We fix $t>0$
such that $\mathbb{P}(\Delta X_t\neq0)=0$. For any $m\geq1$ we set
\[
I(m,n,\varepsilon) = \inf(i\dvtx i\Delta_n\geq T_m(\varepsilon)).
\]
Let $\Omega_n(t,\varepsilon)$ be the set on which two successive
jumps of $X(\varepsilon
)$ in
$[0,t]$ are more than $u_n$ apart, and also $[0,u_n)$ and $[t-u_n,t]$ contain
no jump. Then
$u_n\to0$ and $\mathbb{P}(\Delta X_t\neq0)=0$ yield $\Omega
_n(t,\varepsilon)\to\Omega$
a.s. as
$n\to\infty$. On the set $\Omega_n(t,\varepsilon)$ we have for
$i\leq[t/\Delta
_n]-k_n+1$,
\[
\overline{X(\varepsilon)}^n_i=\cases{
g^n_{I(m,n,\varepsilon)-i} \Delta X_{T_m(\varepsilon)},
& if $I(m,n,\varepsilon)-k_n+1\leq i\leq I(m,n,\varepsilon)-1$\cr
& for some $m$,\cr
0,& otherwise.}
\]
Hence
\[
V(X(\varepsilon),p)^n_t = \overline{g}(p)_n\sum_{s\leq t}|\Delta X_s|^p
1_{\{|\Delta X_s|>\varepsilon\}} \qquad\mbox{on } \Omega
_n(t,\varepsilon),
\]
and (\ref{SET4}) yields
%
%
\begin{equation}\label{LLN7}
\frac1{k_n} V(X(\varepsilon),p)^n_t \to\overline{g}(p)\sum
_{s\leq t}|\Delta X_s|^p
1_{\{|\Delta X_s|>\varepsilon\}}.
\end{equation}

\textit{Step} 5. In this step we study $V(X^c,2)^n_t$. Set
$X^c(n,i)_s=\int_{i\Delta_n}^sg_n(r-i\Delta_n)\,dX^c_r$ when
$s>i\Delta_n$. Using
(\ref{LLN14}) and It\^o's formula, we get
$(\overline{X}{}^{c,n}_i)^2=\zeta^n_i+\zeta'^n_i$ where
\[
\zeta^n_i=\int_{i\Delta_n}^{i\Delta_n+u_n}g_n(s-i\Delta
_n)^2\,dC_s,\qquad
\zeta'^n_i=2\int_{i\Delta_n}^{i\Delta_n+u_n}X^c(n,i)_s\,dX^c_s.
\]
On one hand, $\sum_{i=0}^{[t/\Delta_n]-k_n}\zeta^n_i$ is equal to
$\overline{g}(2)_n C_t$ plus a
term smaller in absolute value than $KC_{u_n}$ and another term
smaller than $K(C_t-C_{t-u_n})$. Then, obviously,
%
%
\begin{equation}\label{LLN15}
\frac1{k_n} \sum_{i=0}^{[t/\Delta_n]-k_n}\zeta^n_i \to\overline
{g}(2) C_t.
\end{equation}

On the other hand, we have $\mathbb{E}(\zeta'^n_i\zeta'^n_{i+j})=0$
when $j\geq
k_n$, and
\[
\mathbb{E}((\zeta'^n_i)^2) \leq4\mathbb{E} \Bigl((C_{i\Delta
_n+u_n}-C_{i\Delta_n})
\sup_{s\in[i\Delta_n,i\Delta_n+u_n]}X^c(n,i)_s^2 \Bigr).
\]
By Doob's inequality,
$\mathbb{E} (\sup_{s\in[i\Delta_n,i\Delta
_n+u_n]}X^c(n,i)_s^4 )\leq K
\mathbb{E}((C_{i\Delta_n+u_n}-C_{i\Delta_n})^2)$, hence the
Cauchy--Schwarz inequality
yields
\[
\mathbb{E}((\zeta'^n_i)^2) \leq K\mathbb{E}\bigl((C_{i\Delta
_n+u_n}-C_{i\Delta_n})^2\bigr)
\leq K\mathbb{E} \bigl((C_{i\Delta_n+u_n}-C_{i\Delta_n})\theta
(C,u_n,t) \bigr),
\]
whenever $i\leq[t/\Delta_n]-k_n+1$. At this point, the same argument
used in (\ref{LLN11}) gives
\[
\mathbb{E} \Biggl( \Biggl(\sum_{i=0}^{[t/\Delta_n]-k_n}\zeta'^n_i
\Biggr)^2 \Biggr) \leq Kk_n^2\mathbb{E}(C_t\theta(C,u_n,t))
\leq Kk_n^2\mathbb{E}(\theta(C,u_n,t)).
\]
But $\theta(C,u_n,t)$ tends to $0$ and is smaller uniformly in $n$ than
a square-integrable variable. We then deduce that $\frac1{k_n}
\sum_{i=0}^{[t/\Delta_n]-k_n}\zeta'^n_i\stackrel{\mathbb
{P}}{\longrightarrow}0$ which, combined with (\ref{LLN15}), yields
%
%
\begin{equation}\label{LLN16}
\frac1{k_n} V(X^c,2)^n_t \stackrel{\mathbb{P}}{\longrightarrow
} \overline{g}(2) C_t.
\end{equation}

\textit{Step} 6. It remains to put all the previous partial results
together. For this we use the following obvious property: for any
$p\geq2$ and $\eta>0$ there is a constant $K_{p,\eta}$ such that
%
%
\begin{equation}\label{LLN17}
x,y\in\mathbb{R} \quad\Rightarrow\quad\bigl||x+y|^p-|x|^p \bigr|
\leq K_{p,\eta}|y|^p+\eta|x|^p.
\end{equation}

Suppose first that $p>2$. Applying (\ref{LLN17}) and
(\ref{LLN3}), we get
\begin{eqnarray*}
&&|V(X,p)^n_t-V(X(\varepsilon),p)^n_t | \\
&&\qquad\leq\eta
V(X(\varepsilon),p)^n_t
+ K_{p,\eta} \bigl(V(B(\varepsilon),p)^n_t+
V(X^c,p)^n_t+V(M(\varepsilon),p)^n_t \bigr).
\end{eqnarray*}
Then by (\ref{LLN9}), (\ref{LLN12}), (\ref{LLN13}) and (\ref{LLN7}),
plus $\sum_{s\leq t}|\Delta X_s|^p
1_{\{|\Delta X_s|>\varepsilon\}}\to\break\sum_{s\leq t}|\Delta X_s|^p$ as
$\varepsilon\to0$, and
by taking $\eta$ arbitrarily small in the above,
we obtain the first part of (\ref{LLN1}).

Next suppose that $p=2$. The same argument shows that
it is enough to prove that
%
%
\begin{equation}\label{LLN18}
\frac1{k_n} V\bigl(X^c+X(\varepsilon),2\bigr)^n_t \stackrel{\mathbb
{P}}{\longrightarrow}
\overline{g}(2) \biggl(C_t+\sum_{s\leq t}|\Delta X_s|^21_{\{|\Delta
X_s|>\varepsilon
\}} \biggr).
\end{equation}
On the set $\Omega_n(t,\varepsilon)$, one easily sees that
\[
V\bigl(X^c+X(\varepsilon),2\bigr)^n_t = V(X^c,2)^n_t+V(X(\varepsilon),2)^n_t
+\sum_{m\geq1\dvtx T_m(\varepsilon)\leq t}\zeta^n_m,
\]
where
\begin{eqnarray*}
\zeta^n_m &=& \sum_{i=I(m,n,\varepsilon
)-k_n+1}^{I(m,n,\varepsilon)-1}
\zeta(m,n,i),\\
\zeta(m,n,i) &=& \bigl|g^n_{I(m,n,\varepsilon)-i}\Delta
X_{T_m(\varepsilon)}
+\overline{X}{}^{c,n}_i \bigr|^2- \bigl|g^n_{I(m,n,\varepsilon)-i}\Delta
X_{T_m(\varepsilon)}
\bigr|^2-|\overline{X}{}^{c,n}_i|^2.
\end{eqnarray*}
In view of (\ref{LLN10}), we deduce from (\ref{LLN17}) that for all
$\eta>0$,
\[
|\zeta(m,n,i)| \leq K_{\eta}\theta(X^c,u_n,t)^2+K\eta\bigl|\Delta
X_{T_m(\varepsilon)}\bigr|^2,
\]
if $I(m,n,\varepsilon)-k_n<i<I(m,n,\varepsilon)$ and $T_m(\varepsilon
)\leq t$.
Since $\eta$ is arbitrarily small, $\zeta^n_m/k_n\to0$ for all
$m$ with $T_m(\varepsilon)\leq t$. Hence (\ref{LLN18}) follows from
(\ref{LLN16}) and (\ref{LLN7}), and we are finished.

\subsection[Proof of Theorem 3.2]{Proof of Theorem \protect\ref{TLLN2}}\label{ssec:T2}

Here $X$ is still an
arbitrary semimartingale, and as for the previous theorem we can
assume by localization that (\ref{LLN4}) holds.
We first prove (a), and we assume Hypothesis \ref{hypoSNq} with $q=p$ for proving
(\ref{LLN111}) and $q=p+2r$ for proving (\ref{LLN119}), so
(\ref{LLN141}) implies
%
%
\begin{equation}\label{LLN54}\quad
\mathbb{E}(V(\chi,g,q,0)^n_t)+\mathbb{E}\bigl(V(\chi,g,0,q/2)^n_t\bigr)
\leq\frac{Kt}{\Delta_nk_n^{q/2}} \leq Ktk_n^{2-q/2}.
\end{equation}
We deduce from (\ref{LLN17}) that, for all $\eta>0$,
\[
|V(Z,g,q,0)^n_t-V(X,g,q,0)^n_t | \leq\eta V(X,g,q,0)^n_t
+K_{q,\eta}V(\chi,g,q,0)^n_t
\]
and thus (\ref{LLN111}) follows from (\ref{LLN1}) and (\ref{LLN54}).

Next, H\"older's inequality yields, when $p,r>0$ with $p+2r=q>2$,
\[
V(Z,g,p,r)^n_t \leq(V(Z,g,q,0)^n_t)^{p/q}
\bigl(V(Z,g,0,q/2)_t^n\bigr)^{2r/q}.
\]
By (\ref{LLN111}), applied with $q$ instead of $p$, we see that the
sequence $k_n^{-1}V(Z,g,q,0)^n_t$ is tight, so for (\ref{LLN119}) it is
enough to show that the sequence $k_n^{q/2-2}V(Z,g,0,q/2)^n_t$ is
also tight.
To see this we first deduce from $|g'^n_j|\leq K/k_n$ that
%
%
\begin{equation}\label{LLN549}
\widehat{X}(g)^n_i \leq \frac K{k_n^2} \sum_{j=i}^{i+k_n-1}(\Delta^n_jX)^2,
\end{equation}
implying\vspace*{-2pt} by H\"older's inequality (recall $q>2$) that
$(\widehat{X}(g)^n_i)^{q/2}\leq\frac
K{k_n^{1+q/2}}\times\break{\sum_{j=i}^{i+k_n-1}}
|\Delta^n_jX|^q$, and hence by (\ref{LLN50}) the sequence
$k_n^{q/2}V(X,g,0,q/2)_t^n$ is tight. Second, (\ref{LLN54})
yields that the sequence $k_n^{q/2-2}V(\chi,g,0,q/2)^n_t$ is tight,
and (\ref{LLN119}) follows because $V(Z,g,0,q/2)^n_t\leq
K_q(V(X,g,0,q/2)^n_t+V(\chi,g,0,q/\break 2)^n_t)$.

Now we turn to (b), so we assume (SN-$2$).
The left-hand side of (\ref{LLN112}) can be written as
\[
\frac1{k_n} V(X,g,2,0)^n_t+\frac1{k_n}\sum_{l=1}^4U(l)^n_t,
\]
where
\[
U(l)^n_t = \cases{
\displaystyle 2\sum_{i=0}^{[t/\Delta_n]-k_n}\overline{X}(g)^n_i \overline{\chi
}(g)^n_i, &\quad if $l=1$,\cr
\displaystyle -\sum_{i=0}^{[t/\Delta_n]-k_n} \sum_{j=1}^{k_n}(g'^n_j)^2\Delta
^n_{i+j}X \Delta^n_{i+j}\chi,
&\quad if $l=2$,\cr
\displaystyle -\frac12 V(X,g,0,1)^n_t, &\quad if $l=3$,\cr
\displaystyle V(\chi,g,2,0)^n_t-\frac12 V(\chi,g,0,1)^n_t,
&\quad if $l=4$,}
\]
and by (\ref{LLN1}) it is enough to prove that for $l=1,2,3,4$,
%
%
\begin{equation}\label{LLN56}
\frac1{k_n} U(l)^n_t \stackrel{\mathbb{P}}{\longrightarrow} 0.
\end{equation}
Equation (\ref{LLN549}) yields $|U(3)^n_t|\leq\frac K{k_n} \sum
_{i=1}^{[t/\Delta_n]}(\Delta^n_i
X)^2$, so (\ref{LLN56}) for $l=3$
follows from (\ref{LLN50}). Next, (\ref{SET11}) implies
$\mathbb{E}(U(l)^n_t\mid\mathcal F^{(0)})=0$ for $l=1,2$; hence (\ref
{LLN56}), for
$l=1,2$ will be implied by
%
%
\begin{equation}\label{LNN569}
\mathbb{E} \biggl( \biggl(\frac1{k_n} U(l)^n_t \biggr)^2\biggm|\mathcal
F^{(0)} \biggr) \stackrel{\mathbb{P}}{\longrightarrow} 0.
\end{equation}
By (\ref{SET11}) and (\ref{RR2}) and (\ref{LLN141}), the variables
$|\mathbb{E}(\overline{\chi}(g)^n_i\overline{\chi}(g)^n_j\mid
\mathcal F^{(0)})|$ vanish if $j\geq k_n$
and are smaller than $K/k_n$ otherwise, whereas the variables
$|\mathbb{E}(\Delta^n_i\chi\Delta^n_{i+j}\chi\mid\mathcal
F^{(0)})|$ are bounded, and
vanish if
$j\geq2$. Then we get
\begin{eqnarray*}
\mathbb{E}\bigl((U(1)^n_t)^2\mid\mathcal F^{(0)}\bigr)&\leq&
\frac K{k_n}\sum_{i=0}^{[t/\Delta_n]-k_n} \sum_{j=1}^{k_n}\overline
{X}(g)^n_i\overline{X}(g)^n_{i+j} \leq K
V(X,g,2,0)^n_t,\\
\mathbb{E}\bigl((U(2)^n_t)^2\mid\mathcal F^{(0)}\bigr)&\leq&
\frac K{k_n^4}\sum_{i,i''=0}^{[t/\Delta_n]-k_n} \sum_{j,j'=0}^{k_n-1}
|\Delta^n_{i+j}X\Delta^n_{i'+j'}X|1_{\{|i'+j'-i-j|\leq2\}}\\
&\leq& \frac K{k_n^2} \sum_{i=1}^{[t/\Delta_n]}(\Delta^n_iX)^2
\end{eqnarray*}
and (\ref{LNN569}) follows from (\ref{LLN1}) when $l=1$ and from
(\ref{LLN50}) when $l=2$.

Finally, an easy calculation shows that $U(4)^n_t=U(5)^n_t+U(6)^n_t$
where
\begin{eqnarray*}
U(5)^n_t &=& \sum_{i=0}^{[t/\Delta_n]}\chi^n_i\sum_{j=1}^{k_n}
\alpha^n_{ij}\chi^n_{i+j},\\
U(6)^n_t &=& \sum_{i=0}^{k_n} \bigl(\alpha'^n_i(\chi^n_i)^2
+\alpha''^n_i\bigl(\chi^n_{i+[t/\Delta_n]-k_n}\bigr)^2 \bigr)
\end{eqnarray*}
for some coefficients $\alpha^n_{ij},\alpha'^n_i,\alpha''^n_i$, all smaller
than $K/k_n$. Then, obviously, $\mathbb{E}(|U(6)^n_t|)\leq
K$ and $E(U(5)^n_t)=0$, and, since $\mathbb{E}(\chi^n_i\chi
^n_{i+j}\chi^n_{i'}
\chi^n_{i'+j'})$ vanishes unless $i=i'$ and $j=j'$ when $j,j'\geq1$,
we also have $\mathbb{E}((U(5)^n_t)^2)\leq Kt/k_n\Delta_n\leq Ktk_n$. Then
(\ref{LLN56}) and (\ref{LNN569}) hold for $l=6$ and $l=5$,
respectively, and thus (\ref{LLN56}) finally holds for $l=4$.

\subsection{A key lemma}\label{ssec-KL}

In this section we prove a key result, useful for deriving the other LLNs
when the process $X$ is continuous and for all CLTs. Before that,
we prove Lemma \ref{LL3N11}.
\begin{pf*}{Proof of Lemma \protect\ref{LL3N11}}
By virtue of (\ref{L3N611}) we have
\[
\mu_p(g;\eta,\zeta) = \sum_{v=0}^{p/2}m_{2v}(\eta^2\overline
{g}(2))^v
(\zeta^2\overline{g}{}'(2))^{p/2-v}\sum_{r=0}^{p/2-v}C_{p-2r}^{2v}
\rho
_{p,r} 2^r
m_{p-2r-2v}.
\]
By (\ref{L3N2}) the last sum above vanishes if $v<p/2$ and equals $1$
when $v=p/2$, hence (\ref{L3N613}). Next, we put $a_i=\mu_p(g_i;\eta
,\zeta)$
and $U^i_t=\eta L(g_i)_t+\zeta L'(g_i)_t$, and, for $T\geq2$,
\[
V^i_T = \sum_{r=0}^{p/2}\rho_{p,r}(2\zeta^2\overline
{g}'_i(2))^r\int_0^T
|U^i_t|^{p-2r}\,dt.
\]
The process $(L(g_i),L'(g_i))$ is stationary, and hence
$\mathbb{E}'(V^i_T)=Ta_i$ for some constant $a_i$. Moreover, the functions
\[
f_{ij}(s,t) = \sum_{r,r'=0}^{p/2}\rho_{p,r} \rho_{p,r'}
(2\zeta^2\overline{g_i}'(2))^r(2\zeta^2\overline{g_j}'(2))^{r'}
\mathbb{E}'(|U^i_s|^{p-2r}
|U^j_{t}|^{p-2r'})-a_ia_j,
\]
satisfy $f_{ij}(s,t)=f_{ij}(s+u,t+u)$ and
$f_{ij}(s,t)=0$ if $|s-t|>1$. Thus if $T>2$,
\begin{eqnarray*}
\operatorname{Cov}(V^i_T,V^j_T)&=&\int_{[0,T]^2}f_{ij}(s,t)\,ds \,dt\\
&=&\int_0^1\,ds\int_0^{s+1}f_{ij}(s,t)\,dt
+\int_{T-1}^T\,ds\int_{s-1}^Tf_{ij}(s,t)\,dt\\
&&{} + \int_1^{T-1}\,ds\int_{s-1}^{s+1}f_{ij}(s,t)\,dt.
\end{eqnarray*}
Therefore $\frac1T\operatorname{Cov}(V^i_T,V^j_T)$ converges to
$\int_0^2f_{ij}(1,u)\,du$ as $T\to\infty$, and this limit equals
$\overline{\mu}_{2p}(g_i,g_j;\eta,\zeta)$. Since the limit of a
sequence of
covariance matrices is symmetric nonnegative, we have the result.
\end{pf*}

Now, we come to the aforementioned key result which consists of
proving the convergence we hinted at after the definition (\ref{KL1})
of the processes $L(g)$ and $L'(g)$. For a precise statement, we fix a
sequence $i_n$ of integers, and we associate the following
processes with $g$, an arbitrary function satisfying (\ref{1}):
%
%
\begin{eqnarray}\label{KL3}
\overline{L}(g)^n_t &=& \sqrt{k_n}\, \overline{W}(g)^n_{i_n+[k_nt]},
\qquad
\overline{L}{}'(g)^n_t=\sqrt{k_n} \overline{\chi
}(g)^n_{i_n+[k_nt]},\nonumber\\[-8pt]\\[-8pt]
\widehat{L}'(g)^n_t &=& k_n \widehat{\chi}(g)^n_{i_n+[k_nt]}.\nonumber
\end{eqnarray}
We do not mention the sequence $i_n$ in this notation, but those
processes clearly depend on it. In view of the ``approximation''
(\ref{AU1}), these processes (and in particular their conditional moments
of various orders) will play a central role in the sequel.

We fix a family $(g_l)_{1\leq l\leq d}$ of weight functions
satisfying (\ref{1}). We denote by $\overline{L}{}^n_t$ and $\overline
{L}{}'^n_t$ and $\widehat{L}
'^n_t$ the
$d$-dimensional processes with respective components, $\overline
{L}(g_l)^n_t$ and
$\overline{L}{}'(g_l)^n_t$ and $\widehat{L}'(g_l)^n_t$. These processes
can be
considered as
variables with values in the Skorokhod space $\mathbb{D}^d$ of all c\`
adl\`ag
functions from $\mathbb{R}_+$ into $\mathbb{R}^d$. The processes
$L_t$ and $L'_t$
with components $L(g_l)_t$ and $L'(g_l)_t$, defined by (\ref{KL1})
\textit{with the same Wiener processes $W^1$ and $W^2$} for all components,
are also $\mathbb{D}^d$-valued variables, and the probability on
$\mathbb{D}^{2d}=\mathbb{D}^d\times\mathbb{D}^d$ which is the law
of the pair, $(L,L')$ is
denoted by $R=R_{(g_v)}=R(dx,dy)$.

We also have a sequence $(f_n)$ of functions on $\mathbb{D}^{3d}$, all
depending on $w\in\mathbb{D}^{3d}$ only through their restrictions to
$[0,m+1]$ for
some $m\geq0$ and which satisfy the following property for some
$q'\geq2$ [below, $x,y,z\in\mathbb{D}^d$, so $v=(x,y)\in\mathbb
{D}^{2d}$ and
$(x,y,z)=(v,z)\in\mathbb{D}^{3d}$, and the same for $x',y',z'$ and $v'$;
moreover for any multidimensional Borel function $u$ on $\mathbb{R}_+$ we
put $u^\star_{m,n}=\|u(0)\|+\|u(m)\|+{\frac1{k_n} \sum_{i=1}^{(m+1)k_n}}\|u(i/k_n)\|$]:
%
%
%
\begin{equation}\label{KL4}
\left.
\begin{array}{l}
|f_n(v,z)| \leq  K
\bigl(1+(v_{m,n}^{\star})^{q'}+(z_{m,n}^{\star})^{q'/2}\bigr),\\
|f_n(v,z)-f_n(v',z')| \\
\qquad\leq
K\bigl((v-v')^\star_{m,n}+(z-z')^*_{m,n}\bigr)\\
\qquad\quad{}\times\bigl(1+(v_{m,n}^{\star})^{q'-1}
 + (v_{m,n}'^{\star})^{q'-1}+(z_{m,n}^{\star})^{q'/2-1}+
(z_{m,n}'^{\star})^{q'/2-1}\bigr).
\end{array}\right\}\hspace*{-35pt}
\end{equation}
\begin{lem}\label{LNN10} Assume Hypothesis \ref{hypoSNq} for some $q>4$ and
that $\sigma$ is bounded. Let $\Gamma$
be the set of all times $s\geq0$ such that both $\sigma$ and $\alpha$
are almost
surely continuous at time $s$. Let $z_0\in\mathbb{D}^d$ be the
constant function
with components $(\overline{g}_l'(2))_{1\leq l\leq d}$. Take any
sequence $(i_n)$ of
integers such that $s_n=i_n\Delta_n$ converges to some $s\in\Gamma$.
If the sequence $(f_n)$ satisfies (\ref{KL4}) for some $q'<q$ and
converges pointwise to a limit $f$, we have the almost sure convergence
%
%
\begin{equation}\label{KL5}\qquad
\mathbb{E} (f_n(\sigma_{s_n}\overline{L}{}^n,
\overline{L}{}'^n,\widehat{L}'^n)\mid\mathcal F_{s_n} ) \to
\int f(\theta\sigma_sx,\alpha_sy,2(\alpha_s)^2z_0) R(dx,dy).
\end{equation}
\end{lem}
\begin{pf}
(1) We first prove an auxiliary result. Let $\Omega
^{(0)}_s$
be the set of all $\omega^{(0)}$ such that both $\sigma(\omega
^{(0)})$ and
$\alpha(\omega^{(0)})$ are continuous at time $s$. We have
$\mathbb{P}^{(0)}(\Omega^{(0)}_s)=1$ because $s\in\Gamma$, and we fix
$\omega^{(0)}\in\Omega_s^{(0)}$. We consider the
probability space $(\Omega^{(1)},\mathcal F^{(1)},\mathbb{Q})$ where
$\mathbb{Q}=\mathbb{Q}(\omega^{(0)},\cdot)$, and our aim in this step is
to show that
under $\mathbb{Q}$,
%
%
\begin{equation}\label{KL6}
\overline{L}{}'^n \stackrel{\mathcal L}{\longrightarrow} \alpha_s\bigl(\omega^{(0)}\bigr)L',\qquad
\mathbb E_{\mathbb Q} ((\overline{L}{}'^n)^{*q}_{m,n}) \leq K_m
\end{equation}
(functional convergence in law in $\mathbb D^d$). In view of the definition of
$(\overline{L}{}'^n)^{*}_{m,n}$ (which is the norm of $t\rightarrow \overline{L}{}'^n_t$ described
above), the second property immediately follows from (\ref{LLN141}).

We first prove the finite-dimensional convergence. Let $0<t_1<\cdots
<t_r$. By (\ref{KL3}) and (\ref{RR2}) the $rd$-dimensional variable
$Z_n=(\overline{L}{}'^{n,l}_{t_i}\dvtx1\leq l\leq d,1\leq i\leq r)$ is
%
%
\begin{equation}\label{KL6676}
\left.
\begin{array}{l}
\displaystyle  Z_n = \sum_{j=1}^\infty z^n_j, \qquad\mbox{where }
z^n_j = \zeta^n_j a_j^n,
\zeta^n_j=\frac1{\sqrt{k_n}} \chi_{i_n+j-1}^n \quad\mbox{and}\\
a_j^{n,l,i} = \cases{
-k_n(g_l)'^n_{j-[k_nt_i]}, &\quad if
$1+[k_nt_i]\leq j\leq k_n+[k_nt_i]$,\cr
0, &\quad otherwise.}
\end{array}
\right\}
\end{equation}
Under $\mathbb{Q}$\vspace*{-1pt} the variables $\zeta^n_j$ are independent centered,
with $\mathbb{E}_Q(|\zeta^n_j|^4)\leq Kk_n^{-2}$ by Hypothesis \ref{hypoSNq}; recall
$q>4$. The
numbers $a_j^{n,l,i}$ being uniformly bounded and equal to $0$ when
$j>k_n+[k_ntr]$, we deduce that
under $\mathbb{Q}$ again the variables $z^n_j$ are independent with
\[
\mathbb{E}_Q(z^n_j) = 0,\qquad
\mathbb{E}_Q(\|z^n_j\|^4)\leq Kk_n^{-2},\qquad
\sum_{j=1}^\infty\mathbb{E}_Q(\|z^n_j\|^4) \to0.
\]
Next,
\[
\sum_{j=1}^\infty\mathbb{E}_Q(z^{n,l,i}_jz_j^{n,l',i'})
= \frac1{k_n} \sum_{j=1}^\infty\alpha_{(i_n+j-1)\Delta_n}\bigl(\omega^{(0)}\bigr)^2
a^{n,l,i}_j a_j^{n,l',i'}.
\]
On one hand $\alpha_{(i_n+j-1)\Delta_n}(\omega^{(0)})^2$ converges
uniformly in $j\leq k_n+[t_rk_n]$ to $\alpha_s(\omega^{(0)})^2$ because
$s\mapsto\alpha_s(\omega^{(0)})$ is continuous at $s$. On the other hand,
since $g_l=0$ outside $[0,1]$,
\begin{eqnarray*}
&&\frac1{k_n} \sum_{j=1}^\infty a^{n,l,i}_j a_j^{n,l',i'}\\
&&\qquad =
k_n\sum_{j=1}^\infty\int_{(j-1)/k_n}^{j/k_n}
g'_l\biggl(u-\frac{[k_nt_i]}{k_n}\biggr)\,du \int_{(j-1)/k_n}^{j/k_n}
g'_{l'}\biggl(u-\frac{[k_nt_{i'}]}{k_n}\biggr)\,du,
\end{eqnarray*}
which clearly converges to $c^{l,i,l',i'}=\int g'_l(v-t_i)g'_{l'}(v-t_{i'})\,dv$
by the mean value theorem, the piecewise continuity of each $g'_l$ and
Riemann approximation. Hence
%
%
\begin{equation}\label{KL6678}
\sum_{j=1}^\infty\mathbb{E}_Q(z^{n,l,i}_jz_j^{n,l',i'})
\to c^{l,i,l',i'} \alpha_s\bigl(\omega^{(0)}\bigr)^2.
\end{equation}
Then a standard limit theorem on row-wise independent triangular arrays
of infinitesimal variables yield that $Z_n$ converges\vspace*{2pt} in law under
$\mathbb{Q}
$ to a
centered Gaussian variable with covariance matrix
$(c^{l,i,l',i'} \alpha_s (\omega^{(0)})^2)$ (see, e.g., Theorem VII-2-36 of \cite{JS}). Now,
in view of (\ref{KL1}), this matrix is the covariance of the
centered Gaussian vector $(L'^{,l}_{t_i}\dvtx1\leq l\leq d,1\leq i\leq q)$,
and the finite-dimensional convergence in (\ref{KL6}) is proved.


To obtain the first property in (\ref{KL6}) it remains to
prove that for each $l$ the sequence of processes  $\overline{L}{}'(g_l)^n$ is
C-tight. Equivalently, we can prove that the sequence of processes
$G^n$ is C-tight, where $G^n$ is continuous, coincides with
$\overline{L}{}'(g_l)^n$ at all times $i/k_n$ and is piecewise linear between these
times. For
this we use a criterion given in \cite{IK} for example. Namely, since
$q>2$, the C-tightness of the sequence $G^n$ is implied by
\begin{equation}\label{KL7}
0\leq v\leq1 \quad\Rightarrow\quad
\mathbb{E}_\mathbb{Q} (|G^n_{t+v}-G^n_t|^q )\leq
Kv^{q/2}.
\end{equation}
A simple computation shows that
$G^n_{t+v}-G^n_t=\sum_j\delta^n_j\chi^n_j$
for suitable coefficients~$\delta^n_j$, such that at most $2[k_nv]$ are
smaller that $K_1/\sqrt{k_n}$, and
at most $k_n$ of them are smaller than $K_2v/\sqrt{k_n}$, and all others
vanish. Then the Burkholder--Davis--Gundy inequality yields
\begin{eqnarray*}
\mathbb{E}_\mathbb{Q} (|G^n_{t+v}-G^n_t|^q )
&\leq& K\mathbb{E}_\mathbb{Q} \biggl( \biggl(
\sum_j(\delta^n_j\chi^n_j)^2 \biggr)^{q/2} \biggr)\\
&\leq&
K\beta(q)\bigl(\omega^{(0)}\bigr) \bigl(K_1^q(2v)^{q/2}+K_2^qv^q \bigr),
\end{eqnarray*}
and (\ref{KL7}) follows. Then (\ref{KL6}) is completely proved.

(2) In exactly the same setting as in the previous step, we
prove here that
%
%
\begin{equation}\label{KL9}
\widehat{L}'^n \stackrel{\mathrm{u.c.p.}}{\longrightarrow} 2\bigl(\alpha
_s\bigl(\omega^{(0)}\bigr)\bigr)^2 z_0,\qquad
\mathbb{E}_\mathbb{Q} \Bigl(\sup_{v\leq t} \|\widehat{L}'^n_v\|
^{q/2} \Bigr) \leq K_t
\end{equation}
(under $\mathbb{Q}$ again, and with $z_0$ as in the statement of the
lemma). These
are componentwise properties, so we may assume $d=1$ here and
$g_1=g$.
The second property again follows from (\ref{LLN141}). For the first one, we see that under $\mathbb Q$
the variable $\zeta^n_{t,j}=k_n(g'(j/k_n) \Delta^n_{i_n+[k_nt]+j}\chi)^2$ satisfies
\begin{eqnarray*} a^n_{t,j}:\!&=&
\mathbb{E}_\mathbb{Q}(\zeta^n_{t,j})\\
&=& k_n\bigl(g'(j/k_n)\bigr)^2 \bigl(\bigl(\alpha
\bigl(\omega^{(0)}\bigr)_{(i_n+[k_nt]+j)\Delta_n}\bigr)^2\\
&&\hspace*{64.8pt}{} +\bigl(\alpha\bigl(\omega^{(0)}\bigr)_{(i_n+[k_nt]+j-1)\Delta_n}\bigr)^2 \bigr),\\
\mathbb{E}_\mathbb{Q}(|\zeta^n_{t,j}|^{q/2}) &\leq& K/k_n^{q/2}.
\end{eqnarray*}
%
%
%
In view of the continuity of $\alpha(\omega^{(0)})$ at time $s$
and of (\ref{SET4}), and since $\widehat{L}{}'^n_t=\sum_{j=1}^{k_n}\zeta^n_{t,j}$,
we see that $B^n_t=\mathbb{E}_{\mathbb{Q}}(\widehat{L}{}'^n_t)=\sum_{j=1}^{k_n}a^n_{t,j}$
converges locally uniformly to the ``constant'' $2(\alpha_s(\omega^{(0)}))^2 z_0$.
Hence it remains to prove that $V^n_t=\widehat{L}{}'^n_t-\break B^n_t
\stackrel{\mathrm{u.c.p.}}{\longrightarrow}0$.
For this it suffices to show $V^n_t\stackrel{\mathbb{P}}{\longrightarrow}0$ for each $t$, and the
C-tightness of both sequences $(\widehat{L}^n)$ and $(B^n)$, and the
latter follows from $B^n_t\stackrel{\mathrm{u.c.p.}}{\longrightarrow}2(\alpha_s(\omega^{(0)}))^2 z_0$.

Now, $V^n_t$ is the sum of the $k_n$ centered variables $\zeta^n_{t,j}-
a^n_{t,j}$, with $(q/2)$th absolute moment smaller than $K/k_n^{q/2}$, and
$\zeta^n_{t,j}$ is independent of $(\zeta^n_{t,l}\dvtx|l-j|\geq2)$. Then
obviously $\mathbb{E}_Q((V^n_t)^2)\leq K/k_n\to0$. For the C-tightness of
$(\widehat{L}^n)$ it suffices as in the end of Step 1 to prove the C-tightness
of the linearized versions $(G'^n)$ of~$(\widehat{L}{}'^n)$. We have
$G'^n_{t+v}-G'^n_t=\sum_i\delta^n_i(\Delta^n_i\chi)^2$ for suitable
coefficients\vspace*{1pt} $\delta^n_j$, such that at most $2[k_nv]$ are smaller that $K_1/k_n$,
and at most $k_n$ of them are smaller than $K_2v/k_n$, and all others
vanish. Then by the Burkholder--Davis--Gundy inequality (applied separately
for the sum of even indices and the sum of odd indices, to ensure the
independence of the summands), we have
\[
\mathbb{E}_{\mathbb{Q}}(|G'^n_{t+v}-G'^n_t|^{q/2})\leq K\mathbb{E}_{\mathbb{Q}} \biggl( \biggl(
\sum_j(\delta^n_j\chi^n_j)^2 \biggr)^{q/4} \biggr)\leq Kv^{q/4}.
\]
Since $q>4$, the C-tightness of $(G'^n)$ follows as in Step 1, and
(\ref{KL9}) holds.

(3) Now we draw some consequences of the previous facts. We set for
$y,z\in\mathbb{D}^d$,
\begin{eqnarray*}
f^{n}_{\omega^{(0)}}(y,z) &=&
f_n\bigl(\sigma_{s_n}\bigl(\omega^{(0)}\bigr)L^n\bigl(\omega^{(0)}\bigr),y,z\bigr),
\\
A^n_{j}\bigl(\omega^{(0)}\bigr)&=&\cases{
\displaystyle \int\mathbb{Q}\bigl(\omega^{(0)},d\omega^{(1)}\bigr) f^n_{\omega
^{(0)}}\bigl(\overline{L}{}
'^n\bigl(\omega^{(1)}\bigr),
\widehat{L}{}'^n\bigl(\omega^{(1)}\bigr)\bigr),&\quad $j=1$,\cr
\displaystyle \int f^n_{\omega^{(0)}}\bigl(\alpha_{s_n}\bigl(\omega^{(0)}\bigr)y,2\alpha_s^2
z_0\bigr) R(dx,dy), &\quad $j=2$.}
\end{eqnarray*}

The $\mathcal F^{(0)}$-measurable variables,
\[
\Phi_n = 1+\sup_{v\in[0,(m+1)u_n]} \sqrt{k_n} |W_{s_n+v}-W_{s_n}|,
\]
satisfy $\mathbb{E}(\Phi_n^u)\leq K_u$ for any $u>0$, by scaling of
the Brownian
motion $W$ whereas $\|\overline{L}{}^n_t\|\leq K\Phi_n$ if $t\leq m$.
Then we
deduce from (\ref{KL4}) and from the boundedness of $\sigma$ and
$\alpha
$ that if
$y,y',z,z'$ are in $\mathbb{D}^d$ and $u=(y,z)$ and $u'=(y',z')$.
%
%
\[
\left.
\begin{array}{l}
\bigl|f^n_{\omega^{(0)}}(u)\bigr| \leq K\Phi_n\bigl(\omega^{(0)}\bigr)^{q'}
\bigl(1+(y_{m,n}^{\star})^{q'}+(z_{m,n}^{\star})^{q'/2}\bigr),\\
\bigl|f^n_{\omega^{(0)}}(u)-f^n_{\omega^{(0)}}(u')\bigr|\leq K\Phi_n\bigl(\omega^{(0)}\bigr)^{q'}
(u-u')^\star_{m,n}\bigl(1+(y_{m,n}^{\star})^{q'-1}+(y_{m,n}'^{\star})^{q'-1}\\
\hspace*{207pt}{}+(z_{m,n}^{\star})^{q'/2-1}
+(z_{m,n}'^{\star})^{q'/2-1}\bigr).
\end{array}
\right\}
\]

Moreover $\alpha_{s_n}(\omega^{(0)})\to\alpha_{s}(\omega^{(0)})$,
so by the
Skorokhod representation theorem according to which, in case of
convergence in law, one can replace the original variables by variables
having the same laws and converging pointwise, one deduces from
(\ref{KL6}) and (\ref{KL9}) [these
imply that the variables $f^n_{\omega^{(0)}}(\overline{L}{}'^n,\widehat
{L}'^n)$ are
uniformly integrable, since $q'<q$], that
%
%
\begin{equation}\label{KL11}
\left.
\begin{array}{l}
\omega^{(0)}\in\Omega_s^{(0)} \quad\Rightarrow\quad
A^n_{1}\bigl(\omega^{(0)}\bigr)-A^n_{2}\bigl(\omega^{(0)}\bigr) \to0,\\
\mathbb{E} (|A^n_{j}|^{q/q'} ) \leq K.
\end{array}
\right\}
\end{equation}

Next, we make the following observation: due to the
$\mathcal F^{(0)}$-conditional independence of the $\chi_t$'s, a
version of
the conditional expectation in (\ref{KL5}) is
$\mathbb{E}(A^n_{1}\mid\mathcal F_{s_n})$.
Therefore in view of (\ref{KL11}) (which ensures the
uniform integrability and the a.s. convergence to $0$ of the sequence
$A^n_{1}-A^n_{2}$), (\ref{KL5}) is implied by
%
%
\begin{equation}\label{KL12}
\mathbb{E}(A^n_{2}\mid\mathcal F_{s_n}) \to F(\sigma_s,\alpha
_s) \qquad\mbox{a.s.},
\end{equation}
where
\[
F(\eta,\zeta) = \int f(\theta\eta x,\zeta y,2(\zeta)^2z_0) R(dx,dy).
\]

(4) For proving (\ref{KL12}) we start again with an auxiliary
result, namely
%
%
\begin{equation}\label{KL13}
\overline{L}{}^n \stackrel{\mathcal L}{\longrightarrow} \theta L.
\end{equation}
For this, we see that
$Z_n=(\overline{L}{}^{n,l}_{t_i}\dvtx1\leq l\leq d,1\leq i\leq r)$ is
given by
(\ref{KL6676}), except that
\[
\zeta^n_j = \sqrt{k_n} \Delta^n_{i_n+j}W,\qquad
a_j^{n,l,i} = \cases{
(g_l)^n_{j-[k_nt_i]}, &\quad if
$1+[k_nt_i]\leq j\leq k_n+[k_nt_i]$,\cr
0, &\quad otherwise.}
\]
Then the proof of (\ref{KL13}), both for the finite-dimensional
convergence and the C-tightness, is exactly the same as for (\ref{KL6})
[note that the right-hand side of (\ref{KL6678}) is now
$\theta^2\int g_l(v-t_i)g_{l'}(v-t_{i'})\,dv$ which is the covariance matrix
of $(\theta L^{l}_{t_i}\dvtx1\leq l\leq d,1\leq i\leq r)$]. 
Further, since $\Vert\overline L{}^n_t\Vert\leq K \Phi_n$ if $t\leq m$,
%
%
\begin{equation}\label{KL14}
\mathbb E \Bigl( {\sup_{v\leq t}} \Vert\overline L{}^n_v\Vert^q \Bigr)\leq K_t.
\end{equation}

(5) Now we introduce some functions on $\mathbb{R}^2$:
\begin{eqnarray*}
F_n(\eta,\zeta)&=&\int\mathbb{E} (f_n(\eta\overline{L}{}^n,
\zeta y,2(\zeta)^2z_0) ) R(dx,dy), \\
F'_n(\eta,\zeta)&=&\int\mathbb{E} (f_n(\theta\eta L,\zeta
y,2(\zeta)^2z_0)
) R(dx,dy).
\end{eqnarray*}
Under $R$ the canonical process is locally in time bounded in each
$\mathbb{L}^r$. Then in view of (\ref{KL4}) we deduce from (\ref
{KL13}) and
(\ref{KL14}), and exactly as for (\ref{KL11}), that
$F_n-F'_n\to0$ locally uniformly in $\mathbb{R}^2$. We also deduce from
(\ref{KL4}) that $F'_n(\eta_n,\zeta_n)-F'_n(\eta,\zeta)\to0$ whenever
$(\eta_n,\zeta_n)\to(\eta,\zeta)$, and also that $F'_n\to F$ pointwise
because $f_n\to f$ pointwise, and hence we have $F_n(\eta_n,\zeta_n)
\to F(\eta,\zeta)$.

At this point it remains to observe that, because
$(W_{s_n+t}-W_{s_n})_{t\geq0}$ is independent of $\mathcal F_{s_n}$,
we have
$\mathbb{E}(A^n_{2}\mid\mathcal F_{s_n})=F_n(\sigma_{s_n},\alpha
_{s_n})$. Since
$(\sigma_{s_n},\alpha_{s_n})\to(\sigma_s,\alpha_s)$ a.s., we
readily deduce
(\ref{KL5}), and we are done.
\end{pf}
\begin{rem}\label{RKL1} In the previous lemma, suppose that all
$f_n$ (hence $f$ as well) only depend on $(x,y)$ and not on $z$;
that is, the processes $\widehat{L}'^n$ do not enter the picture. Then
it is easily
seen from the
previous proof that we do not need $q>4$, but only $q>2$.
\end{rem}

\subsection{Asymptotically negligible arrays}\label{ssec-ANA}

An array $(\delta^n_i)$ of \textit{nonnegative} variables is called AN
(for ``asymptotically negligible'') if
%
%
\begin{equation}\label{LNN2}
\sqrt{\Delta_n} \sup_{0\leq j\leq k_n} \mathbb{E} \Biggl(\sum_{i=0}^{[t/u_n]}
\delta^n_{ik_n+j} \Biggr) \to0,\qquad |\delta^n_i| \leq K,
\end{equation}
for all $t>0$. With any process $\gamma$ (in the sequel, $\gamma$
will usually be
$\gamma=\sigma$ or $\gamma=\alpha$) and any integer $m$ we
associate the variables
\[
\Gamma(\gamma,m)^n_i = {\sup_{t\in[i\Delta_n,i\Delta_n+(m+1)u_n]}}
|\gamma_t-\gamma_{i\Delta_n}|,\qquad
\Gamma'(\gamma,m)^n_i = \mathbb{E}(\Gamma(\gamma,m)^n_i\mid
\mathcal F^n_i).
\]
\begin{lem}\label{LLNN1} \textup{(a)} If $(\delta^n_i)$ is an AN array, we have
%
%
\begin{equation}\label{LNN4}
\Delta_n\mathbb{E} \Biggl(\sum_{i=1}^{[t/\Delta_n]}\delta^n_i
\Biggr) \to0
\end{equation}
for all $t>0$, and the array $((\delta^n_i)^q)$ is also AN for each $q>0$.

\textup{(b)} If $\gamma$ is a c\`adl\`ag bounded process, then for all $m\geq1$ the
two arrays $(\Gamma(\gamma,m)^n_i)$ and $(\Gamma'(\gamma,m)^n_i)$
are AN.
\end{lem}
\begin{pf}
(a) The left-hand side of (\ref{LNN4}) is smaller than
a constant times the left-hand side of (\ref{LNN2}), hence the first
claim. The second claim follows from H\"older's inequality if $q<1$, and
from $\sum_{i\in I}(\delta^n_i)^q\leq K\sum_{i\in I}\delta^n_i$ if $q>1$
(recall that $|\delta^n_i|\leq K$).

(b) Let $\delta^n_i=\Gamma(\gamma,m)^n_i$. If $\varepsilon>0$,
denote by $N(\varepsilon)_t$ the number of jumps of $\gamma$ with
size bigger
than $\varepsilon$ on the interval $[0,t]$, and by $v(\varepsilon
,t,\eta)$
the supremum of $|\gamma_s-\gamma_r|$ over all pairs $(r,s)$ with
$s\leq
r\leq s+\eta$ and $s\leq t$ and such that $N(\varepsilon
)_s-N(\varepsilon)_r=0$. Since
$\gamma$ is bounded,
\[
u_n \sup_{0\leq j\leq
k_n} \mathbb{E} \Biggl(\sum_{i=0}^{[t/u_n]}\delta^n_{ik_n+j} \Biggr)
\leq\mathbb{E} \bigl(t v\bigl(\varepsilon,t+1,(m+1)u_n\bigr)+(Kt)\wedge(K
u_n N(\varepsilon
)_{t+1}) \bigr)
\]
as soon as $(m+2)u_n\leq1$. Since $\limsup_{n\to\infty
} v(\varepsilon,t+1,(m+1)u_n)
\leq\varepsilon$, Fatou's lemma implies that the $\limsup$ of the
left-hand side
above is smaller than $Kt\varepsilon$, so we have (\ref{LNN2})
because $\varepsilon$
is arbitrarily small. Since $\mathbb{E}(\Gamma'(\gamma
,m)^n_i)=\mathbb{E}(\Gamma(\gamma
,m)^n_i)$,
the second claim follows.
\end{pf}

\subsection{Some estimates}\label{ssec-ES}

Here we provide a number of
estimates under the following assumption for some $q>2$:
%
%
\begin{equation}\label{L4N0}
\begin{tabular}{p{328pt}}
$\bullet$ we have (\ref{CIS}) and Hypothesis \ref{hypoSNq}, and $b$ and $\sigma
$ are
bounded, and $\sigma$ \\
\hspace*{8pt}and $\alpha$ are c\`adl\`ag.
\end{tabular}\hspace*{-42pt}
\end{equation}
This list of estimates is quite tedious, but unfortunately they
play a central role in many places in the sequel.
We first introduce some
notation where $i$ and $j$ are integers, $Y$ is an arbitrary process
and $\rho_{p,l}$ is given by (\ref{L3N2}) and $i+j\geq1$ in the first
line below, and $p$ an even integer in (\ref{L3N401}):
\begin{eqnarray}\label{L3N401}\qquad
&&\hspace*{15.30pt}\kappa^n_{i,j} =
\sigma^n_i\Delta^n_{i+j}W+\Delta^n_{i+j}\chi,\nonumber\\
&&\hspace*{14.65pt}\lambda^n_{i,j} = \Delta^n_{i+j}Z-\kappa^n_{i,j}
 = \Delta
^n_{i+j}X-\sigma
^n_i\Delta^n_{i+j}W,
\nonumber\\
&&\overline{\kappa}(g)^n_{i,j}=\sum_{l=1}^{k_n-1}g^n_l\kappa
^n_{i,j+l},\qquad
\overline{\lambda}(g)^n_{i,j}=\sum_{l=1}^{k_n-1}g^n_l\lambda
^n_{i,j+l},\nonumber\\[-8pt]\\[-8pt]
&&\hspace*{0.26pt}\widehat{\lambda}(g)^n_{i,j}=\sum_{l=1}^{k_n}(g'^n_l\lambda^n_{i,j+l})^2,
\nonumber\\
&&\left.
\begin{array}{l}
\displaystyle\phi(Y,g,p)^n_i = \sum_{l=0}^{p/2}\rho_{p,l} (\overline{Y}(g)^n_i)^{p-2l}
(\widehat{Y}(g)^n_i)^l,\nonumber\\
\phi(g,p)^n_{i,j} =
\displaystyle\sum_{l=0}^{p/2}\rho_{p,l} (\overline{\kappa
}(g)^n_{i,j})^{p-2l} (\widehat{\chi}(g)^n_{i,j})^l.
\end{array}
\right\}\nonumber
\end{eqnarray}
Recalling (\ref{AU1}), we see that $\overline{\kappa}(g)^n_{i,j}$ is
an approximation
of $\overline{Z}(g)^n_i$, and its asymptotic behavior is described in Lemma
\ref{LNN10}. Then $\phi(g,p)^n_{i,j}$ is an approximation of
$\phi(Y,g,p)^n_i$ whereas by (\ref{L3N4}) we have
\[
\overline{V}(Y,g,p)^n_t = \sum_{i=0}^{[t/\Delta_n]-k_n}\phi(Y,g,p)^n_i.
\]
One of the aims of the estimates below is to prove that these
approximations induce a negligible error.

In the forthcoming inequalities, we have $0\leq j\leq
mk_n$ where $m$ is a fixed integer.
First, if we use (\ref{LLN14}) and the boundedness of $g$, and also
(\ref{LLN141}), we obtain for $u>0$
%
%
\begin{eqnarray}\label{L4N2}\qquad
&&\left.
\begin{array}{l}
\mathbb{E}\bigl(|\overline{X}(g)^n_i|^u+|\overline{W}(g)^n_i|^u\mid
\mathcal F^n_i\bigr)\leq K_u\Delta_n^{u/4},\\
u\leq q \quad\Rightarrow\quad
\mathbb{E}\bigl(|\overline{Z}(g)^n_i|^u+|\overline{\kappa}(g)^n_i|^u\mid
\mathcal F^n_i\bigr)\leq K_u\Delta_n^{u/4},
\end{array}
\right\}
\\
\label{L4N3}
&&\left.
\begin{array}{l}
\displaystyle\lambda^n_{i,j} = \int_{(i+j-1)\Delta_n}^{(i+j)\Delta_n}
\bigl(b_s\,ds+(\sigma_s-\sigma^n_i)\,dW_s \bigr),\\[8pt]
\displaystyle\overline{\lambda}(g)^n_{i,j} = \int_{(i+j)\Delta
_n}^{(i+j+k_n)\Delta_n}
g_n\bigl(s-(i+j)\Delta_n\bigr) \bigl(b_s\,ds+(\sigma_s-\sigma^n_i)\,dW_s \bigr).
\end{array}
\right\}
\end{eqnarray}
Hence we obtain for $u\geq1$, and
recalling that $\Gamma(\sigma,m)^n_i\leq K$,
%
%
\begin{equation}\label{L4N4}
\left.
\begin{array}{l}
\mathbb{E}(|\lambda^n_{i,j}|^u\mid\mathcal F^n_i )\leq
K_u\Delta_n^{u/2} \bigl(\Delta_n^{u/2}+\Gamma'(\sigma,m)^n_i
\bigr),\\
\mathbb{E}(|\overline{\lambda}(g)^n_{i,j}|^u\mid\mathcal F^n_i
)\leq
K_u\Delta_n^{u/4} \bigl(\Delta_n^{u/4}+\Gamma'(\sigma,m)^n_i \bigr).
\end{array}
\right\}
\end{equation}
If $u$ is an \textit{odd integer}, (\ref{L4N2}), (\ref{L4N4})
and an expansion of $(\sigma^n_i\overline{W}(g)^n_i+\overline
{\lambda}(g)^n_i)^u$ yield
%
%
\begin{eqnarray}\label{L4N5}
\mathbb{E}((\overline{W}(g)^n_i)^u\mid\mathcal
F^n_i)&=&0,\nonumber\\[-8pt]\\[-8pt]
\mathbb{E}((\overline{X}(g)^n_i)^u\mid\mathcal F^n_i )&\leq&
K_u\Delta_n^{u/4} \bigl(\Delta_n^{1/4}+\sqrt{\Gamma'(\sigma
,m)^n_i} \bigr).\nonumber
\end{eqnarray}

Next, using $|g'^n_i|\leq K/k_n$ and (\ref{LLN14}) and the first part
of (\ref{L4N4}), plus H\"older's inequality and the definition of
$\widehat{Y}(g)^n_i$, plus the obvious fact that
$\mathbb{E}(|\kappa^n_{i,j}|^u\mid\mathcal F^n_i)\leq K_u$ if $u\leq
q$, and after some
calculations, we get for $u\geq1$
%
%
\begin{equation}\label{L4N6}\qquad
\left.
\begin{array}{l}
\mathbb{E}\bigl(|\widehat{X}(g)^n_i|^u+|\widehat{W}(g)^n_i|^u\mid
\mathcal F^n_i\bigr)\leq K_u\Delta_n^{3u/2},\\
u\leq q/2 \quad\Rightarrow\quad
\mathbb{E}\bigl(|\widehat{Z}(g)^n_{i+j}|^u+|\widehat{\chi
}(g)^n_{i+j}|^u\mid\mathcal F^n_i\bigr)\leq
K_u\Delta_n^{u/2},\\
u\leq q \quad\Rightarrow\quad
\mathbb{E}\bigl(|\widehat{Z}(g)^n_{i+j}-\widehat{\chi}(g)^n_{i+j}|^u\mid
\mathcal F^n_i\bigr)\leq K_u\Delta_n^u.
\end{array}
\right\}
\end{equation}

Then, if we combine (\ref{L4N2}), (\ref{L4N4}) and (\ref{L4N6}),
and use again H\"older's inequality, we obtain for all reals $l,u\geq1$
and $r\geq0$,
%
%
\begin{equation}\label{L4N7}
\left.
\begin{array}{l}
(l+2r)u\leq q \\
\quad\mbox{$\Rightarrow$}\quad
\mathbb{E} \bigl(
|(\overline{Z}(g)^n_{i+j})^l(\widehat{Z}(g)^n_{i+j})^r\\[2pt]
\hspace*{47.6pt}{}- (\overline{\kappa}(g)^n_{i,j})^l(\widehat{\chi}(g)^n_{i,j})^r
|^u\mid\mathcal F^n_i \bigr)\\[4pt]
\qquad\quad\qquad\leq K_{u,l,r}\Delta_n^{ul/4+ur/2} \bigl(\Delta_n^{u/4}+
(\Gamma'(\sigma,m)^n_i)^{1-u(l+2r-1)/q} \bigr),\\[2pt]
2ru\leq q \quad\Rightarrow\quad
\mathbb{E} \bigl( |(\widehat{Z}(g)^n_{i+j})^r-(\widehat{\chi
}(g)^n_{i,j})^r |^u\mid\mathcal F
^n_i \bigr)
\leq K_{u,r}\Delta_n^{ru/2+u/2}.
\end{array}
\right\}\hspace*{-32pt}
\end{equation}
Finally, by (\ref{L3N401}), this readily gives for $p\geq2$ an even
integer and
$u\geq1$ a real, such that $pu\leq q$,
%
%
\begin{eqnarray}\label{L4N8}
&&\mathbb{E}\bigl(|\phi(Z,g,p)^n_{i+j}|^u+|\phi(g,p)^n_{i,j}|^u\mid
\mathcal F^n_i\bigr)\nonumber\\
&&\qquad
\leq K_{u,p}\Delta_n^{pu/4},\nonumber\\[-8pt]\\[-8pt]
&&\mathbb{E}\bigl(|\phi(Z,g,p)^n_{i+j}-\phi(g,p)^n_{i,j}|^u\mid\mathcal
F^n_i\bigr)\nonumber\\
&&\qquad\leq K_{u,p}\Delta_n^{pu/4} \bigl(\Delta_n^{u/4}+
(\Gamma'(\sigma,m)^n_i)^{1-u(p-1)/q} \bigr).\nonumber
\end{eqnarray}

\subsection[Proof of Theorem 3.3]{Proof of Theorem \protect\ref{TLNN1}}\label{ssec:T3}

By localization we can and will assume (\ref{L4N0}). We set
\begin{eqnarray*}
\mu^n_i&=&\Delta_n^{-p/4} |\overline{Z}(g)^n_i|^p,\qquad
\zeta^n_i=\Delta_n^{-p/4}|\overline{\kappa}(g)^n_{i,0}|^p,\\
\gamma_t&=&m_p \biggl|\theta
\overline{g}(2)\sigma_t^2+\frac{\overline{g}{}'(2)}{\theta} \alpha
_t^2 \biggr|^{p/2}.
\end{eqnarray*}
The left-hand side of (\ref{LNN5}) is $\sum_{i=0}^{[t/\Delta
_n]-k_n}\mu^n_i$ whereas
we deduce from (\ref{L4N7}) with $r=0$ and Lemma \ref{LLNN1} that
$\Delta_n\sum_{i=0}^{[t/\Delta_n]-k_n}|\mu^n_i-\zeta^n_i|\stackrel
{\mathrm{u.c.p.}}{\longrightarrow}0$. Then it remains to prove
%
%
\begin{equation}\label{LLN25}
\Delta_n\sum_{i=0}^{[t/\Delta_n]-k_n}\zeta^n_i \stackrel{\mathrm
{u.c.p.}}{\longrightarrow} \int_0^t\gamma_s \,ds.
\end{equation}

Set $\zeta'^n_i=\mathbb{E}(\zeta^n_i\mid\mathcal F^n_i)$. By (\ref{L4N2}),
$\mathbb{E}((\zeta^n_i)^2\mid\mathcal F^n_i)\leq K$, and in
particular $\zeta'^n_i\leq K$. Moreover $\zeta^n_i$ is
$\mathcal F^n_{i+k_n}$-measurable, hence $\mathbb{E}((\zeta
^n_{i}-\zeta'^n_{i})
(\zeta^n_{j}-\zeta'^n_{j}))=0$ if $|j-i|\geq k_n$, and
\begin{eqnarray*}
\mathbb{E} \Biggl( \Biggl|\Delta_n\sum_{i=0}^{[t/\Delta_n]-k_n}(\zeta
^n_i-\zeta'^n_i) \Biggr|^2 \Biggr) &=&
\Delta_n^2\sum_{i,j=1}^{[t/\Delta_n]-k_n}\mathbb{E} \bigl((\zeta
^n_i-\zeta'^n_i)
(\zeta^n_j-\zeta'^n_j) \bigr)\\
&\leq& K\Delta_nk_n \to0.
\end{eqnarray*}
Thus it is enough to prove (\ref{LLN25}) with $\zeta^n_i$ substituted
with $\zeta'^n_i$. Since $\gamma_t+\zeta'^n_i\leq K$,
\[
\Biggl|\Delta_n\sum_{i=0}^{[t/\Delta_n]-k_n}\zeta'^n_i-\int
_0^t\gamma_s \,ds \Biggr| \leq
\int_{\Delta_n}^{([t/\Delta_n]-k_n)\Delta_n}|\gamma^n_s-\gamma
_s|\,ds+Kk_n\Delta_n,
\]
where $\gamma^n_s=\zeta'^n_i$ when $(i-1)\Delta_n\leq s<i\Delta_n$.
Therefore,
since $|\gamma_s^n-\gamma_s|\leq K$, in order to obtain (\ref
{LLN25}) it is
enough to prove that for Lebesgue-almost all $s$ we have
$\gamma^n_s\to\gamma_s$ a.s.
In particular it is enough to prove that, for all $s\in\Gamma$ (cf. Lemma
\ref{LNN10}), we have
%
%
\begin{equation}\label{LLN59}
\zeta'^n_{[s/\Delta_n]+1} \to\gamma_s \qquad\mbox{a.s.}
\end{equation}

With the notation of Lemma \ref{LNN10}, we take $d=1$ and the
weight function \mbox{$g_1=g$}, and the functions $f_n=f$ on $\mathbb{D}^3$ as
$f(x,y,z)=|x(0)+y(0)|^p$, so (\ref{KL4}) is satisfied with $q'=p<q$ and $m=0$.
Moreover we fix $s\in\Gamma$ and set
$i_n=[s/\Delta_n]+1$, so $s_n=i_n\Delta_n\to s$. The left-hand
side of (\ref{KL5}) is $\Delta_n^{p/4}k_n^{p/2}\zeta'^n_{i_n}$ whereas
its right-hand side is $\mathbb{E}' (|\theta\eta L(g)_0+\eta'
L'(g)_0|^p )$ [recall
(\ref{KL1})] evaluated at $\eta=\sigma_s$ and $\eta'=\alpha_s$.
Since $L(g)_0$
and $L'(g)_0$ are independent centered normal with respective variances
$\overline{g}(2)$ and $\overline{g}{}'(2)$, this right-hand side is
$m_p (\theta^2\overline{g}(2)\sigma^2_s
+\overline{g}{}'(2)\alpha^2_s )^{p/2}=\theta^{p/2}\gamma_s$. Since
$\Delta_n^{p/4}k_n^{p/2}\to\theta^{p/2}$, we get (\ref{LLN59}).

\subsection[Proof of Theorem 3.4]{Proof of Theorem \protect\ref{TL3N1}}\label{ssec:T4b}

As we said already, (a) is a particular case of (\ref{LLN111}) when
$p\geq4$, and of (\ref{LLN112}) when $p=2$. For (b), we can again
assume (\ref{L4N0}). We set
\begin{eqnarray*}
\mu^n_i&=&\Delta_n^{-p/4}\phi(Z,g,p)^n_i,\qquad
\zeta^n_i=\Delta_n^{-p/4}\phi(g,p)^n_{i,0},\\
\gamma_t&=&m_p (\theta\overline{g}(2))^{p/2} |\sigma_t|^p,
\end{eqnarray*}
and $\zeta'^n_i=\mathbb{E}(\zeta^n_i\mid\mathcal F^n_i)$. We deduce
from (\ref{L4N8}) and
Lemma \ref{LLNN1} that $\Delta_n\times\break\sum_{i=0}^{[t/\Delta_n]-k_n}|\mu
^n_i-\zeta^n_i|\stackrel{\mathrm{u.c.p.}}{\longrightarrow}0$. Then
it is
enough to prove (\ref{LLN25}).

By (\ref{L4N8}) we have $\mathbb{E}((\zeta^n_i)^2\mid\mathcal
F^n_i)\leq K$; hence
$\zeta'^n_i\leq K$. Then, exactly as in the previous proof, it remains
to show (\ref{LLN59}) when $s\in\Gamma$. For this, we use Lemma
\ref{LNN10} with $d=1$ and $g_1=g$ and the functions $f_n=f$ given by
\[
f(x,y,z) =
\sum_{l=0}^{p/2}\rho_{p,l}|x(0)+y(0)|^{p-2l} |z(0)|^l.
\]
The left-hand side of (\ref{KL5}) is again
$\Delta_n^{p/4}k_n^{p/2}\zeta'^n_{i_n}$. Its right-hand side is
$\mu_p(g$;\break $\theta\sigma_s,\alpha_s)$, as given by (\ref{L3N612}),
and by
(\ref{L3N613}) this is also $\theta^{p/2}\gamma_s$. Then (\ref{LLN9})
holds.

\subsection[Proof of Theorem 3.6]{Proof of Theorem \protect\ref{TL3N7}}\label{ssec:T6}

The proof is basically the same as in the previous subsection, using
again Lemma \ref{LNN10} and the fact that we deal with asymptotically
$k_n$-dependent variables. We can assume (\ref{L4N0}), and we have
\begin{eqnarray*}
&&\overline{\mu}_{2p}(g,h;\eta,\zeta)\\
&&\qquad=\sum_{r,r'=0}^{p/2}\rho
_{p,r}\rho_{p,r'}
(2\zeta^2\overline{g}{}'(2))^r(2\zeta^2\overline{h}{}'(2))^{r'}
\\
&&\qquad\quad\hspace*{20.5pt}{}\times\bigl(m_{p-2r,p-2r'}(g,h;\eta,\zeta)
-2m_{p-2r}(g;\eta,\zeta)m_{p-2r'}(h;\eta,\zeta) \bigr).
\end{eqnarray*}
Therefore is is enough to prove that for $r,r'$ between $0$ and $p/2$,
and with the notation
\begin{eqnarray*}
\mu^n_i&=&\Delta_n^{-p/2}(\widehat{Z}(g)^n_i)^r (\widehat{Z}(h)^n_i)^{r'}
\Biggl(|\overline{Z}(g)^n_{i+k_n}|^{p-2r}
\frac1{k_n}\sum_{j=1}^{2k_n}|\overline{Z}(h)^n_{i+j}|^{p-2r'}\\
&&\hspace*{141.1pt}{}
-2|\overline{Z}(g)^n_i|^{p-2r}|\overline{Z}(h)^n_{i+k_n}|^{p-2r'}
\Biggr),\\
\gamma_t&=&\theta^{-p/2}(2\alpha_t^2\overline{g}{}'(2))^{r+r'}
\bigl(m_{p-2r,p-2r'}
(g,h;\theta\sigma_t,\alpha_t)\\
&&\hspace*{92.8pt}{}
-2m_{p-2r}(g;\theta\sigma_t,\alpha_t)m_{p-2r'}(h;\theta\sigma
_t,\alpha_t) \bigr),
\end{eqnarray*}
we have
\[
\Delta_n\sum_{i=0}^{[t/\Delta_n]-3k_n}\mu^n_i \stackrel{\mathrm
{u.c.p.}}{\longrightarrow} \int_0^t\gamma_s\,ds.
\]

By (\ref{L4N6}) and (\ref{L4N7}) we have ${\Delta_n\sum
_{i=0}^{[t/\Delta_n]-3k_n}}|\mu^n_i-\zeta^n_i|
\stackrel{\mathrm{u.c.p.}}{\longrightarrow}0$ where
\begin{eqnarray*}
\zeta^n_i&=&\Delta_n^{-p/2}(\widehat{\chi}(g)^n_i)^r (\widehat
{\chi}(h)^n_i)^{r'}
\Biggl(|\overline{\kappa}(g)^n_{i,k_n}|^{p-2r} \frac1{k_n}\sum_{j=1}^{2k_n}
|\overline{\kappa}(h)^n_{i,j}|^{p-2r'}\\
&&\hspace*{131.1pt}{}
-2|\overline
{\kappa}(g)^n_{i,0}|^{p-2r}
|\overline{\kappa}(h)^n_{i,k_n}|^{p-2r'} \Biggr),
\end{eqnarray*}
so we are left to prove
\[
\Delta_n\sum_{i=0}^{[t/\Delta_n]-3k_n}\zeta^n_i \stackrel{\mathrm
{u.c.p.}}{\longrightarrow} \int_0^t\gamma_s\,ds.
\]

We set $\zeta'^n_i=\mathbb{E}(\zeta^n_i\mid\mathcal F^n_i)$, so as in
the proof of Theorem \ref{TLNN1} it is enough to prove (\ref{LLN59})
when $s\in\Gamma$. We apply Lemma \ref{LNN10} with $d=2$ and $g_1=g$
and $g_2=h$ and the functions $f_n$ and $f$ on $\mathbb{D}^6$ defined by
\begin{eqnarray*}
&&f_n((x,x'),(y,y'),(z,z'))\\
&&\qquad=z(0)^r z'(0)^{r'} \Biggl(
|x(1)+y(1)|^{p-2r}\frac1{k_n}\sum_{j=1}^{2k_n}
\biggl|x'\biggl(\frac j{k_n}\biggr)+y'\biggl(\frac j{k_n}\biggr)\biggr|^{p-2r'}\\
&&\qquad\quad\hspace*{93.7pt}{}-2|x(0)+y(0)|^{p-2r}|x'(1)+y'(1)|^{p-2r'} \Biggr),\\
&&f((x,x'),(y,y'),(z,z'))\\
&&\qquad=z(0)^r z'(0)^{r'}
\biggl(|x(1)+y'(1)|^{p-2r}\int_0^2|x'(t)+y'(t)|^{p-2r'}\,dt\\
&&\qquad\quad\hspace*{69.4pt}{} -2|x(0)+y(0)|^{p-2r} |x'(1)+y'(1)|^{p-2r'} \biggr)
\end{eqnarray*}
and again $i_n=[s\Delta_n]+1$.
Then (\ref{KL4}) is satisfied with $q'=2p<q$ and \mbox{$m=1$}, and $f_n\to f$
pointwise. The left-hand side of (\ref{KL5}) is
$\Delta_n^{p/2}k_n^{p}\zeta'^n_{i_n}$, whereas its right-hand side is
$\theta^{p/2}\gamma_s$ [recall that\vspace*{1pt} $(L(g)_0,L'(g)_0)$ and
$(L(h)_1,L'(h)_1)$ are independent]. Since
$\Delta_n^{p/2}k_n^p\to\theta^{p/2}$, we get (\ref{LLN59}) by the lemma,
and the proof is finished.

\subsection{Auxiliary results on the noise process}\label{ssec:ARN}
At this stage we start the proof of our CLTs, and this is done
through a large number of steps. In the first, crucial step
we derive some estimates on the (conditional) moments of the noise process
$\chi$. Recall that $\mathcal G^n_i$ denotes the $\sigma$-field
generated by
$\mathcal F^{(0)}$ and $\mathcal F^n_i$. Set
%
%
\begin{equation}\label{L3N8}
A(g)^n_i = \sum_{j=1}^{k_n}(g'^n_j)^2(\alpha^n_{i+j-1})^2.
\end{equation}
For random variables $U_\gamma$ and $V_\gamma$ indexed by a parameter
$\gamma$ [for example, $\gamma=(n,i)$ just below], with $V_\gamma>0$,
we write $U_\gamma=\mathrm{O}_u(V_\gamma)$ if the family $U_\gamma
/V_\gamma$ is bounded
in probability.
\begin{lem}\label{LL3N1} Assume Hypothesis \ref{hypoSNq} for some $q\geq2$, and
let $v$ and $r$ be integers such that $2\leq v+2r\leq q$. Let also
$m\geq0$ and $j$ be arbitrary in $\{0,1,\ldots,mk_n\}$.

\textup{(a)} When $v$ is even we have
%
%
\begin{eqnarray}\label{L3N10}\qquad
&&\mathbb{E} ((\overline{\chi}(g)^n_{i+j})^v(\widehat{\chi
}(g)^n_{i+j})^r\mid\mathcal G^n_i )\nonumber\\[-8pt]\\[-8pt]
&&\qquad=
m_v 2^r (A(g)^n_{i+j})^{r+v/2}+\mathrm{O}_u(\Delta_n^{r/2+v/4+1/2})
\nonumber\\
\label{L3N100}
&&\qquad= m_v 2^r\frac{\overline
{g}'(2)^{r+v/2}}{k_n^{r+v/2}} (\alpha^n_i)^{2r+v}+
\mathrm{O}_u \bigl(\Delta_n^{r/2+v/4} \bigl(\Delta_n^{1/2}+\Gamma
(\alpha,m)^n_i\bigr)\bigr).
\end{eqnarray}

\textup{(b)} When $v$ is odd we have
%
%
\begin{equation}\label{L3N11}
\mathbb{E} ((\overline{\chi}(g)^n_{i+j})^v(\widehat{\chi
}(g)^n_{i+j})^r\mid\mathcal G^n_i ) =
\mathrm{O}_u(\Delta_n^{r/2+v/4+1/4}),
\end{equation}
and also, for some suitable numbers $\gamma_{v,r}$, depending on $g$,
%
%
\begin{eqnarray}\label{L3N12}
&&\mathbb{E} ((\overline{\chi}(g)^n_{i+j})^v(\widehat{\chi
}(g)^n_{i+j})^r\mid\mathcal G^n_i )\nonumber\\
&&\qquad=
\frac{\gamma_{v,r}}{k_n^{r+v/2+1/2}} (\alpha^n_i)^{2r+v-3} \beta(3)^n_i
\\
&&\qquad\quad{}
+\mathrm{O}_u \bigl(\Delta_n^{r/2+v/4+1/4} \bigl(\Delta_n^{1/4}
+\Gamma(\alpha,m)^n_i+\Gamma(\beta(3),m)^n_i \bigr) \bigr).\nonumber
\end{eqnarray}
\end{lem}
\begin{pf}
Equations (\ref{L3N100}) and (\ref{L3N11}) are simple
consequences of (\ref{L3N10})\break and (\ref{L3N12}), respectively, upon
observing\vspace*{-2pt} that $A(g)^n_{i+j}=\overline{g}{}'(2)(\alpha^n_i)^2/k_n+\break
\mathrm{O}_u(\Delta_n^{1/2}(\Delta_n^{1/2}+\Gamma(\alpha
,m)^n_i))$. As for
(\ref{L3N10}) and (\ref{L3N12}), and up to taking a further
conditional expectation, it is enough to prove them when $j=0$, so
in the rest of the proof we take $j=0$, and thus $m=0$ as well.
The product $(\overline{\chi}(g)^n_i)^v(\widehat{\chi}(g)^n_i)^r$
is the sum of
all the terms of the form
%
%
\begin{equation}\label{L3N103}
\left.
\begin{array}{l}
\displaystyle\Phi(J,n)=(-1)^v\prod_{l=1}^vg'^n_{j_l}\chi^n_{i+j_l-1}
\prod_{l=1}^s(g'^n_{j'_l}\chi^n_{i+j'_l+\overline{j}{}'_l-1})^2\\
\hspace*{45.9pt}{}\times
\displaystyle\prod_{l=1}^{r-s} (-2(g'^n_{j''_l})^2\chi^n_{i+j''_l}\chi^n_{i+j''_l-1}
),\\
\displaystyle J=\{s,j_1,\ldots,j_v,j'_1,\ldots,j'_s,\overline{j}{}'_1,\ldots
,\overline{j}{}'_s,
j''_1,\ldots,j''_{r-s}\},\\
\qquad\mbox{where } s\in\{0,\ldots,r\},
j_l,j'_l,j''_l\in\{1,\ldots,k_n\},
\overline{j}{}'_l\in\{0,1\}.
\end{array}
\right\}
\end{equation}
We denote by $I(J)$ the family of all indices of the
variables $\chi^n_j$ occurring in (\ref{L3N103}), the index $j$
appearing $l$ times if $\chi^n_j$ is taken at the power $l$, so that
$I(J)$ contains $v+2r$ indices. We also denote by $D(u)^n$ the class of
all $J$'s such that among the $v+2r$ indices in $I(J)$, there are
exactly $u$ different indices, each one appearing at least twice. Note
that $D(u)^n$ is the disjoint union over $s'=0,\ldots,r$ of the set
$D(u,s')^n$ of all $J\in D(u)^n$ such that $s=s'$. Note also that
$D(u)^n=\varnothing$ if $u>v/2+r$.

By (\ref{SET11}) and the $\mathcal F^{(0)}$-conditional independence of
the $\chi_t$'s, the conditional expectation
$\mathbb{E}(\Phi(J,n)\mid\mathcal G^n_i)$ is always smaller than
$K/k_n^{v+2r}$, and
it vanishes if $J$ is outside $\bigcup_{u\geq1}D(u)^n$; that is,
\[
\mathbb{E} ((\overline{\chi}(g)^n_i)^v(\widehat{\chi
}(g)^n_i)^r\mid\mathcal G^n_i ) =
\sum_{u=1}^{[v/2]+r}\overline{\Phi}^n_u,
\]
where
\[
\overline{\Phi}^n_u = \sum_{s=0}^r\overline{\Phi}(u,s)^n,\qquad
\overline{\Phi}(u,s)^n=\sum_{J\in D(u,s)^n}\mathbb{E} (\Phi
(J,n)\mid\mathcal G^n_i ).
\]
Now $\# D(u,s)^n\leq Kk_n^u$, so $|\overline{\Phi}(u,s)^n|\leq Kk_n^{u-v-2r}$;
hence $\overline{\Phi}^n_u=\mathrm{O}_u(\Delta_n^{r/2+v/4+1/4})$ as
soon as $u\leq
r-1/2+v/2$. We deduce that for proving (\ref{L3N10}), so $v$ is even,
it is enough to show that $\overline{\Phi}^n_u$ equals the right-hand
side of
(\ref{L3N10}), for $u=r+v/2$. In the same way, for proving (\ref{L3N12}),
so $v$ is odd, it is enough to show that $\overline{\Phi}^n_u$ equals the
right-hand side
of (\ref{L3N12}) for $u=r+(v-1)/2$.

(a) Suppose that $v$ is even and $u=r+v/2$. The definition of $D(u)^n$ and
the property $u=r+v/2$ yield that, if $J\in D(u)^n$, there is a nonnegative
integer $w\leq\frac v2\wedge\frac{r-s}2$ such that $\Phi(J,n)$ is the
product of $\frac{v+s+r-w}2$ terms, of three types, all for different
indices for $\chi^n$:

(1) $s-w+\frac v2$ terms of the form $(g'^n_j\chi^n_{i+j-1})^2$
or $(g'^n_j\chi^n_{i+j})^2$;

(2) $w$ terms of the form $-2(g'^n_j)^3g'^n_{j+1}(\chi^n_{i+j-1}
\chi^n_{i+j})^2$;

(3) $\frac{r-s-w}2$ terms of the form $4(g'^n_j)^4(\chi^n_{i+j-1}
\chi^n_{i+j})^2$.

Hence $\# D(u,s)^n\leq Kk_n^{({v+s+r})/2}$ because the
number of terms for a particular $J$ is smaller than $\frac{v+s+r}2$,
and the
indices range from $1$ to $k_n$. Moreover, since $\alpha$ is bounded and
$|g'^n_j|\leq K/k_n$, we have
$\mathbb{E}(|\Phi(J,n)|\mid\mathcal G^n_i)\leq K/k_n^{v+2r}$. We
then deduce that
%
%
\begin{equation}\label{XX0}
|\overline{\Phi}(u,s)^n| \leq Kk_n^{({v+s+r})/2-v-2r} \leq
K\Delta_n^{r/2+v/4+({r-s})/2}.
\end{equation}

In particular, $\overline{\Phi}(u,s)^n=$ O$(\Delta_n^{r/2+v/4+1/2})$ when
$s<r$, and it thus remains to prove that $\overline{\Phi}(u,r)^n$ is equal
to the right-hand side of (\ref{L3N10}). If $J\in D(u,r)^n$, then
$\Phi(J,n)$ contains only terms of type (1). In fact $D(u,r)^n$
contains exactly the families $J$ for which $s=r$, and among
$j_1,\ldots,j_v$ there are $v/2$ distinct indices, each one
appearing twice (we then denote by $J_1$ the set of the $v/2$
distinct indices), and the sets $J_2=\{j'_l+\overline{j}{}'_l\dvtx1\leq
l\leq
r,\overline{j}{}'_l=0\}$ and $J_3=\{j'_l+\overline{j}{}'_l\dvtx1\leq l\leq
r,\overline{j}{}'_l=1\}$
have distinct indices, and $J_1$, $J_2$ and $J_3$ are pair-wise
disjoint. With this notation, we have (with $u$ terms all together
in the products)
%
%
\begin{equation}\label{XX6}
\mathbb{E}(\Phi(J,n)\mid\mathcal G^n_i) =
\prod_{j\in J_1\cup J_2}(g'^n_j\alpha^n_{i+j-1})^2 \prod_{j\in
J_3}(g'^n_{j-1}\alpha^n_{i+j-1})^2.
\end{equation}
The assumption (\ref{1}) on
$g$ yields that $|g'^n_j-g'^n_{j-1}|\leq K/k_n^2$, except for
$j$ belonging to the set $Q_n$ of indices for which $g'$
fails to exist or to be Lipschitz on $[(j-1)/k_n,jk_n]$, so
$\# Q_n\leq K$. Since
$\alpha^n_i\leq K$, we thus have
\[
\mathbb{E}(\Phi(J,n)\mid\mathcal G^n_i)=\cases{
\displaystyle\prod_{j\in J_1\cup J_2\cup J_3}(g'^n_j\alpha^n_{i+j-1})^2+
\mathrm{O}_u(k_n^{-2u-1}),\cr
\hspace*{70.2pt} \mbox{if $Q_n\cap(J_1\cup J_2\cup
J_3)=\varnothing$},
\cr
\mathrm{O}_u(k_n^{-2u}), \qquad \mbox{otherwise.}}
\]
Consider
now $L=\{l_1,\ldots,l_u\}$ in the set $\mathcal L_n$ of all families
of indices
with $1\leq l_1<\cdots<l_u\leq k_n$, and let $w_n(L)$ be the number of
$J\in D(u,r)^n$ such that the associated sets $J_1,J_2,J_3$ satisfy
$J_1\cup J_2\cup J_3=L$. Then since $\# D(u,r)^n\leq Kk_n^u$ and
$\sup_n\# Q_n<\infty$, we deduce from the above that
%
%
\begin{equation}\label{XX1}
\overline{\Phi}(u,r)^n = \sum_{L\in\mathcal L_n}w_n(L)
\prod_{j\in L}(g'^n_j\alpha^n_{i+j-1})^2+\mathrm{O}_u(\Delta_n^{u/2+1/2}).
\end{equation}
Now we have to evaluate $w_n(L)$. There are $C_u^r$ many ways of
choosing the two complementary subsets, $J_1$ and $J_2\cup J_3$, of
$L$. Next,
with $J_1$ given, there are $(v/2)! (v-1)(v-3)\cdots3 \cdot1$ ways
of choosing
the indices $j_l$ so that $j_1,\ldots,j_v$ has $v/2$ paired distinct
indices which are the indices in $J_1$, and
we recall that $(v-1)(v-3)\cdots3 \cdot1=m_v$ (if $v=0$ then $J_1$ is
empty and there is $m_0=1$ ways again of choosing $J_1$). Finally
with $J_2\cup J_3$ fixed, there are $2^r r!$ ways of choosing the indices
$j'_l+\overline{j}{}'_l$, all of them different, when the smallest index in
$J_2\cup J_3$ is
bigger than $1$, and $2^{r-1} r!$ ways if this smallest index is $1$.
Summarizing, we get
%
%
\begin{equation}\label{XX2}
w_n(L) = \cases{
m_v 2^r u!, &\quad if
$1\notin L$,\cr
m_v 2^{r-1} u!, &\quad if $1\in L$.}
\end{equation}

On the other hand, we have by (\ref{L3N8})
\[
(A(g)^n_i)^u = u! \sum_{L\in\mathcal L_n}\prod_{j\in L}(g'^n_j\alpha
^n_{i+j-1})^2
+\mathrm{O}_u(k_n^{-1-u}).
\]
Therefore, by (\ref{XX1}) and (\ref{XX2}), we deduce that
\[
m_v 2^r (A(g)^n_i)^u-\overline{\Phi}(u,r)^n = m_v 2^{r-1}
\sum_{L\in\mathcal L_n\dvtx1\in L} \prod_{j\in L}(g'^n_j\alpha^n_{i+j-1})^2
+\mathrm{O}_u(\Delta_n^{u/2+1/2}).
\]
Since\vspace*{-1pt} $|g'^n_j|\leq K/k_n$ and since the number of $L\in\mathcal L_n$ such
that $1\in L$ is smaller than~$k_n^{u-1}$, the right-hand\vspace*{1pt} side above is
smaller than $K\Delta_n^{u/2+1/2}$, and we deduce that $\overline
{\Phi}(u,r)^n$
is equal to the right-hand side of (\ref{L3N10}). In view of (\ref{XX0}),
this completes the proof of (\ref{L3N10}).

(b) Suppose that $v$ is odd and $u=r+v/2-1/2$, and recall that we
need to prove that $\overline{\Phi}^n_u$ equals the right-hand side
of (\ref{L3N12}).
Again, the definition of $D(u)^n$ and the property $u=r+v/2-1/2$ yield that,
if $J\in D(u)^n$, there is a number $a$ in $\{0,1\}$ and a nonnegative
integer $w\leq\frac{v-1}2\wedge\frac{r-s-2a}2$ such that $\Phi
(J,n)$ is
the product of $\frac{v+s+r-w-1}2$ terms, all for different indices for
$\chi^n$ with $s-w+a+\frac{v-3}2$ terms of type 1, $w$ terms of type 2,
$\frac{r-s-w-2a}2$ terms of type 3 and $1-a$ and $a$ term, respectively,
of the types (4) and (5) described below:

(4) terms of the form $(g'^n_j\chi_{i+j-1}^n)^3$ or
$(g'^n_j)^2g'^n_{j+1}(\chi^n_{i+j})^3$,

(5) terms of the form $-2(g'^n_j)^4g'^n_{j+1}(\chi_{i+j-1}^n)^3
(\chi^n_{i+j})^2$ or $-2(g'^n_j)^3(g'^n_{j+1})^2\times\break
(\chi^n_{i+j-1})^2(\chi^n_{i+j})^3$,
the whole product being multiplied by $-1$.
It follows that $\# D(u,s)^n\leq Kk_n^{({v+s+r-1})/2}$ by the
same argument as in (a) whereas\break
$\mathbb{E}(|\Phi(J,n)|\mid\mathcal G^n_i)\leq K/k_n^{v+2r}$ still
holds. Hence, instead
of (\ref{XX0}) we get\break
$|\overline{\Phi}(u,s)^n|\leq K\Delta_n^{r/2+v/4+1/4+({r-s})/2}$.
In particular, $\overline{\Phi}(u,s)^n= \mathrm{O}(\Delta_n^{r/2+v/4+1/2})$
when $s<r$,
and it thus remains to prove that $\overline{\Phi}(u,r)^n$ is equal
to the
right-hand side of (\ref{L3N12}).

If $J\in D(u,r)^n$ then $\Phi(J,n)$ has $u-1$ terms of type (1) and
one of type (4), and there is exactly one common index among
$j_1,\ldots,j_v$
and $j'_1+\overline{j}{}'_1,\ldots,j'_s+\overline{j}{}'_s$. In other
words, we can associate
with $J$ three sets, $J_1,J_2,J_3$, pairwise disjoint [with the same
description than when $v$ is even, except that $\# J_1=\frac{v-1}2$ and
$\#(J_2\cup J_3)=r-1$], plus an index $l$ outside $J_1\cup J_2\cup J_3$ and
an integer
$\overline{l}$ equal to $0$ or $1$, such that instead of (\ref{XX6})
we have
\begin{eqnarray*}
\mathbb{E}(\Phi(J,n)\mid\mathcal
G^n_i) &=& -(g'^n_l)^2g'^n_{l+\overline{l}} \beta
(3)^n_{i+l+\overline{l}-1}
\\
&&{}\times\prod_{j\in J_1\cup J_2}(g'^n_j\alpha^n_{i+j-1})^2
\prod_{j\in J_3}(g'^n_{j-1}\alpha^n_{i+j-1})^2.
\end{eqnarray*}
This is equal to
\[
-\beta(3)^n_i (\alpha^n_i)^{2u-2} (g'^n_l)^3
\prod_{j\in J_1\cup J_2\cup J_3}(g'^n_j)^2,
\]
up to $\mathrm{O}_u(k_n^{-2u}(k_n^{-1}+\Gamma(\alpha,m)^n_i+\Gamma
(\beta
(3),m)^n_i))$ when
$Q_n\cap(\{k\}\cup J_1\cup J_2\cup J_3)=\varnothing$ and to
$\mathrm{O}_u(k_n^{-2u})$, otherwise. Therefore, since $\#
D(u,r)^n\leq Kk_n^u$,
we deduce that
\[
\overline{\Phi}(u,r)^n = -\beta(3)^n_i (\alpha^n_i)^{2u-2}
\sum_{l,J_1,J_2,J_3}(g'^n_l)^3
\prod_{j\in J_1\cup J_2\cup J_3}(g'^n_j)^2+R_n,
\]
where the remainder term $R_n$ is like the last term in (\ref{L3N12}),
and the sum is extended over all $l,J_1,J_2,J_3$ such that
$\{l\},J_1,J_2,J_3$ are pairwise disjoint in the set $\{1,\ldots,k_n\}$.
Then with $R'_n$ as $R_n$ above, we have
\[
\overline{\Phi}(u,r)^n = -\beta(3)^n_i (\alpha^n_i)^{2u-2}
\Biggl(\sum_{j=1}^{k_n}(g'^n_j)^3 \Biggr)
\Biggl(\sum_{j=1}^{k_n}(g'^n_j)^2 \Biggr)^{u-1}+R'_n.
\]
Then by an estimate similar to (\ref{SET4}) (without the absolute value),
we deduce (\ref{L3N12}), with $\gamma_{v,r}=-\overline{g}(2)^{r+v/2-1/2}
\int_0^2(g'(s))^3\,ds$.
\end{pf}
\begin{lem}\label{LL3N3} Assume Hypothesis \ref{hypoSNq} for some $q\geq2$, and
let $p$ be an even integer. With the notation (\ref{L3N401}) for
$\phi(g,p)$, the variables
%
%
\begin{equation}\label{AU2}
\Psi(g,p)^n_{i,j} = \mathbb{E} (\phi(g,p)^n_{i,j}\mid\mathcal
G^n_i )
-(\sigma^n_i\overline{W}(g)^n_{i+j})^p
\end{equation}
satisfy, for all $u\leq q/p$ and $m\geq0$ and
$0\leq j\leq mk_n$,
%
%
\begin{eqnarray}\label{L3N210}\hspace*{43pt}
|\mathbb{E}(\Psi(g,p)^n_{i,j}\mid\mathcal F^n_i)| &\leq& K\Delta
_n^{p/4+1/4} \bigl(\Delta_n^{1/4}+
\Gamma'(\alpha,m)^n_i+\Gamma'(\beta(3),m)^n_i \bigr),
\\
\label{L3N211}
\mathbb{E}(|\Psi(g,p)^n_{i,j}|^u\mid\mathcal F^n_i) &\leq& K\Delta
_n^{up/4+u/4}.
\end{eqnarray}
\end{lem}
\begin{pf}
In view of (\ref{L3N401}), and recalling that
$\sigma^n_i\overline{W}(g)^n_{i+j}$ is $\mathcal G^n_i$-measurable,
we see that
\begin{eqnarray*}
\mathbb{E}(\phi(g,p)^n_{i,j}\mid\mathcal G^n_i)&=&
\sum_{r=0}^{p/2}\sum_{w=0}^{p-2r}C_{p-2r}^w\rho_{p,r}
(\sigma^n_i\overline{W}(g)^n_{i+j})^w \\
&&\hspace*{35.11pt}{}\times\mathbb{E} ((\overline
{\chi}(g)_{i+j}^n)^{p-2r-w}
(\widehat{\chi}(g)^n_{i+j})^r\mid\mathcal G^n_i ).
\end{eqnarray*}
By (\ref{L3N2}) and a change of the order of summation, we easily get
\begin{eqnarray*}
&&\sum_{r=0}^{p/2}\sum_{v=0}^{p/2-r}C_{p-2r}^{2v} \rho_{p,r}
2^r m_{p-2r-2v} (\sigma^n_i\overline{W}(g)^n_{i+j})^{2v}
(A(g)^n_{i+j})^{p/2-v} \\
&&\qquad= (\sigma^n_i\overline{W}(g)^n_{i+j})^p;
\end{eqnarray*}
hence
\begin{eqnarray*}
\Psi(g,p)^n_{i,j}&=&\sum_{r=0}^{p/2}\sum_{v=0}^{p/2-r}C_{p-2r}^{2v}
\rho_{p,r} (\sigma^n_i\overline{W}(g)^n_{i+j})^{2v}\\
&&\hspace*{39.1pt}{}\times \bigl(\mathbb
{E} (
(\overline{\chi}(g)_{i+j}^n)^{p-2r-2v} (\widehat{\chi
}(g)^n_{i+j})^r\mid\mathcal G^n_i )\\
&&\hspace*{84.1pt}{} - 2^r m_{p-2r-2v} (A(g)^n_{i+j})^{p/2-v} \bigr)\\
&&\hspace*{0pt}{} + \sum_{r=0}^{p/2}\sum_{v=0}^{p/2-r-1}C_{p-2r}^{2v+1}
\rho_{p,r} (\sigma^n_i\overline{W}(g)^n_{i+j})^{2v+1}\\
&&\hspace*{62.3pt}{}\times  \mathbb
{E} ((\overline{\chi}(g)_{i+j}^n)^{p-2r-2v-1} (\widehat{\chi
}(g)^n_{i+j})^r\mid\mathcal G^n_i ).
\end{eqnarray*}
Now (\ref{L3N210}) is a simple consequence of (\ref{L4N2}) and
(\ref{L3N10}) applied to the terms in the first sum above and of
(\ref{L4N5}) and (\ref{L3N12}) for those in the second sum. Finally,
(\ref{L3N211}) follows from (\ref{L4N2}), (\ref{L3N10}) and (\ref{L3N11}),
plus H\"older's inequality.
\end{pf}

\subsection{Block splitting}\label{BS}

In this subsection we start the proof of Theorem~\ref{TCC1}. Due to
overlapping intervals the summands involved in the definition
of $V(Z,g$, $p,r)^n_t$ are asymptotically $k_n$-dependent variables, and
we will use the (classical) block splitting method to ensure some
``conditional'' independence. Namely, we split the sum over $i$ in
the definition of $\overline{V}(Z,g,p)^n_t$ into big blocks of
size $mk_n$ ($m$~is an integer which will eventually go to
$\infty$) which are separated by small blocks of size $k_n$. The big
blocks become asymptotically conditionally independent, and the
small blocks become negligible as $m\to\infty$. In a second step we prove
a CLT for big blocks, for any fixed $m$. We then obtain the result by
standard methods.

Here we fix the integer $m\geq2$. Recalling (\ref{L3N401}),
the $i$th block of size $mk_n$ contains $\phi(Z,g,p)^n_j$
for all $j$ between $I(m,n,i)=(i-1)(m+1)k_n$ and $I(m,n,i)+mk_n-1$. In
a similar way, the $i$th block of size $k_n$ corresponds to indices
$j$ between $\overline{I}(m,n,i)=i(m+1)k_n-k_n$ and $\overline
{I}(m,n,i)+k_n-1$. The number
of pairs of blocks which can be accommodated without using data
after time $t$ is then
$i_n(m,t)= [\frac{t-(k_n-1)\Delta_n}{(m+1)k_n\Delta_n} ]$.
The ``real'' times corresponding to the beginnings of the $i$th big
and small blocks are then $t(m,n,i)=I(m,n,i)\Delta_n$ and
$\overline{t}(m,n,i)=\overline{I}(m,n,i)\Delta_n$.

At this stage, we need some more notation. The summands in $\overline
{V}(Z,g,p)_t$
are the $\phi(Z,g,p)^n_i$, but we will indeed show that they can be
replaced by $\phi(g,p)_{i,j}^n$ [see (\ref{L3N401}) for suitable choice
of $i$]. This leads us to consider the partial sums
(we drop the mention of $p$, but we keep the function $g$)
\begin{eqnarray*}
\zeta(g,m)^n_i &=& \sum_{j=0}^{mk_n-1}\phi(Z,g,p)^n_{I(m,n,i)+j},\\
\overline{\zeta}(g,m)^n_i &=& \sum_{j=0}^{k_n-1}\phi
(Z,g,p)^n_{\overline{I}(m,n,i)+j},\\
\delta(g,m)^n_i &=& \sum_{j=0}^{mk_n-1}\phi(g,p)^n_{I(m,n,i),j},\\
\overline{\delta}(g,m)^n_i &=& \sum_{j=0}^{k_n-1}\phi(g,p)^n_{\overline
{I}(m,n,i),j},\\
U(g,m)^n_t &=& \sum_{i=1}^{i_n(m,t)}\bigl(\zeta(g,m)^n_i-\delta(g,m)^n_i\bigr),\\
\overline{U}(g,m)^n_t &=& \sum_{i=1}^{i_n(m,t)}\bigl(\overline{\zeta
}(g,m)^n_i-\overline{\delta}(g,m)^n_i\bigr),
\\
U'(g,m)^n_t &=& \sum_{i=i_n(m,t)(m+1)k_n}^{[t/\Delta_n]-k_n}\phi(Z,g,p)^n_i.
\end{eqnarray*}
Consider the discrete time filtrations $\mathcal F(m)^n_i=\mathcal
F_{I(m,n,i+1)}^n$
and $\overline{\mathcal F}(m)^n_i=\break\mathcal F_{\overline
{I}(m,n,i+1)}^n$. Observe
that $\delta(g,m)^n_i$ is $\mathcal F(m)^n_i$-measurable and
$\overline{\delta}(g,m)^n_i$ is $\overline{\mathcal
F}(m)^n_i$-measurable, and set
\begin{eqnarray*}
\gamma(g,m)^n_i &=& \mathbb{E}(\delta(g,m)^n_i\mid\mathcal
F(m)^n_{i-1}),\qquad
\overline{\gamma}(g,m)^n_i = \mathbb{E}(\overline{\delta
}(g,m)^n_i\mid\overline{\mathcal F}(m)^n_{i-1}),
\\
D(g,m)^n_t &=& \sum_{i=1}^{i_n(m,t)}\gamma(g,m)^n_i,\qquad
N(g,m)^n_t=\sum_{i=1}^{i_n(m,t)}\bigl(\delta(g,m)^n_i-\gamma(g,m)^n_i\bigr),\\
\overline{D}(g,m)^n_t &=&
\sum_{i=1}^{i_n(m,t)}\overline{\gamma}(g,m)^n_i,\qquad
\overline{N}(g,m)^n_t=\sum_{i=1}^{i_n(m,t)}\bigl(\overline{\delta
}(g,m)^n_i-\overline{\gamma}(g,m)^n_i\bigr).
\end{eqnarray*}
The key point is the following obvious relation, for any $m\geq1$:
%
%
\begin{eqnarray}\label{CC15}
\overline{V}(Z,g,p)^n_t&=&N(g,m)^n_t+D(g,m)^n_t+\overline{N}(g,m)^n_t+
\overline{D}(g,m)^n_t\nonumber\\[-8pt]\\[-8pt]
&&{} + U(g,m)^n_t+\overline{U}(g,m)^n_t+U'(g,m)^n_t;\nonumber
\end{eqnarray}
the contribution of the big blocks being $N(g,m)^n_t+D(g,m)^n_t$,
whereas $\overline{N}(g,\break m)^n_t+\overline{D}(g,m)^n_t$ accounts for
the small blocks,
and $U(g,m)^n_t$ and $\overline{U}(g,m)^n_t$ being asymptotically
negligible, whereas $U'(g,m)^n_t$ is a border term.
Note that $D(g,m)^n_t$ is a sort
of drift which asymptotically cancels with the centering term in
(\ref{CC1}). To be more specific, the leading term for the CLT is the
martingale $N(g,m)^n_t$, and we will eventually prove a CLT for it
and the negligibility of the rest in the sense that
%
%
\begin{eqnarray}\label{AXX}\qquad
&&\varepsilon,t>0 \quad\Rightarrow\quad \lim_{m\to\infty}\limsup_{n\to
\infty}
\mathbb{P} \Bigl(\sup_{s\leq t}
|\widetilde{V}(Z,g,p)^n_s\nonumber\\[-8pt]\\[-8pt]
&&\hspace*{160pt}{}-\Delta
_n^{3/4-p/4}N(g,m)^n_s |
>\varepsilon\Bigr) = 0.\nonumber
\end{eqnarray}
\begin{lem}\label{LCC2} Under (SN-$p$) we have
$\Delta_n^{3/4-p/4} U'(g,m)^n_t\stackrel{\mathbb{P}}{\longrightarrow
}0$ as $n\to\infty$.
\end{lem}
\begin{pf}
The variable
$U'(g,m)^n_t$ is the sum of at most $(m+1)k_n$ terms $\phi(Z,g,p)^n_i$,
all of them satisfying (\ref{L4N8}). Then the expectation of the
absolute value of $\Delta_n^{3/4-p/4} U'(g,m)^n_t$ is less than
$K_{m,p}k_n\Delta_n^{3/4}$ which clearly goes to $0$.
\end{pf}
\begin{lem}\label{LCC3} Under (SN-$2p$) we have, as $n\to\infty$,
and for each fixed $m$,
\[
\Delta_n^{3/4-p/4} U(g,m)^n \stackrel{\mathit
{u.c.p.}}{\longrightarrow} 0,\qquad
\Delta_n^{3/4-p/4} \overline{U}(g,m)^n \stackrel{\mathit
{u.c.p.}}{\longrightarrow} 0.
\]
\end{lem}
\begin{pf}
(1) The proofs of both claims
are the same, and we prove, for example, the first one. With the
notation $\eta^n_i=\zeta(g,m)^n_i-\delta(g,m)^n_i$ and $\eta'^n_i=
\mathbb{E}(\eta^n_i\mid\mathcal F(m)^n_{i-1})$, we have
\[
U(g,m)^n_t = U^{n,1}_t+U^{n,2}_t,
\]
where
\[
U^{n,1}_t = \sum_{i=1}^{i_n(m,t)}(\eta^n_i-\eta'^n_i),\qquad
U^{n,2}_t=\sum_{i=1}^{i_n(m,t)}\eta'^n_i.
\]
Then we need to prove
%
%
\begin{equation}\label{EXT}
\Delta_n^{3/4-p/4} U^{n,k} \stackrel{\mathrm
{u.c.p.}}{\longrightarrow} 0,\qquad k=1,2.
\end{equation}

(2) By the inequalities of Doob and Cauchy--Schwarz,
\begin{eqnarray*}
\mathbb{E} \Bigl(\sup_{s\leq t}|U_s^{n,1}|^2 \Bigr)&\leq&
4\sum_{i=1}^{i_n(m,t)}\mathbb{E} (|\eta^n_i|^2 )\\
&\leq&
4mk_n\sum_{i=1}^{i_n(m,t)} \sum_{j=0}^{mk_n-1}
\mathbb{E}\bigl(|\phi(Z,g,p)^n_{i+j}-\phi(g,p)^n_{i,j}|^2\bigr).
\end{eqnarray*}
By (\ref{L4N8}) and Lemma \ref{LLNN1}, the right-hand side above, multiplied
by $\Delta_n^{3/2-p/2}$, goes to $0$, so (\ref{EXT}) for $k=1$ follows.

For (\ref{EXT}) with $k=2$, and by virtue of (\ref{L3N401}), and
dropping $g$ from the notation,
we see that it is enough to show that, for all integers $l$
between $0$ and $p/2$, we can find an AN array $(\delta^n_i)$ (depending
on $l$) such that
\begin{eqnarray*}
&&0\leq j\leq mk_n \quad\Rightarrow\quad
\bigl|\mathbb{E} \bigl((\overline{Z}{}^n_{i+j})^{p-2l}(\widehat
{Z}^n_{i+j})^l-(\overline{\kappa}^n_{i,j})^{p-2l}
(\widehat{\chi}^n_{i+j})^l\mid\mathcal F^n_i \bigr) \bigr|\\
&&\hspace*{95.2pt}\qquad \leq
K\Delta
_n^{p/4+1/4} \delta^n_{i+j}.
\end{eqnarray*}
When $l=p/2$ the second estimate (\ref{L4N7}) with $u=1$ gives the
result, but otherwise the first estimate (\ref{L4N7}) with $u=1$ is
not quite
enough. Below we fix $l$ between $0$ and $p/2-1$, and the result will
be true if we have the following:
%
%
\begin{equation}\label{CC352}
|\mathbb{E}(F^n_{i,j}\mid\mathcal F^n_i)| \leq K\Delta
_n^{p/4+1/4} \delta^n_{i+j},
\end{equation}
where
\[
F^n_{i,j} = \cases{
(\overline{\kappa}^n_{i,j})^{p-2l}\bigl((\widehat{Z}^n_{i+j})^l-(\widehat
{\chi}^n_{i+j})^l\bigr),
&\quad (called Case A),\cr
(\widehat{\chi}^n_{i+j})^l\bigl((\overline{Z}{}^n_{i+j})^{p-2l}-(\overline
{\kappa}^n_{i,j})^{p-2l}\bigr),
&\quad (called Case B),\cr
\bigl((\overline{Z}{}^n_{i+j})^{p-2l}-(\overline{\kappa}^n_{i,j})^{p-2l}\bigr)
\bigl((\widehat{Z}^n_{i+j})^l-(\widehat{\chi}^n_{i,j})^l\bigr), &\quad
(called Case C),}
\]
and where again $(\delta^n_i)$ is an AN array (perhaps different for
each case).

In Cases\vspace*{-1pt} A and C we have $F^n_{i,j}=0$ when $l=0$, and we have
$|\mathbb{E}(F^n_{i,j}\mid\mathcal F^n_i)|\leq K\Delta_n^{p/4+1/2}$
when $l\geq1$
[apply (\ref{L4N2}) and the second part of (\ref{L4N7}), plus the fact
that $\Gamma'(\sigma,m+1)^n_i\leq K$, and the Cauchy--Schwarz inequality],
hence (\ref{CC352}) with $\delta^n_i=\Delta_n^{1/4}$.

(3) Now we consider Case B. Recall that
$\overline{Z}{}^n_{i+j}=\overline{\kappa}^n_{i,j}+\overline{\lambda
}{}^n_{i,j}$; hence
\[
F^n_{i,j} = \sum_{u=1}^{p-2l}C_{p-2l}^u G^{u,n}_{i,j},\qquad
G^{u,n}_{i,j} = (\widehat{\chi}^n_{i+j})^l
(\overline{\kappa}^n_{i,j})^{p-2l-u} (\overline{\lambda}{}^n_{i,j})^u
\]
and we will prove (\ref{CC352}) separately for each
$G^{u,n}_{i,j}$. For this, we begin with a decomposition of $\overline{\lambda}{}^n_{i,j}$.
Recall (\ref{CC002}) and the boundedness of the coefficients. By
(\ref{L4N3}) we have $\overline{\lambda}{}^n_{i,j}=\xi^n_{i,j}+\xi
'^n_{i,j}$ where,
with the simplifying notation $S=i\Delta_n$ and $T=(i+j)\Delta_n$,
\begin{eqnarray*}
\xi^n_{i,j} &=& \int_{T}^{T+u_n} g_n\bigl(s-(i+j)\Delta_n\bigr)\\
&&\hspace*{29.8pt}{}\times \biggl(
(b_s-b^n_i)\,ds+ \biggl(\int_S^s \bigl(\widetilde{b}_r\, dr
+(\widetilde{\sigma}_r-\widetilde{\sigma}^n_i)\,dW_r \bigr)
\biggr) \,dW_s \biggr),
\\
\xi'^n_{i,j} &=& \int_{T}^{T+u_n}g_n\bigl(s-(i+j)\Delta_n\bigr) \bigl(
b^n_i\,ds+\widetilde{\sigma}^n_i(W_s-W_S)\,dW_s+(M_s-M_S)\,dW_s \bigr).
\end{eqnarray*}
Then for $v\geq1$ and $j\leq mk_n$, we have
%
%
\begin{equation}\label{CC36}
\left.
\begin{array}{l}
\mathbb{E}(|\xi^n_{i,j}|^v\mid\mathcal F^n_i) \leq K_{m,v}\Delta
_n^{v/2} \bigl(\Delta_n^{v/2}+
\Gamma'(b,m+1)^n_i+\Gamma'(\widetilde{\sigma},m+1)^n_i \bigr),\\
\mathbb{E}(|\xi'^n_{i,j}|^v\mid\mathcal F^n_i) \leq
K_{m,v}\Delta_n^{v/4+((1/2)\wedge(v/4))},\\
\mathbb{E}(|\overline{\lambda}{}^n_{i,j}|^v\mid\mathcal F^n_i) \leq
K_{m,v}\Delta_n^{v/4+((1/2)\wedge(v/4))}.
\end{array}
\right\}\hspace*{-32pt}
\end{equation}

(4) Next we prove that, for $u$, an odd integer,
%
%
\begin{equation}\label{CC37}
\mathbb{E} ((\overline{W}{}^n_{i+j})^u \xi'^n_{i,j}\mid\mathcal
F^n_i ) = 0.
\end{equation}
We prove this separately for each of the three terms constituting
$\xi'^n_{i,j}$. Since $x\mapsto x^u$ is an odd function, this is
obvious for the term involving $b^n_i$ and also for the
term involving $\widetilde{\sigma}^n_i$. For the term involving $M$,
we have $(\overline{W}{}^n_{i+j})^u=
Y+\int_S^{T+u_n}\rho_s\,dW_s$ for some $\mathcal F_S$-measurable
variable $Y$
and process $\rho$ adapted to the filtration $(\mathcal F^W_t)$
generated by
the Brownian motion. Since this
term is a martingale increment we are left to prove
$\mathbb{E}(U_{T+u_n}\mid\mathcal F_S)=0$ where
\[
U_t = \biggl(\int_S^t\rho_s\,dW_s \biggr) \biggl(\int
_T^{T+u_n}g_n\bigl(s-(i+j)\Delta_n\bigr)
(M_s-M_S)\,dW_s \biggr).
\]
It\^o's formula yields
$U_t=M'_t+\int_T^tg_n(s-(i+j)\Delta_n)
\rho_s(M_s-M_S)\,ds$ for $t\geq T$ where $M'$
is a martingale with $M'_S=0$, so it is enough to prove that
%
%
\begin{equation}\label{CC38}
\mathbb{E}\bigl(\rho_t(M_t-M_S)\mid\mathcal F_S\bigr)=0.
\end{equation}
But for any fixed $t\geq T$ we again have $\rho_t=Y'_t+\int_S^t
\rho'_s\,dW_s$ where $Y'_t$ is $\mathcal F_S$-measurable. Hence (\ref{CC38})
follows from the orthogonality of $W$ and $M$, and we have (\ref{CC37}).

(5) Now we\vspace*{2pt} use (\ref{L4N2}), (\ref{L4N6}) and (\ref{CC36}),
and the form of $G^{u,n}_{i,j}$ as a product of three terms at the
respective powers $l$, $v=p-2l-u$ and $u$. Then
H\"older's inequality with the respective exponents $l'=2p/l$ and
$v'=4p/(p-2l-u)$ [so $2ll'=vv'=4p$ and (\ref{L4N2}) and (\ref{L4N6})
apply] and $u'=4p/(3p+u)$ yields $\mathbb{E}(|G^{u,n}_{i,j}|\mid
\mathcal F^n_i)\leq
K\Delta_n^{p/4+((u/4)\wedge(1/2u'))}$. Observing\vspace*{-2pt} that
$(u/4)\wedge(1/2u')>1/4$ when $u\geq2$, we deduce that
(\ref{CC352}) holds for $G^{u,n}_{i,j}$ when $u\geq2$. It remains
to study~$G^{1,n}_{i,j}$, which is the sum $G^n_{i,j}+G'^n_{i,j}$, where
\[
G^n_{i,j} = (\widehat{\chi}^n_{i+j})^l (\overline{\kappa
}^n_{i,j})^{p-2l-1} \xi
^n_{i,j},\qquad
G'^n_{i,j} = (\widehat{\chi}^n_{i+j})^l (\overline{\kappa
}^n_{i,j})^{p-2l-1} \xi'^n_{i,j}.
\]

By (\ref{L4N2}), (\ref{L4N6}) and (\ref{CC36}), and by
H\"older's inequality as above, we get
\[
\mathbb{E}(|G^n_{i,j}|\mid\mathcal F^n_i) \leq K\Delta
_n^{p/4+1/4} \bigl(
\sqrt{\Delta_n}+\sqrt{\Gamma'(b,m+1)^n_i+\Gamma'(\widetilde
{\sigma},m+1)^n_i} \bigr).
\]
Then by Lemma \ref{LLNN1} we deduce that $G^n_{i,j}$ satisfies
(\ref{CC352}).

(6) It remains to study $G'^n_{i,j}$, which is also
$G'^n_{i,j}=\sum_{w=0}^{p-2l-1}C_{p-2l-1}^wa(n,w,i,j)$, where
$a(n,w,i,j)=(\sigma^n_i\overline{W}{}^n_{i+j})
^{p-2l-1-w} \xi'^n_{i,j} (\widehat{\chi}^n_{i+j})^l (\overline
{\chi}^n_{i+j})^w$.
By successive conditioning, (\ref{L3N11}) and (\ref{CC36}) yield
$\mathbb{E}(|a(n,w,i,j)|\mid\mathcal F^n_i))\leq K\Delta
_n^{p/4+1/2}$ when $w$ is odd.
When $w$ is even,
the same argument with (\ref{L3N10}), plus (\ref{CC37}) and the
fact that $p-2l-1-w$ is then odd yield
\[
|\mathbb{E}(a(n,w,i,j)\mid\mathcal F^n_i)| = \mathrm{O}_u
\bigl(\Delta_n^{p/4+1/4}
\bigl(\Delta_n^{1/2}+\Gamma(\alpha,m)^n_i \bigr) \bigr)
\]
and by Lemma \ref{LLNN1}, the proof is complete.
\end{pf}
\begin{lem}\label{LCC5} Under (SN-$p$) we have, as $n\to\infty$,
%
%
\begin{equation}\label{CC40}\qquad\quad
\left.
\begin{array}{l}
\displaystyle\frac1{\Delta_n^{1/4}} \biggl(\Delta_n^{1-p/4}D(g,m)^n_t-\frac m{m+1}
m_p (\theta \overline{g}(2))^{p/2}\int_0^t|\sigma_s|^p \,ds
\biggr) \stackrel{\mathit{u.c.p.}}{\longrightarrow} 0,\\[12pt]
\displaystyle\frac1{\Delta_n^{1/4}} \biggl(\Delta_n^{1-p/4}\overline
{D}(g,m)^n_t-\frac1{m+1}
m_p (\theta\overline{g}(2))^{p/2}\int_0^t|\sigma_s|^p \,ds
\biggr) \stackrel{\mathit{u.c.p.}}{\longrightarrow} 0.
\end{array}
\right\}
\end{equation}
\end{lem}
\begin{pf}
By (\ref{RR2}), $\overline{W}(g)^n_{i+j}$ is
independent of $\mathcal F^n_i$, and $\mathcal N(0,\overline
{g}(2)_n\Delta_n)$.
So by virtue of
(\ref{SET3}) and (\ref{SET4}) we have $\mathbb{E}((\overline
{W}(g)^n_{i+j})^p\mid\mathcal F^n_i)
=m_p(\theta\overline{g}(2))^{p/2}\Delta_n^{p/4}+\mathrm
{O}_u(\Delta
_n^{p/4+1/2})$. Therefore
by (\ref{L3N211}), the left-hand side of the first expression in
(\ref{CC40}) is smaller in absolute value than
\begin{eqnarray*}
&&\frac K{\Delta_n^{1/4}} \Biggl|
(m+1)k_n\Delta_n\sum_{i=1}^{i_n(m,t)}\bigl|\sigma_{t(m,n,i)}\bigr|^p-
\int_0^t|\sigma_s|^p \,ds \Biggr|
+Kt(m+1)\Delta_n^{1/4}\\
&&\qquad{} +Kt(m+1)\sqrt{\Delta_n}\sum_{i=1}^{i_n(m,t)}
\bigl(\Gamma'(\alpha,m)^n_{I(m,n,i)}+\Gamma'(\beta
(3),m)^n_{I(m,n,i)} \bigr).
\end{eqnarray*}
The second term above goes to $0$, as the last term (locally uniformly
in $t$, in probability) by Lemma \ref{LLNN1}. The first term
goes to $0$ locally uniformly in $t$ in probability as well because
of our Hypothesis \ref{hypoK} (see, e.g., \cite{J2}). Therefore the first
assertion in (\ref{CC40}) holds, and the second one is proved in
the same way.
\end{pf}
\begin{lem}\label{LCC7} Under (SN-$2p$) we have for all $m\geq2$
and $t>0$,
\[
\mathbb{E} \Bigl(\sup_{s\leq t} (\Delta_n^{3/4-p/4}\overline
{N}(g,m)^n_s
)^2 \Bigr)
\leq\frac{Kt}m.
\]
\end{lem}
\begin{pf}
By Doob's inequality, the
left-hand side above is smaller than
\[
4\Delta_n^{3/2-p/2} \sum_{i=1}^{i_n(m,t)}\mathbb{E} ((\overline
{\delta}(g,m)^n_i)^2 ),
\]
whereas (\ref{L4N8}) yields
$\mathbb{E} ((\overline{\delta}(g,m)^n_i)^2 )\leq K\Delta
_n^{p/2-1}$. Since
$i_n(m,t)\leq Kt/m\sqrt{\Delta_n}$, we readily deduce the result.
\end{pf}

\subsection{An auxiliary CLT}\label{ssec:ACLT}

From what precedes, the leading processes for the behavior of
$\overline{V}(g,p)^n$ are the processes $N(g,m)^n$, and here
we prove a CLT for the vector $(N(g_i,m)^n)_{1\leq i\leq d}$ when
$m\geq2$ is fixed and $(g_i)_{1\leq i\leq d}$ is a family of functions
satisfying (\ref{1}).
We first complement the notation (\ref{L3N61}). For $\zeta,\eta\in
\mathbb{R}$
and $p>0$ and $m\geq1$ we set
%
%
\begin{equation}\label{CC70}
\left.
\begin{array}{l}
\displaystyle\mu_{2p}^m(g,h;\eta,\zeta) = \sum_{r,r'=0}^{p/2}\rho_{p,r}\rho_{p,r'}
(2\zeta^2\overline{g}{}'(2))^r(2\zeta^2\overline{h}{}'(2))^{r'},\\
\displaystyle\int_0^m\int_0^m\mathbb{E}' \bigl(\bigl(\eta L(g)_s+\zeta L'(g)_s\bigr)^{p-2r}
\bigl(\eta L(h)_t+\zeta L'(h)_t\bigr)\bigr)^{p-2r'} \,ds \,dt,\\
\displaystyle\overline{\mu}^m_{2p}(g,h;\eta,\zeta) = \frac1{m+1}\bigl(
\mu^m_{2p}(g,h;\eta,\zeta)-m^2 \mu_p(g;\eta,\zeta) \mu
_p(h;\eta
,\zeta) \bigr).
\end{array}
\right\}\hspace*{-36pt}
\end{equation}
Exactly as in Lemma \ref{LL3N11}, the matrix with entries
$\overline{\mu}^m_{2p}(g_i,g_j;\eta,\zeta)$ is symmetric nonnegative.
\begin{prop}\label{PCC1} Assume (SN-$4p$), and let $m\geq2$.
The sequence of $d$-dimensional processes with components
$\Delta_n^{3/4-p/4}N(g_i,m)$ converges stably in law to a process of the
following form:
%
%
\begin{equation}\label{CC41}
\Biggl(\theta^{1/2-p/2}\sum_{j=1}^d\int_0^t\psi_{ij}^m(\theta
\sigma
_s,\alpha_s)\, dB^j_s
\Biggr)_{1\leq i\leq d},
\end{equation}
where $B$ is as in Theorem \ref{TCC1} and $\psi^m$ is a measurable
$d\times d$ matrix-valued function such that
$(\psi^m\psi^{m\star})(\eta,\zeta)$ is the matrix with entries
$\overline{\mu}_{2p}^m(g_i,g_j;\eta,\zeta)$, as defined by (\ref{CC70}).
\end{prop}

We begin with a lemma, for which we use the notation $\Gamma$ of Lemma
\ref{LNN10}.
\begin{lem}\label{LCC1} Let $m\geq2$ and $s\in\Gamma$ and
$i_n=\min(i\dvtx
I(m,n,i)\Delta_n\geq s)$. Then under (SN-$4p$)
we have the following almost sure convergences:
%
%
\begin{eqnarray}\label{CC85}
&&\Delta_n^{1/2-p/4}\gamma(g,m)^n_{i_n} \nonumber\\[-8pt]\\[-8pt]
&&\qquad\to m m_p \theta^{1+p/2}
\overline{g}(2)^{p/2}
|\sigma_s|^p = m\theta^{1-p/2}\mu_p(g;\theta\sigma_s,\alpha_s),\nonumber
\\
\label{CC86}
&&\Delta_n^{1-p/2}\mathbb{E} (\delta(g,m)^n_{i_n}\delta(h,m)^n_{i_n}
\mid\mathcal F(m)^n_{i_n-1} ) \nonumber\\[-8pt]\\[-8pt]
&&\qquad\to \theta^{2-p} \mu
^m_{2p}(g,h;\theta\sigma
_s,\alpha_s).\nonumber
\end{eqnarray}
\end{lem}
\begin{pf}
We set $i'_n=I(m,n,i_n)$ and $s_n=i'_n\Delta_n$,
which converges to $s$. Both results are consequences of Lemma
\ref{LNN10}: first, by (\ref{L3N211}) with $u=1$, (\ref{CC85})
follows from
%
%
\begin{equation}\label{CC851}\qquad\quad
\Delta_n^{1/2-p/4} \mathbb{E} \Biggl(\sum_{j=0}^{mk_n-1}|\sigma
_{s_n}\overline{W}
(g)^n_{i'_n+j}|^p
\mid\mathcal F_{s_n} \Biggr) \to m m_p \theta^{1+p/2} \overline
{g}(2)^{p/2} |\sigma_s|^p.
\end{equation}
Then we apply Lemma \ref{LNN10} with $d=1$ and $g_1=g$ and with the functions
\[
f_n(x,y,z)=\frac1{k_n}\sum_{j=0}^{mk_n-1}|x(j/k_n)|^p,\qquad
f(x,y,z)=\int_0^m|x(s)|^p \,ds,
\]
which satisfy (\ref{KL4}) and $f_n\to f$ pointwise. The left-hand (right)
side of (\ref{CC851}) is equal to $\Delta_n^{1/2-p/4}/k_n^{p/2-1}$ times
($\theta^{1-p/2}$ times) the left-hand (right) side of
(\ref{KL5}); hence (\ref{CC85}) holds.

For (\ref{CC86}) we apply Lemma \ref{LNN10} with $d=2$ and $g_1=g$
and $g_2=h$ and the functions
\begin{eqnarray*}
&&f_n((x,x'),(y,y'),(z,z'))\\
&&\qquad=\sum_{r,r'=0}^{p/2}\rho_{p,r}\rho_{p,r'}
\frac1{k_n^2}\sum_{j,j'=0}^{mk_n-1}
\biggl(x \biggl(\frac j{k_n} \biggr)+y \biggl(\frac j{k_n} \biggr) \biggr)^{p-2r}\\
&&\qquad\quad\hspace*{101.4pt}{}\times
\biggl(x' \biggl(\frac{j'}{k_n} \biggr)+y' \biggl(\frac{j'}{k_n} \biggr) \biggr)^{p-2r'}
z \biggl(\frac j{k_n} \biggr)^r z' \biggl(\frac{j'}{k_n} \biggr)^{r'},\\
&&
f((x,x'),(y,y'),(z,z'))\\
&&\qquad=\sum_{r,r'=0}^{p/2}\rho_{p,r}\rho_{p,r'}
\int_0^m\int_0^m\bigl(x(s)+y(s)\bigr)^{p-2r}
\bigl(x'(t)+y'(t)\bigr)^{p-2r'}\\
&&\qquad\quad\hspace*{97.6pt}{}\times z(s)^r z'(t)^{r'} \,ds \,dt,
\end{eqnarray*}
which satisfy (\ref{KL4}) and $f_n\to f$ pointwise. The left-hand (right)
side of (\ref{CC86}) is equal to $\Delta_n^{1-p/2}/k_n^{p-2}$ times
($\theta^{2-p}$ times) the left-hand (right) side of
(\ref{KL5}); hence (\ref{CC86}) holds.
\end{pf}
\begin{pf*}{Proof of Proposition \protect\ref{PCC1}}
(1) As is well known,
and with
the $d$-dimensional variables with components
$\xi^{n,k}_i=\Delta_n^{1/2-p/4} (\delta(g_k,m)^n_i-\gamma
(g_k,m)^n_i)$ (which
are martingale differences),
it suffices to prove the following three convergences, for all $t>0$ and
all bounded martingales $N$:
%
%
\begin{eqnarray}\label{CC71}\qquad\quad
&\displaystyle\sqrt{\Delta_n} \sum_{i=1}^{i_n(m,t)}\mathbb{E}(\xi^{n,k}_i \xi
^{n,l}_i\mid\mathcal F(m)^n_{i-1}) \stackrel{\mathbb
{P}}{\longrightarrow}
\theta^{1-p}\int_0^t\overline{\mu}_{2p}^m(g_k,g_l;\theta\sigma
_s,\alpha_s) \,ds,&
\\
%
%
\label{CC72}
&\displaystyle\Delta_n \sum_{i=1}^{i_n(m,t)}\mathbb{E}(\|\xi^n_i\|^4\mid\mathcal
F(m)^n_{i-1}) \stackrel{\mathbb{P}}{\longrightarrow} 0,&
%
%
\\
\label{CC73}
&\displaystyle\Delta_n^{1/4} \sum_{i=1}^{i_n(m,t)}\mathbb{E}\bigl(\xi
^n_i \bigl(N_{i(m+1)u_n}-N_{(i-1)(m+1)u_n}\bigr)
\mid\mathcal F(m)^n_{i-1}\bigr) \stackrel{\mathbb{P}}{\longrightarrow} 0&
\end{eqnarray}
(we use Theorem IX.7.28 of \cite{JS}, with $Z$ being a
bounded martingale of the form $Z_t=\int_0^tu_s\,dW_s$ for some
predictable process $u$ with values in $(0,1]$).

(2) Equation (\ref{L4N8}) and H\"older's inequality imply
$\mathbb{E}(|\delta(g_k,m)^n_i|^4\mid\mathcal F(m)^n_{i-1})\leq
K_m\Delta_n^{p-2}$.
Then the
expected value of the left-hand side of (\ref{CC72}) is smaller than
$K_m\sqrt{\Delta_n}$, yielding (\ref{CC72}).
The proof of (\ref{CC71}) is similar to the proof of Theorem
\ref{TLNN1}. Set $\zeta^n_i=\mathbb{E} (\xi^{n,k}_i
\xi^{n,l}_i\mid\mathcal F(m)^n_{i-1} )$, and $\gamma_s=
\overline{\mu}_{2p}^m(g_k,g_l;\theta\sigma_s,\alpha_s)$. Since
$k_n\sqrt{\Delta_n}\to\theta$, we need to show that
\[
(m+1)k_n\Delta_n \sum_{i=1}^{i_n(m,t)}\zeta^n_i \stackrel{\mathbb
{P}}{\longrightarrow} (m+1)\theta^{2-p}\int_0^t\gamma_s \,ds.
\]
Note that $|\zeta^n_i|\leq K_m$. Then, as for Theorem
\ref{TLNN1}, the above will follow from the fact that for any
$s\in\Gamma$, and with the notation $i_n$ of
Lemma \ref{LCC1}, we have [similar to (\ref{LLN59})]
%
%
\begin{equation}\label{CC740}
\zeta^n_{i_n} \to\theta^{2-p}\gamma_s \qquad\mbox{a.s.}
\end{equation}
Then (\ref{CC740}) readily follows from Lemma \ref{LCC1} and (\ref{CC70}),
once observed that
\[
\zeta^n_{i_n} = \Delta_n^{1-p/2} \bigl(\mathbb{E} (\delta(g_k,m)^n_{i_n}
\delta(g_l,m)^n_{i_n}\mid\mathcal F(m)^n_{i_n-1} )-\gamma(g_k,m)^n_{i_n}
\gamma(g_l,m)^n_{i_n} \bigr).
\]

(3) Now we turn to (\ref{CC73}), which we prove for the first
component, say with \mbox{$g=g_1$}. For simplicity we write
$D^n_i(Y)=Y_{i(m+2)u_n}-Y_{(i-1)(m+2)u_n}$ for any process $Y$. In
view of the definition of $\xi^n_i$, and since $N$ is a martingale, it
is enough to prove that
\[
\Delta_n^{3/4-p/4} \sum_{i=1}^{i_n(m,t)}\mathbb{E}(\delta
(g,m)^n_iD^n_i(N)\mid\mathcal F
(m)^n_{i-1}) \stackrel{\mathbb{P}}{\longrightarrow} 0.
\]

Observe that (\ref{AU2}) and (\ref{L3N211}) yield $\delta
(g,m)^n_i=\delta'^n_i+\Psi'^n_i$ where
\[
\delta'^n_i = \sum_{j=0}^{mk_n-1}\bigl(\sigma^n_{I(m,n,i)}\overline{W}
(g)^n_{I(m,n,i)+j}\bigr)^p,\qquad
\mathbb{E}(|\Psi'^n_i|^2\mid\mathcal F(m)^n_{i-1}) \leq
K\Delta_n^{p/2-1/2}.
\]
Since $N$ is bounded, $\sum_{i=1}^{i_n(m,t)}\mathbb
{E}((D^n_i(N))^2)\leq K$ and the Cauchy--Schwarz
inequality yields
\begin{eqnarray*}
\mathbb{E} \Biggl(\Delta_n^{3/4-p/4}\sum_{i=1}^{i_n(m,t)}|\Psi'^n_i|
|D^n_i(N)| \Biggr) &\leq&
K\Delta_n^{3/4-p/4} \Biggl(\mathbb{E} \Biggl(\sum_{i=1}^{i_n(m,t)}(\Psi
'^n_i)^2 \Biggr)
\Biggr)^{1/2}
\\
&\leq& K\Delta_n^{1/4}.
\end{eqnarray*}
Therefore it remains to prove that
%
%
\begin{equation}\label{CC76}
\Delta_n^{3/4-p/4} \sum_{i=1}^{i_n(m,t)}\mathbb{E}(\delta
'^n_iD^n_i(N)\mid\mathcal F
(m)^n_{i-1}) \stackrel{\mathbb{P}}{\longrightarrow} 0.
\end{equation}

(4) Observe that, by $\mathbb{E}(|\delta'^n_i|^2)\leq K\Delta
_n^{p/2-1}$ and the
Cauchy--Schwarz
inequality,
\[
\Delta_n^{3/4-p/4} \sum_{i=1}^{i_n(m,t)}\mathbb{E} (|\delta
'^n_iD^n_i(N)| ) \leq
K\sqrt{\mathbb{E}(N^2_t)}.
\]
Then the set of all square-integrable martingales $N$ satisfying
(\ref{CC76}) is closed under $\mathbb{L}^2$-convergence, and thus for proving
(\ref{CC76}) we can use the following scheme:

(a) Prove (\ref{CC76}) when $N$ is $(\mathcal F^{(0)}_t)$-adapted and
orthogonal to $W$.

(b) Prove (\ref{CC76}) when $N_t=\int_0^t\gamma_s\,dW_s$ where
$\gamma$ is $(\mathcal F^{(0)}_t)$-adapted and constant in time over
intervals $(t_{i-1},t_i]$ with $t_0=0$ and $t_q=\infty$ for some
$q$.

(c) Conclude from the closeness proved before that (\ref{CC76}) holds
for all $N\in\mathcal N^0$, the set of all bounded
$(\mathcal F^{(0)}_t)$-martingales.

(d) Prove (\ref{CC76}) when $N$ is in the set $\mathcal N^1$ of all
martingales
having $N_\infty=f(\chi_{t_1},\ldots,\chi_{t_q})$ where $f$ is any
Borel bounded on $\mathbb{R}^q$ and $t_1<\cdots<t_q$ and $q\geq1$.

(e) Since $\mathcal N^0\cup\mathcal N^1$ is a total subset of the
set of
all square-integrable $(\mathcal F_t)$-martingales, conclude once more
from the
closeness that (\ref{CC76}) holds for all such~$N$.

We are thus left to prove (a), (b) and (d), and for these we can
reproduce Step 5 of the proof of Lemma 5.7 in \cite{JLMPV}.
\end{pf*}

\subsection[Proof of Theorem 4.1]{Proof of Theorem \protect\ref{TCC1}}\label{ssec:41}

By localization we may assume (SN-$4p$) and Hypothesis \ref{hypoSK}. Then, upon applying
Lemmas \ref{LCC2}, \ref{LCC3}, \ref{LCC5} and \ref{LCC7}, we readily
deduce (\ref{AXX}) from (\ref{CC1}) and (\ref{CC15}).

On the other hand, we fix the $d$-dimensional
Brownian motion $B$ in (\ref{CC41}) and (\ref{CC5}) (the same in
both). Proposition \ref{PCC1} yields, for each fixed $m$, that
$\Delta_n^{3/4-p/4} N(g,m)$ stably converges in law to the right-hand
side of (\ref{CC41}). Next, the following property is implicitly proved
in the proof of Lemma \ref{LL3N11} in Section~\ref{ssec-KL}
(with $T$ playing the role of $m$ here):
\[
\overline{\mu}^m_{2p}(g,h;\eta,\zeta) \to\overline{\mu
}_{2p}(g,h;\eta,\zeta)\qquad
\mbox{as } m\to\infty.
\]
Then we see that we can choose suitable versions for the square-roots
$\psi$ and $\psi^m$ in such a way that $\psi^m(\eta,\zeta)
\to\psi(\eta,\zeta)$ for all $\eta ,\zeta $. Then
(\ref{CC41}) converges in probability toward (\ref{CC5}). The result
then follows from (\ref{AXX}) in a standard way.

\subsection[Theorem 4.4: A key decomposition]{Theorem \protect\ref{TCD1}: A key decomposition}\label{ssec:43}

Here we start the proof of Theorem~\ref{TCD1}, by providing a
decomposition for the processes $\widetilde{V}{}^*(g,p)^n$ of
(\ref{CD1}). So we fix $p>3$, and assume $\alpha$ c\`adl\`ag.
By localization we can and will assume (SN-$2p$) and Hypothesis \ref{hypoSK} without
special mention.

The choice of the exhausting sequence $(T_m)$ in (\ref{CD6}) is arbitrary,
but a convenient choice is as
follows: for $q\geq1$ we consider the successive
jump times $(T(q,m)\dvtx m\geq1)$ of the Poisson process
$\underline{\mu}((0,t]\times\{z\dvtx1/q<\gamma(z)\leq1/(q-1)\})$ where
$\gamma$
is the
function occurring in Hypothesis \ref{hypoSH}. Those stopping times have
pairwise disjoint graphs as $m$ and $q$ vary, and $(T_m)_{m\geq1}$
denotes any reordering of the double sequence $(T(q,m)\dvtx q,m\geq1)$. We
complete this sequence by setting $T_0=0$.

Let $P_q$ be the set of all $m\geq1$ such that $T_m=T(q',m')$ for some
$m'\geq1$ and some $q'\leq q$. We consider
the following processes [compare with (\ref{LLN2})]:
%
%
\begin{equation}\label{CD307}\qquad
\left.
\begin{array}{l}
X^q = \bigl(\delta1_{\{z\dvtx\gamma(z)>1/q\}}\bigr)*\underline{\mu},
\qquad
M^q = \bigl(\delta1_{\{z\dvtx\gamma(z)\leq1/q\}}\bigr)*(\underline{\mu
}-\underline{\nu}),\\[2pt]
X'^q = X-X^q,\qquad X''^q = X'^q-M^q,\\
Z'^q = X'^q+\chi,\qquad Z''^q = X''^q+\chi.
\end{array}
\right\}
\end{equation}
So $X''^q$ satisfies (\ref{CIS}) with the same $\sigma$ as
in (\ref{CD300}) and a bounded drift given by
%
%
\begin{equation}\label{CD3079}
b^q_t = b_t-\int_{\{z\dvtx\gamma(z)>1/q,|\delta(t,z)|\leq1\}}\delta
(t,z) \lambda(dz).
\end{equation}
Here, $X^q$ is the sum of ``big'' jumps, and this is the part of
$X$ which essentially imports for our CLT: more precisely, we
single out the summands in $V(X,g,p,0)^n_t$ which involve at least
one jump of $X^q$ [after centering this is the process $Y(q,g)$
defined below in (\ref{CD123})]. We obtain a CLT for these
processes $Y(q,g)$ in a relatively simple way, and then prove that
the contribution of the other summands is negligible, when $n$
and $q$ are large (so the cut-off level $1/q$ for the ``big'' jumps
is small).

We denote by $\Omega_n(t,q)$ the set of all $\omega$ such that for
any $m,m'
\in P_q$ with $T_m(\omega)\leq t$, we have $2u_n<T_m(\omega)\leq
t-4u_n$, and
$|T_m(\omega)-T_{m'}(\omega)|>4u_n$, and also $T_m(\omega)/\Delta_n$
is not an
integer. Since the set $\{T_m\dvtx m\in P_q\}$
is locally finite and $\mathbb{P}(T_m=t)=0$ for all $m$ and $t\geq0$,
we have
%
%
\begin{equation}\label{CD121}
\Omega_n(t,q) \to\Omega\qquad\mbox{a.s., as } n\to\infty.
\end{equation}

Next, we\vspace*{1pt} denote by $\widetilde{V}{}^*(X,g,p)^n$ the process defined by
(\ref{CD1}) to emphasize the dependency on $X$, and likewise
we have $\widetilde{V}{}^*(X'^q,g,p)^n$.
Then a (relatively) simple computation shows the
following key property which holds on the set $\Omega_n(t,q)$:
%
%
\begin{eqnarray}\label{CD123}
\widetilde{V}{}^*(X,g,p)^n_t &=& \widetilde
{V}{}^*(X'^q,g,p)^n_t+Y(q,g)^n_t,
Y(q,g)^n_t \nonumber\\[-8pt]\\[-8pt]
&=& \sum_{m\in P_q\dvtx T_m\leq t}\zeta(q,g)^n_m,\nonumber
\end{eqnarray}
where, with the random integer $I^n_m=[T_m/\Delta_n]$, we have set
\begin{eqnarray*}
\zeta(q,g)^n_m &=& \frac1{\Delta_n^{1/4}k_n} \Biggl(\sum_{j=1}^{k_n-1}
\bigl( |\overline{Z'^q}(g)^n_{I^n_m+1-j}+g^n_j\Delta X_{T_m} |^p
\\
&&\hspace*{62.66pt}{} - |\overline{Z'^q}(g)^n_{I^n_m+1-j} |^p
- |g^n_j\Delta X_{T_m} |^p \bigr)
\\
&&\hspace*{85pt}{}
+\bigl(\overline
{g}(p)_n-k_n\overline
{g}(p)\bigr)|\Delta
X_{T_m}|^p \Biggr).
\end{eqnarray*}
[Note that $\overline{Z'^q}(g)^n_{I^n_m-j}$ possibly involves
$\Delta^n_lZ'(q)$ for negative integers $l$, although this does not occur
on the set $\Omega_n(t,q)$ when $m\in P_q$ and $j\leq2k_n$; however,
to have
such variables defined everywhere, we make the convention $\Delta^n_iY=0$
for any process $Y$ when $i\leq0$.]

\subsection{The processes $Y(q,g)^n$}\label{ssec:YQ}

The aim of this subsection is to prove the following proposition.
\begin{prop}\label{P37} If $q\geq1$ and $t\geq0$ are fixed, and in the
same setting as before, we have (with $\stackrel{\mathcal
L-(s)}{\longrightarrow}$ denoting the
stable convergence in law)
%
%
\begin{equation}\label{CD135}
(Y(q,g_l)^n_t)_{1\leq l\leq d} \stackrel{\mathcal
L-(s)}{\longrightarrow} U(p,q)_t,
\end{equation}
where $U(p,q)$ is the $d$-dimensional process associated with the
functions $(g_l)$ by (\ref{CD6}), except that the sums are taken over
$m\in P_q$ only.
\end{prop}

We start with the following lemma which describes the behavior of the
variables (with $m\in P_q$):
%
%
\begin{equation}\label{CD11}
\eta(q,g)^n_m = \frac1{\Delta_n^{1/4}k_n} \sum_{j=1}^{k_n-1}
\{g^n_j\}^{p-1} \overline{Z'^q}(g) ^n_{I^n_m-j+1}.
\end{equation}
The two key properties for the next lemma are the approximation
(\ref{AU1}) and the fact that the times $T_m$ are independent of $W$
and with an absolutely continuous law. Recall that we
have the family $(g_l)_{1\leq l\leq d}$ of $d$ weight
functions with $(U_{m-},U_{m+},\overline{U}_{m-},\overline{U}_{m+})$
associated as before
Theorem \ref{TCD1}.
\begin{lem}\label{LCD1} For any $q\geq1$, the $(\mathbb
{R}^d)^{\mathbb{N}^\star}$-valued
variables $(\eta(q$,\break $g_l)^n_m)_{1\leq l\leq d,m\in P_q}$
converge stably in law, as $n\to\infty$, to
$(\eta_m)_{m\in P_q}$, where $\eta_m$\break is the $d$-dimensional
variable given by
%
%
\begin{equation}\label{CD9}
\eta_m = \sqrt{\theta} \sigma_{T_m-}U_{m-}
+\sqrt{\theta} \sigma_{T_m}U_{m+}+\frac{\alpha_{T_m-}}{\sqrt
{\theta}}
\overline{U}_{m-}+\frac{\alpha_{T_m}}{\sqrt{\theta}} \overline{U}_{m+}.
\end{equation}
\end{lem}
\begin{pf}
As is well known, it is enough to prove the result for
any finite subset of $m$'s, say in a finite subset $P'_q$ of $P_q$. Since
$q$ is fixed, we drop it from the notation, writing $Z'=Z'^q$,
$Z''=Z''^q$, $M=M^q$, $X'=X'^q$ and $X''=X''^q$.

(1) The times $(T_m\dvtx m\in P_q)$ are independent of
$W$, and also of the restriction $\underline{\mu}(q)$ of the Poisson measure
$\underline{\mu}$ to the set $\mathbb{R}_+\times\{z\dvtx\gamma(z)\leq
1/q\}$. Hence if
$\mathcal H_t=\mathcal F_t\vee\sigma(T_m\dvtx m\in P_q)$, the process $W$
is a
Brownian motion and the measure $\underline{\mu}(q)$ a Poisson
measure with
compensator $\underline{\nu}(q)$, the restriction of $\underline{\nu
}$ to
$\mathbb{R}_+\times\{z\dvtx\gamma(z)\leq1/q\}$ again, relative
to the filtration $(\mathcal H_t)$. Thus $X''$ admits the same representation
(\ref{CIS}) and $M$ has the same form (\ref{CD307}) relative to the two
filtrations $(\mathcal F_t)$ and~$(\mathcal H_t)$. With the random
integers $I_m^n$
being $\mathcal H_0$-measurable, we deduce from (\ref{LLN141}) and
(\ref{L4N2})
and (\ref{LLN14}) for $M$ together with $|(g_l)_n(s)|\leq K$
that, for $v\in(0,2p]$ and $j\in\mathbb{Z}$ and $i=l,\ldots,d$,
%
%
\begin{equation}\label{CD320}\quad
\left.
\begin{array}{l}
\mathbb{E}(|\Delta^n_{I^m_n+j}X''|^v) \leq K_{v,q}\Delta
_n^{v/2},\\[4pt]
\mathbb{E}(|\overline{M}(g_l)^n_{I^n_m+j}|^2) \leq K\sqrt{\Delta
_n},\\[4pt]
\mathbb{E} \bigl(|\overline{X}{}''(g_l)^n_{I^n_m+j}|^v+
|\overline{Z}{}''(g_l)^n_{I^n_m+j}|^v \bigr)
\leq K_{v,q}\Delta_n^{v/4}.
\end{array}
\right\}
\end{equation}

Now if $f$ is a bounded function on $\mathbb{R}$, arguments similar to
the one
giving (\ref{L4N4}) [relative to the filtration $(\mathcal H_t)$ and using
that $\sigma$ is c\`adl\`ag and bounded and the drift $b^q$ is also
bounded], we obtain that if $2\leq k'_n\leq2k_n$,
%
%
\begin{equation}\label{CD319}
\left.
\begin{array}{l}
\displaystyle\mathbb{E} \Biggl( \Biggl|\sum_{j=0}^{k'_n}f(j/k_n)\Delta
^n_{I^m_n-j}X''-\sigma_{T_m-}
\sum_{j=0}^{k'_n}f(j/k_n)\Delta^n_{I^m_n-j}W \Biggr|^v \Biggr)\\
\qquad =
\mathrm{o}_u(\Delta_n^{v/4}),\\
\displaystyle\mathbb{E} \Biggl( \Biggl|\sum_{j=2}^{k'_n}f(j/k_n)\Delta
^n_{I^m_n+j}X''-\sigma_{T_m}
\sum_{j=2}^{k'_n}f(j/k_n)\Delta^n_{I^m_n+j}W \Biggr|^v \Biggr)\\
\qquad =
\mathrm{o}_u(\Delta_n^{v/4}).
\end{array}
\right\}
\end{equation}
Moreover $A_n=\sum_{j=0}^{k'_n}f(j/k_n)\Delta^n_{I^m_n-j}M$, say, can be
written as $\delta^n\star(\underline{\mu}(q)-\underline{\nu}(q))_{T_m}-
\delta^n\star(\underline{\mu}(q)-\underline{\nu}(q))_{T_m-2u_n}$
for some predictable function
$\delta^n$ satisfying $|\delta^n(t,z)|\leq K\gamma(z)$. Then a well-known
result (see, e.g., Lemma 5.12 of \cite{J2}, used with $2u_n$ instead of
$\Delta_n$ and $\eta=\sqrt{u_n}$, and relative to the
filtration\vspace*{-2pt}
$(\mathcal H_t)$) says that $A_n/\sqrt{u_n}\stackrel{\mathbb
{P}}{\longrightarrow}0$. The same holds if we take
the indices $I^n_l+j$ instead of $I^n_m-j$, and thus
%
%
\begin{eqnarray}\label{CD3190}
\frac1{\Delta_n^{1/4}} \sum_{j=0}^{k'_n}f(j/k_n)\Delta
^n_{I^m_n-j}M &\stackrel{\mathbb{P}}{\longrightarrow}&
0,\nonumber\\[-8pt]\\[-8pt]
\frac1{\Delta_n^{1/4}} \sum_{j=2}^{k'_n}f(j/k_n)\Delta
^n_{I^m_n+j}M &\stackrel{\mathbb{P}}{\longrightarrow}& 0.\nonumber
\end{eqnarray}

(2) We put for $i\geq0$ and any weight function $g$,
\begin{eqnarray*}
G(g)^n_{i-}&=&\frac1{k_n} \sum_{j=i+2}^{k_n-1} \{g^n_j \}
^{p-1}g^n_{j-i-1},\qquad
G(g)^n_{0}=\frac1{k_n}
\sum_{j=1}^{k_n-1} \{g^n_j \}^{p-1}g^n_j,\\
G(g)^n_{i+} &=& \frac1{k_n}
\sum_{j=1}^{k_n-i} \{g^n_j \}^{p-1}g^n_{j+i-1},\\
\overline{G}(g)^n_{i-} &=& \sum_{j=i+1}^{k_n-1} \{g^n_j \}
^{p-1}g'^n_{j-i},\qquad \overline{G}(g)^n_{i+}
=\sum_{j=1}^{k_n-i} \{g^n_j \}^{p-1}g'^n_{j+i}.
\end{eqnarray*}
Then a (tedious) computation shows that
\begin{eqnarray*}
\eta(q,g)^n_m&=&\frac1{\Delta_n^{1/4}} \Biggl(
\sum_{i=0}^{k_n-3}G(g)^n_{i-} \Delta^n_{I^n_m-i}X'\\
&&\hspace*{28.8pt}{}
+G(g)^n_{0} \Delta^n_{I^n_m+1}X'
+\sum_{i=2}^{k_n-1}G(g)^n_{i+} \Delta^n_{I^n_m+i}X' \Biggr)\\
&&{} - \frac1{\Delta_n^{1/4}k_n} \Biggl(
\sum_{i=0}^{k_n-2}\overline{G}(g)^n_{i-}\chi^n_{I^n_m-i}
+\sum_{i=1}^{k_n-1}\overline{G}(g)^n_{i+}\chi^n_{I^n_m+i} \Biggr).
\end{eqnarray*}
Moreover, if
\begin{eqnarray*}
H_-(g,t) &=& \int_t^1\{g(s)\}^{p-1}g(s-t)\,ds, \\
H_+(g,t) &=& \int_0^{1-t}\{g(s)\}^{p-1}g(s+t)\,ds,\\
\overline{H}_-(g,t) &=& \int_t^1\{g(s)\}^{p-1}g'(s-t)\,ds, \\
\overline{H}_+(g,t) &=& \int_0^{1-t}\{g(s)\}^{p-1}g'(s+t)\,ds,
\end{eqnarray*}
we have
\[
G(g)^n_{i\pm}=H_\pm\biggl(\frac i{k_n},g \biggr)+\mathrm{O}_u\bigl(\sqrt
{\Delta_n}\bigr),\qquad
\overline{G}(g)^n_{i\pm}=\overline{H}_\pm\biggl(\frac i{k_n},g
\biggr)+\mathrm{O}_u\bigl(\sqrt{\Delta_n}\bigr).
\]

Using $|G^n_{0}|\leq K$ and $X'=M+X''$, (\ref{CD320}),
(\ref{CD319}), (\ref{CD3190}) and (SN-$2p$), we deduce
%
%
\begin{equation}\label{CD323}\qquad\quad
\eta(q,g)^n_m = \sigma_{T_m-} \rho(g)^n_{m-}+\sigma_{T_m} \rho(g)^n_{m+}+
\overline{\rho}(g)^n_{m-}+\overline{\rho}(g)^n_{m+}+\mathrm{o}_{Pu}(1),
\end{equation}
where
\begin{eqnarray*}
\rho(g)^n_{m-}&=&\frac{1}{\Delta_n^{1/4}}
\sum_{i=0}^{k_n-3}H_{-} \biggl(\frac i{k_n},g \biggr) \Delta
^n_{I^n_m-i}W,\\
\rho(g)^n_{m+}&=&\frac{1}{\Delta_n^{1/4}}
\sum_{i=2}^{k_n-1}H_{-} \biggl(\frac i{k_n},g \biggr) \Delta
^n_{I^n_m+i}W,\\
\overline{\rho}(g)^n_{m-}&=&-\frac1{\Delta_n^{1/4}k_n}
\sum_{i=0}^{k_n-2}\overline{H}_{-} \biggl(\frac i{k_n},g \biggr) \chi
^n_{I^n_m-i},\\
\overline{\rho}(g)^n_{m+}&=&-\frac1{\Delta_n^{1/4}k_n}
\sum_{i=1}^{k_n-1}\overline{H}_{+} \biggl(\frac i{k_n},g \biggr) \chi
^n_{I^n_m+i}.
\end{eqnarray*}

(3) At this stage, we use the same ideas as in Lemma \ref{LNN10}.
We denote by $\rho^n_{m\pm}$ and $\overline{\rho}^n_{m\pm}$ the
$d$-dimensional
variables with components $\rho(g_i)^n_{m\pm}$ and $\overline{\rho}
(g_i)^n_{m\pm}$.
First we argue that $\omega^{(0)}\in\Omega^{(0)}$ fixed. Under
$\mathbb{Q}=\mathbb{Q}(\omega^{(0)},\cdot)$ the variables $\overline
{\rho}^n_{m-}$ and $\overline{\rho}^{n}_{m+}$
are independent from each other and
also when $m$ varies in $P'_p$ as soon as $n$ is large
enough [so that $\omega^{(0)}\in\Omega_n(t,q)$]. Moreover, they are
sums, normalized by $1/\Delta_n^{1/4}k_n$, of (approximately) $k_n$
centered independent variables with a bounded fourth moment, and their
covariance matrices are (approximately again) $\alpha_{T_m-}^2/\theta$
and $\alpha_{T_m}^2/\theta$ times Riemann approximations of the
integrals defining $\overline{\Psi}_{p-}$ and $\overline{\Psi
}_{p+}$, respectively. Then
we prove
exactly as for (\ref{KL6}) (only the finite-dimensional convergence
is needed here) that under $\mathbb{Q}$,
%
%
\begin{equation}\label{CD324}\quad
(\overline{\rho}^n_{m-},\overline{\rho}^n_{m+} )
_{m\in P'_q} \stackrel{\mathcal L}{\longrightarrow} \biggl(\frac
{\alpha_{T_m-}(\omega^{(0)})}{\sqrt{\theta
}} \overline{U}_{m-},
\frac{\alpha_{T_m}(\omega^{(0)})}{\sqrt{\theta}} \overline
{U}_{m+} \biggr)_{m\in P'_q}.
\end{equation}
[In fact we prove the convergence in (\ref{CD324}) for each
$m$ first, and then we use the fact that the variables in the left-hand
side are independent for different values of~$m$, under $\mathbb{Q}$,
and as soon as $\omega^{(0)}\in\Omega_n(t,q)$.]

Second, exactly as for (\ref{KL13}) [or as above for (\ref{CD324})],
we get
%
%
\begin{equation}\label{CD326}
(\rho^n_{m-},\rho^n_{m+} )
_{m\in P'_q} \stackrel{\mathcal L}{\longrightarrow} \bigl(\sqrt
{\theta} U_{m-},
\sqrt{\theta} U_{m+} \bigr)_{m\in P'_q}.
\end{equation}
Then $(U_{m\pm},\overline{U}_{m\pm})$ are as described after (\ref{CD4}),
and as in Steps 2 and 4 of the proof of Lemma \ref{LNN10}, we
deduce from the convergences (\ref{CD324}) under $\mathbb{Q}(\omega^{(0)},\cdot)$
and (\ref{CD326}) under $\mathbb{P}^{(0)}$, and from (\ref{CD323}),
that $(\rho^n_{m-},\rho^n_{m+},\overline{\rho}^n_{m-},\overline
{\rho}^n_{m+})_{m\in P'_q}$
converges in law to
$(\sqrt{\theta} U_{m-}$, $\sqrt{\theta} U_{m+},\alpha
_{T_m-}\overline{U}_{m-}/\sqrt
{\theta},
\alpha_{T_m}\overline{U}_{m+}/\sqrt{\theta})_{m\in P'_q}$.
This convergence in law is indeed a stable convergence by the
same argument used to obtain a similar result in \cite{JP}. Finally by
(\ref{CD9}) and (\ref{CD323}) and the definition of the stable convergence
in law, we obtain the claim.
\end{pf}
\begin{pf*}{Proof of Proposition \protect\ref{P37}}
With $\overline{g}_l(p)_n={\sum_{i=1}^{k_n}}|g_l(i/k_n)|^p$, we
have $|\overline{g}_l(p)_n-k_n\overline{g}_l(p)|\leq K$ by
(\ref{SET4}). Then, with the notation (\ref{CD11}), a Taylor expansion
and $|\Delta X_{T_m}|\leq K$ yield
\begin{eqnarray*}
&&|\zeta(q,g_l)^n_m-p\{\Delta
X_{T_m}\}^{p-1}\eta(q,g_l)_m^n |\\
&&\qquad\leq
K\Delta_n^{1/4} \Biggl( 1+\sum_{j=1}^{k_n-1}
\bigl((\overline{Z'^q}(g_l)^n_{I_m^n-j+1})^p
+(\overline{Z'^q}(g_l)^n_{I_m^n-j+1})^2 \bigr) \Biggr).
\end{eqnarray*}
If we apply (\ref{CD320}) we see that the expectation of the sum
in the right-hand side above is bounded (recall $p>3$). Therefore Lemma
\ref{LCD1}
implies
\[
(\zeta(q,g_l)^n_m)_{1\leq l\leq d,m\in P_q} \stackrel{\mathcal
L-(s)}{\longrightarrow}
(p\{\Delta X_{T_m}\}^{p-1}\eta_m )_{m\in P_q}
\]
and (\ref{CD135}) readily follows.
\end{pf*}

\subsection{The processes $\widetilde{V}{}^*(X'^q,g,p)^n$}\label{ssec:YT}

The aim of this subsection is to prove the following proposition.
\begin{prop}\label{PP10} Under the same assumptions as
before, and for all $\varepsilon>0$, we have
\[
\lim_{q\to\infty} \limsup_n \mathbb{P} \bigl(|\widetilde
{V}{}^*(X'^q,g,p)^n_t|>\varepsilon
\bigr) = 0.
\]
\end{prop}

The proof is based on the following easy property ($g$ is fixed throughout):
\[
\widetilde{V}{}^*(X'^q,g,p)^n_t = \frac1{\Delta_n^{1/4}k_n} \sum
_{i=0}^{[t/\Delta_n]-k_n}\Gamma(q)^n_i+R(q)^n_t,
\]
where
\begin{eqnarray*}
\Gamma(q)^n_i &=& |\overline{Z'^q}(g)^n_i|^p-\sum_{j=1}^{k_n-1}
|g^n_j|^p \Delta^n_{i+j}\Sigma(q),\qquad
\Sigma(q)_t = \sum_{s\leq t}|\Delta X'(q)_s|^p,\\
|R(q)^n_t| &\leq& \frac K{\Delta_n^{1/4}}
\biggl(\Sigma(q)_t \biggl(\frac{\overline{g}(p)_n}{k_n}-\overline{g}(p)\biggr)
+\Sigma(q)_{u_n}+\bigl(\Sigma(q)_t-\Sigma(q)
_{t-2u_n}\bigr) \biggr).
\end{eqnarray*}
\begin{lem}\label{L45} We can find a sequence $\eta_q$
going to $0$ as $q\to\infty$, with the following property: for
any $q\geq1$ and $i\geq1$ we have a decomposition
$\Gamma(q)^n_i=\Gamma'(q)^n_i+\Gamma''(q)^n_i$ where both $\Gamma
'(q)^n_i$ and
$\Gamma''(q)^n_i$ are $\mathcal F^n_{i+k_n}$-measurable and
%
%
\begin{equation}\label{CD15}
\left.
\begin{array}{l}
\mathbb{E}(|\Gamma'(q)^n_i|) \leq K_q\Delta_n^{1\wedge(p/4)}+\eta
_q\Delta
_n^{3/4},\\
\mathbb{E}(\Gamma''(q)^n_i\mid\mathcal F^n_i)=0,\\
\mathbb{E}(|\Gamma''(q)^n_i|^2) \leq K_q\Delta_n^{3/2}+\eta
_q\Delta_n.
\end{array}
\right\}
\end{equation}
\end{lem}
\begin{pf}
(1) Let us fix $i$, $q$ and $n$ which
will be left out in most notation below. We consider the filtration
$\mathcal F'_t=
\mathcal F_{i\Delta_n+t}$, and associated with this filtration the Brownian
motion $W'_t=W_{i\Delta_n+t}-W'_{i\Delta_n}$ and the Poisson random measure
$\underline{\mu}'((0,t]\times A)=\underline{\mu}((i\Delta
_n,i\Delta_n+t]\times A)$ whose
compensator is still $\underline{\nu}$. Recalling (\ref{CD3079}),
we set $b'_t=b^q_{i\Delta_n+t}$, and observe that $|b'_t|\leq Kq$ because
$b_t$ is bounded and $\int_{\{z\dvtx\gamma(z)>1/q\}}|\delta
(t,z)|\lambda
(dz)\leq
q\int\gamma^2(z)\lambda(dz)$. With all this notation
and (\ref{CD307}), we have
\[
X'^q_{i\Delta_n+t} = X'^q_{i\Delta_n}+\int_0^tb'_s\,ds+\int_0^t\sigma'_s\,dW'_s
+\bigl(\delta'1_{\{\gamma\leq1/q\}}\bigr)\star(\underline{\mu}'-\underline
{\nu})_t.
\]
Recalling $g_n$ in (\ref{LLN14}), we then set
%
%
\begin{eqnarray}\label{CD17}\quad
Y_t &=& \int_0^tb'_s g_n(s) \,ds+\int_0^t\sigma'_s g_n(s) \,dW'_s
+\bigl(\delta'g_n1_{\{\gamma\leq1/q\}}\bigr)\star(\underline{\mu
}'-\underline{\nu})_t
\nonumber\\[-8pt]\\[-8pt]
&&{}-\sum_{j=1}^{[t/\Delta_n]}g'^n_j\chi^n_{i+j-1}.\nonumber
\end{eqnarray}
Then by (\ref{RR2}) and (\ref{LLN14}) and (\ref{CD17}), we see that
$\overline{Z'(g})^n_i=Y_{u_n}$. If we further set
\[
Y'_t = \bigl(|\delta'g_n|^p 1_{\{\gamma\leq1/q\}}\bigr)\star\underline{\mu}'_t,
\]
we obtain
$Y'_{u_n}={\sum_{j=1}^{k_n-1}}|g^n_j|^p \Delta^n_{i+j}\Sigma(q)$. Hence
$\Gamma(q)^n_i=|Y_{u_n}|^p-Y'_{u_n}$.

For simplicity of notation we write $f(x)=|x|^p$ which is $C^2$
(recall $p>3$), and we associate the functions
\begin{eqnarray*}
F(x,y) &=& f(x+y)-f(x)-f'(x)y,\\
G(x,y) &=& f(x+y)-f(x)-f(y),\\
H(x,y) &=& F(x,y)-f(y),
\end{eqnarray*}
which clearly satisfy
%
%
\begin{equation}\label{CD19}
\left.
\begin{array}{l}
|F(x,y)| \leq K(|y|^p+y^2|x|^{p-2}),\\
|G(x,y)| \leq K(|x| |y|^{p-1}+|y| |x|^{p-1}),\\
|H(x,y)| \leq K(|x| |y|^{p-1}+y^2|x|^{p-2}).
\end{array}
\right\}
\end{equation}
Then we apply It\^o's formula and use (\ref{CD17}) to obtain
\[
|Y_t|^p-Y'_t = A_t+A'_t+N_t+N'_t,
\]
where
\begin{eqnarray*}
A_t &=& \int_0^ta_s\,ds,\qquad
A'_t = \sum_{j=1}^{[t/\Delta_n]}F(Y_{j\Delta_n-},-g'^n_j\chi^n_{i+j-1}),
\\
a_t &=& f'(Y_t)g_n(t)b'_t+\frac12 f''(Y_t)g_n(t)^2\sigma'^2_t\\
&&{} + \int_{\{z\dvtx\gamma(z)\leq1/q\}}H(Y_t,g_n(t)\delta'(t,z)) \lambda(dz)
\end{eqnarray*}
and $N_t$ is a martingale with angle bracket $C=\langle N,N\rangle$
given by
\begin{eqnarray*}
C_t &=& \int_0^tc_s\,ds,\\
c_t &=& f'(Y_t)^2g_n(t)^2\sigma'^2_t
+\int_{\{z\dvtx\gamma(z)\leq1/q\}}G(Y_t,g_n(t)\delta'(t,z))^2 \lambda(dz)
\end{eqnarray*}
and, finally,
\[
N_t' = -\sum_{j=1}^{[t/\Delta_n]}f'(Y_{j\Delta_n-})
g'^n_j\chi^n_{i+j-1},
\]
which is another martingale (because the $\chi_t$'s are centered) with
square bracket,
\[
C'_t := [N',N']_t = \sum_{j=1}^{[t/\Delta_n]}f'(Y_{j\Delta_n-})^2
(g'^n_j)^2(\chi^n_{i+j-1})^2.
\]

(2) The decomposition $\Gamma(q)^n_i=\Gamma'(q)^n_i+\Gamma''(q)^n_i$
is given by
\[
\Gamma'(q)^n_i = A_{u_n}+A'_{u_n},\qquad
\Gamma''(q)^n_i = N_{u_n}+N'_{u_n}.
\]
The $\mathcal F^n_{i+k_n}$-measurability of $\Gamma'(q)^n_i$ and
$\Gamma''(q)^n_i$
is obvious, as is the second part of (\ref{CD15}). The rest of (\ref{CD15})
will readily follow if we can find a sequence $\eta_q\to0$ such that
%
%
\begin{equation}\label{CD20}
\left.
\begin{array}{l}
\mathbb{E}(|A_{u_n}|) \leq K_q\Delta_n^{1\wedge(p/4)}+\eta
_q\Delta
_n^{3/4}, \qquad
\mathbb{E}(|A'_{u_n}|) \leq K_q\Delta_n,\\
\mathbb{E}(C_{u_n}) \leq K_q\Delta_n^{3/2}+\eta_q\Delta_n,\qquad
\mathbb{E}(C'_{u_n}) \leq K_q\Delta_n^{3/2}+\eta_q\Delta_n.
\end{array}
\right\}
\end{equation}

For this we need moment estimates for $Y_t$ as defined by
(\ref{CD17}). Recall $|b'|\leq Kq$ and $|\sigma'|\leq K$ and
$|g_n|\leq K$ and $|\delta(\cdot,z)|\leq\gamma(z)$ whereas
$\eta'_q=\int_{\{z\dvtx\gamma(z) \leq1/q\}}\gamma(z)^2 \lambda(dz)$
goes to $0$ as $q\to\infty$. In view of (SN-$2p$) and since
$|g'^n_j|\leq K\sqrt{\Delta_n}$, and using the
Burkholder--Davis--Gundy inequality
for the martingale which is the last term in (\ref{CD17}), we see that for
all $r\in(0,2p]$,
%
%
\begin{equation}\label{CD22}
\mathbb{E}(|Y_t|^r) \leq Kq^rt^r+Kt^{r/2}+K\eta'_qt^{1\wedge(r/2)}.
\end{equation}
By $f(x)=|x|^p$ and (\ref{CD19}), plus $p>3$, we see that
$|a_t|\leq K(q|Y_t|^{p-1}+|Y_t|^{p-2}+\eta'_q|Y_t|)$. Therefore (\ref{CD22})
yields $\mathbb{E}(|a_t|)\leq Kq^pt^{1\wedge(p/2-1)}+K\eta'_qt^{1/2}$.
In a similar way $c_t\leq K|Y_t|^{2p-2}+K\eta'_qY_t^2$; hence
$\mathbb{E}(c_t)\leq K(q^{2p-2}t^2+\eta'_qt)$. Then the
estimate for $A_{u_n}$ and $C_{u_n}$ in (\ref{CD20}) follows
upon taking $\eta_q=K\eta'_q$ for a $K$ large enough.

For the same reasons, plus (SN-$2p$), the $j$th summand in the
definition of $A'_t$ has an expectation smaller than $K\Delta_n^{p/2}+
K\Delta_n(q^{p-2}(j\Delta_n)^{p-2}+(j\Delta_n)^{p/2-1}+\eta
'_q(j\Delta_n))$ whereas
the $j$th summand in the expression for $C'_t$ has an expectation
smaller than
$K\Delta_n(q^{2p-2}(j\Delta_n)^{2p-2}+(j\Delta_n)^{p-1}+\eta
'_q(j\Delta_n))$.
The two other estimates in (\ref{CD20}) follow.
\end{pf}
\begin{pf*}{Proof of Proposition \protect\ref{PP10}}
In view of (\ref{SET4}) and of the fact that
$\mathbb{E}(\Sigma(q)_{s+u}-\Sigma(q)_s)\leq K_qu$, we deduce that
$R(q)^n_t\stackrel{\mathbb{P}}{\longrightarrow}0$
for all $q$. Hence it remains to prove that
%
%
\begin{equation}\label{CD14}
\lim_{q\to\infty} \limsup_n \mathbb{P} \Biggl(\frac1{\Delta
_n^{1/4}k_n}
\Biggl|\sum_{i=0}^{[t/\Delta_n]-k_n}
\Gamma(q)^n_i \Biggr|>\varepsilon\Biggr) = 0.
\end{equation}

We set, with $t$ fixed and the notation of the previous lemma,
\[
L'(q)_n = \frac1{\Delta_n^{1/4}k_n} \sum_{i=0}^{[t/\Delta
_n]-k_n}\Gamma'(q)^n_i,\qquad
L''(q)_n = \frac1{\Delta_n^{1/4}k_n} \sum_{i=0}^{[t/\Delta
_n]-k_n}\Gamma''(q)^n_i.
\]
The first property in (\ref{CD15}) yields $\mathbb{E}(|L'(q)_n|)\leq
K_q\Delta_n^{(1/4)\wedge(p/4-3/4)}+\eta_q$; hence since $p>3$,
\[
\lim_{q\to\infty} \limsup_n \mathbb{E}(|L'_n|) = 0.
\]
Next, the properties of $\Gamma''(q)^n_i$ in the
Lemma \ref{L45} imply that $|\mathbb{E}(\Gamma''(q)^n_i \Gamma''(q)^n_j)|$
vanishes when
$|j-i|>k_n$, and, otherwise, is smaller than
$K_q\Delta_n^{3/2}+\eta_q\Delta_n$. Hence
$\mathbb{E}((L''_n)^2)\leq K_q\Delta_n^{1/2}+\eta_q$, which yields
\[
\lim_{q\to\infty} \limsup_n \mathbb{E}(|L''_n|^2) = 0.
\]
Putting these two results together immediately yields
(\ref{CD14}).
\end{pf*}

\subsection[Proof of Theorem 4.4]{Proof of Theorem \protect\ref{TCD1}}\label{ssec:TCD1}

We start with the first claim, which easily follows from what precedes.
The family $(g_l)$ of weight functions is fixed. Since
$U(p,q)_t\stackrel{\mathbb{P}}{\longrightarrow}U(p)_t$ as $q\to
\infty$, the result is a trivial
consequence of (\ref{CD123}) and Propositions \ref{P37} and \ref{PP10}.

Next, we show that the second claim can be reduced
to the first claim. We take $p\geq4$, an even integer,\vspace*{-2pt} and it is enough
to prove that $\frac1{\Delta_n^{1/4}k_n} (\overline
{V}{}^*(g,p)^n-\widetilde{V}
{}^*(g,p)^n)\stackrel{\mathrm{u.c.p.}}{\longrightarrow}0$
for any weight\vspace*{1pt} function $g$. To see
this we observe that the difference $\overline
{V}{}^*(g,p)^n_t-\widetilde{V}{}^*(g,p)^n_t$
is a
linear combination of the processes (we omit to mention the function
$g$ below)
\[
\frac1{\Delta_n^{1/4}k_n} \sum_{i=0}^{[t/\Delta_n]-k_n}(\overline
{Z}^n_i)^{p-2r} (\widehat{Z}^n_i)^r
\]
for $r=1,\ldots,p/2$. So it enough to prove that, for
some $\rho>3/4$ and all $r=1,\ldots,p/2$,
%
%
\begin{equation}\label{CD304}
\mathbb{E} ((\overline{Z}{}^n_i)^{p-2r} (\widehat{Z}^n_i)^r
) \leq K\Delta_n^{\rho}.
\end{equation}

Hypothesis \ref{hypoSK} yields $\mathbb{E}(|\Delta^n_iX|^v\mid\mathcal F^n_{i-1})\leq
K\Delta_n^{(v\wedge2)/2}$ and (SN-$4p$) holds, so when $v\geq1$ we have
%
%
\begin{equation}\label{CD305}
\left.
\begin{array}{l}
\displaystyle\mathbb{E}(|\overline{X}{}^n_i|^v\mid\mathcal F^n_i) \leq K_v\Delta
_n^{(v/4)\wedge
(1/2)},\qquad
\mathbb{E}(|\widehat{X}^n_i|^v\mid\mathcal F^n_i) \leq K_v\Delta
_n^{1+v/2},\\
\displaystyle v\leq2p \quad\Rightarrow\quad
\mathbb{E} \Biggl( \Biggl|\sum_{j=1}^{k_n}(g^n_j)^2\Delta^n_{i+j}X\Delta
^n_{i+j}\chi
\Biggr|^v
\biggm|\mathcal F^n_i\Biggr) \leq K_v\Delta_n^{v/2+(v/2)\wedge1}.
\end{array}
\right\}\hspace*{-35pt}
\end{equation}
Now $(\overline{Z}{}^n_i)^{p-2r} (\widehat{Z}^n_i)^r$ is a linear
combination of terms of
the form
\[
a(u,v,w,s,t)^n_i = (\overline{X}{}^n_i)^u (\widehat
{X}^n_i)^v (\overline{\chi}^n_i)^w
(\widehat{\chi}_i^n)^s
\Biggl(\sum_{j=1}^{k_n}(g^n_j)^2\Delta^n_{i+j}X\Delta^n_{i+j}\chi
\Biggr)^t,
\]
where $u,v,w,s,t$ are integers with $u+w=p-2r$ and $v+s+t=r$. Using
H\"older's inequality, and taking advantage of (\ref{CD305}) and of
$\mathbb{E}(|\overline{\chi}^n_i|^l)\leq K_r\Delta_n^{l/4}$ and
$\mathbb{E}(|\widehat{\chi}^n_i|^l)\leq K_r\Delta_n^{l/2}$, we see
that for all $u',v',w',s',t'\geq0$ such that $u'+v'+w'+s'+t'=1$ and
$u'=0$ (resp., $v'=0$, $w'=0$, $s'=0$, $t'=0$) if and only if
$u=0$ (resp., $v=0$, $w=0$, $s=0$, $t=0$), and also
$\frac w{w'}\vee\frac{2s}{s'}\vee\frac t{t'}\leq2p$
(which is possible because $w+2s+t\leq p$), we have
$\mathbb{E}(|a(u,v,w,s,t)^n_i|)\leq K\Delta_n^\rho$ where
\begin{eqnarray*}
\rho &=& \frac u4\wedge\frac{u'}2+v'1_{v>0}+\frac v2+\frac w4+\frac s2+
\frac t2+t'\wedge\frac t2 \\
&=&
\frac r2+\frac w4+\frac u4\wedge\frac{u'}2+v'1_{v>0}+ t'\wedge\frac t2.
\end{eqnarray*}
Then $\rho>3/4$ as soon as $r\geq2$, or $r=1$ and $w\geq1$. The only
other case is $r=1$ and $w=0$, so $u=p-2\geq2$ and we have
\[
\rho= \frac12+\frac{u'}2+v'1_{v>0}+t'\wedge\frac t2.
\]
Then we have three sub-cases:

(1) $v=1$, hence $t=t'=s=s'=w'=0$ and $\rho=\frac{1+u'}2+v'$ with
the condition $u'+v'=1$, so $u'=v'=1/2$ yields $\rho>3/4$;

(2) $s=1$, hence $t=t'=v=v'=w'=0$ and $\rho=\frac{1+u'}2$ with
the conditions $u'+s'=1$ and $s'p\geq1$, so $s'=1/3$ yields $\rho>3/4$;

(3) $t=1$, hence $v=v'=s=s'=w'=0$ and $\rho=\frac12+\frac{u'}2+
t'\wedge\frac12$ with the condition $u'+t'11$ and $2t'p\geq1$,
so $u'=t'=1/2$ yield $\rho>3/4$.

Hence in all cases (\ref{CD304}) holds with some $\rho>3/4$,
and the proof is complete.

\subsection[Proof of Theorem 4.6]{Proof of Theorem \protect\ref{TVQ1}}\label{ssec-VQ}

Here again the proof will be divided into several steps, and before
proceeding we observe two preliminary facts. First, that
$\overline{\mu}_4(g,g;\eta,\zeta)$ takes the form (\ref{VQ3})
results from a
tedious but elementary calculation. Second, by localization we may assume
(SN-$4$) and Hypothesis \ref{hypoSH}.

We omit the mention of the function $g$ in $\overline{Y}{}^n_i$ and
$\widehat{Y}^n_i$. We
generally use the notation of the proof of Theorem \ref{TCD1}, and in
particular the stopping times $T(q,m)$ and $T_m$ introduced in Section
\ref{ssec:43}, the processes of (\ref{CD307}), the sets $\Omega_n(t,q)$
satisfying (\ref{CD121}) and the (random) integers $I^n_m$.
In the sequel, we will vary the process $X$ (but not
the noise process $\chi$), so the process $\overline{V}{}^n$ of (\ref
{VQ1}) will
be denoted by $\overline{V}(X)^n$. We also write $U'(\sigma)_t$ and
$U(2,\sigma,\delta)_t$
for the two terms in (\ref{VQ2}), and $\overline{U}(\sigma,\delta
)_t$ for
their sum,
to emphasize their dependency on the process $\sigma$ and the function
$\delta$ (through the jumps of $X$, for the latter).

\textit{Step} 1. In this step we prove the result when, in
addition to
Hypothesis \ref{hypoSH}, we have
%
%
\begin{eqnarray}\label{VQ11}
&&\gamma(z)\leq1/q \quad\Rightarrow\quad\delta(\omega,t,z)=0,
\\
%
%
\label{VQ10}
&&\left.
\begin{array}{l}
\displaystyle b'_s = b_s-\int\delta(t,x)1_{\{|\delta(y,z)|\leq1\}} \lambda(dz) =
\sum_{r\geq1}b_{S_r}1_{[S_r,S_{r+1})}(t),\\
\displaystyle\sigma_s = \sum_{r\geq1}\sigma_{S_r}1_{[S_r,S_{r+1})}(t),
\end{array}
\right\}
\end{eqnarray}
for some $q\geq1$ and
a sequence of stopping times $S_r$, increasing to
$\infty$ and with $S_0=0$.

(1) Under (\ref{VQ11}) and (\ref{VQ10}) we have $X^q_t=\sum_{s\leq
t}\Delta
X_s$, and $X'=X'^q$ is the continuous process given by the right-hand side
of (\ref{CIS}) with $b'$ instead of $b$. Similarly to (\ref{CD123}),
we have on $\Omega_n(t,q)$,
%
%
\begin{eqnarray}\label{VQ4}\qquad
\overline{V}(X)_t^n &=& \overline{V}(X')^n_t+Y^n_t-\frac12 Y'^n_t,
\nonumber\\
Y^n_t &=& \sum_{m\in P_q\dvtx T_m\leq t}\zeta^n_m,\qquad
Y'^n_t = \sum_{m\in P_q\dvtx T_m\leq t}\zeta'^n_m,\nonumber\\
\zeta^n_m &=& \frac1{\Delta_n^{1/4}k_n} \Biggl(\sum_{j=1}^{k_n-1}
\bigl( |\overline{(X'+\chi)}^n_{I^n_m+1-j}+g^n_j\Delta
X_{T_m} |^2
\nonumber\\[-8pt]\\[-8pt]
&&\hspace*{62.3pt}{}  - |\overline{(X'+\chi)}^n_{I^n_m+1-j} |^2
- |g^n_j\Delta X_{T_m} |^2 \bigr)\nonumber\\
&&\hspace*{97.6pt}{} +\bigl(\overline
{g}(2)_n-k_n\overline
{g}(2)\bigr)|\Delta
X_{T_m}|^2 \Biggr),\nonumber\\
\zeta'^n_m &=& \frac1{\Delta_n^{1/4}k_n}\sum_{j=1}^{k_n}(g'^n_j)^2
\bigl((\Delta X_{T_m})^2+2\Delta X_{T_m}\Delta^n_{I^n_m+1-j}(X'+\chi
) \bigr).\nonumber
\end{eqnarray}
Let $(\mathcal H_t)$ be the filtration defined in the proof of Lemma
\ref{LCD1} and associated with our $q$. The same argument used in
that lemma shows $\mathbb{E}(|\Delta^n_{I^n_m+1-j}(X'+\chi)|\mid
\mathcal H_0)\leq K$
whereas $|\Delta X_{T_m}|\leq K$ by Hypothesis \ref{hypoSH}. It follows that
$\mathbb{E}(|\zeta'^n_m|)\leq K\Delta_n^{3/4}$; hence
%
%
\begin{equation}\label{VQ5}
Y'^n_t \stackrel{\mathbb{P}}{\longrightarrow} 0.
\end{equation}

(2) Next we prove the (functional) stable convergence $\overline
{V}(X')^n\stackrel{\mathcal L-(s)}{\longrightarrow}
U'(\sigma)$. This looks the same as Theorem
\ref{TCC1} for $p=2$, however we do not have Hypothesis \ref{hypoK} here. Now a look
at the proof of this theorem shows that Hypothesis \ref{hypoK} [instead of Hypothesis \ref{hypoH}] is used in
two places only, namely for (\ref{EXT}) for $k=2$, and in Lemma \ref{LCC5}.
Here, the proof of Lemma \ref{LCC5} proceeds in an obvious way under
(\ref{VQ10}), and we are left to show that (\ref{EXT})
holds when $k=2$.

The variable $\Delta_n^{3/4-p/2}U_t^{n,2}$ for $p=2$
is the sum $\sum_{i=1}^{i_n(m,t)}\sum_{j=0}^{mk_n-1}\theta^n_{i,j}$, where
\[
\theta^n_{i,j} = \Delta_n^{1/4} \mathbb{E} \bigl(\phi(X'+\chi
,g,2)_{I(m,n,i),j}
-\phi(g,2)_{I(m,n,i)+j}\mid\mathcal F(m)^n_{i-1} \bigr).
\]
Let $J_n$ be the set of all $i$ such that $(i-1)(m+1)u_n<S_r\leq
imu_n$ for some $r\geq1$ (the indices of those ``big blocks''
that contain at least one $S_r$), and consider the two processes
\[
A^n_t = \sum_{i\in\{1,\ldots,i_n(m,t)\}\cap J_n}
\sum_{j=0}^{mk_n-1}\theta^n_{i,j},\qquad
A'^n_t = \sum_{i\in\{1,\ldots,i_n(m,t)\}\cap J_n^c}
\sum_{j=0}^{mk_n-1}\theta^n_{i,j}.
\]
Applying (\ref{L4N8}) with $u=1$ [recall (SN-$4$)], we obtain
$\mathbb{E}(|\theta^n_{i,j}|)\leq K\Delta_n^{3/4}$. Therefore
$\mathbb{E}(\sup_{s\leq
t\wedge S_r}
|A^n_s|)$ is obviously smaller than $Kr\Delta_n^{1/4}$ and, since
$S_r\to\infty$ as $r\to\infty$, we deduce $A^n\stackrel{\mathrm
{u.c.p.}}{\longrightarrow}0$, and it
remains to
prove the same for $A'^n$.

For this, and reproducing the proof of Lemma \ref{LCC3}, we observe that
Hypothesis~\ref{hypoK} comes in only to decompose the variables $\overline{\lambda
}^n_{i+j}$ as
$\xi^n_{i,j}+\xi'^n_{i,j}$. We easily~deduce from
(\ref{VQ10}) that when $i\notin J_n$ such a decomposition holds with
$\xi^n_{i,j}=0$ and $\xi'^n_{i,j}=b'_{i\Delta_n}\Delta_n$. Then
the\vspace*{1pt}
original proof
goes through to show that $A'^n\stackrel{\mathrm
{u.c.p.}}{\longrightarrow}0$, and thus (\ref{EXT}) for $k=2$
holds here.

(3) We have $\overline{V}(X')^n\stackrel{\mathcal
L-(s)}{\longrightarrow}U'(\sigma)$ from what precedes, and this gives
the result (functional stable convergence in law) when $X$ is
continuous, in addition to satisfying (\ref{VQ10}). When $X$\vspace*{-1pt} has
jumps, the proof of Proposition \ref{P37} is valid
when $p=2$ (it only supposes the $C^2$ property of $x\mapsto|x|^p$), so
$Y^n_t\stackrel{\mathcal L-(s)}{\longrightarrow}U(2,\sigma,\delta
)_t$ (for $t$ fixed, not
functional convergence).

Now, exactly as in the proof of Lemma 5.8 in \cite{J2}, one can
show that we have the \textit{joint} stable convergence in law in
Proposition \ref{PCC1} and Lemma \ref{LCD1} which results in the
joint convergence
\[
(\overline{V}(X')^n_t,Y^n_t) \stackrel{\mathcal
L-(s)}{\longrightarrow} (U'(\sigma)_t,U(2,\sigma,\delta)_t).
\]
Then we easily deduce from (\ref{CD121}),
(\ref{VQ4}) and (\ref{VQ5}) that $\overline{V}(X)^n_t\stackrel
{\mathcal L-(s)}{\longrightarrow}\overline{U}(\sigma
,\delta)_t$.

\textit{Step} 2. We turn to the general case, and we begin by
constructing an approximation of $X$ satisfying (\ref{VQ11}) and
(\ref{VQ10}).

For $q\geq1$ we recall the process $b^q$ of (\ref{CD3079}).
If further $r\geq1$ we denote by
$S(q,r)_r$ the strictly increasing rearrangement of the points
in the set $\{k2^{-r}\dvtx k\geq0\}\cup\{T(q,m)\dvtx m\geq1\}$.
By a classical density argument there are adapted processes
$b(q,r)$ and $\sigma(q,r)$ with the following properties: they are
bounded by the same bounds as $b^q$ and $\sigma$,
respectively, constant over each
interval $[(k-1)2^{-r},k2^{-r})$ for $b(q,r)$ and each interval
$[S(q,r)_{k-1},S(q,r)_q)$ for $\sigma$ and such that for all
$q,m\geq1$ and $t\geq0$,
%
%
\begin{equation}\label{VQ15}\qquad
r\to\infty\quad\Rightarrow\quad\cases{
\displaystyle\varepsilon(q,r)_t = \mathbb{E} \biggl(\int_0^t\bigl(|b(q,r)_s-b^q_s|^2
\cr
\hspace*{77.8pt}{}+|\sigma(q,r)_s-\sigma_s|^2\bigr)\,ds \biggr) \to0,\cr
\displaystyle\sigma(q,r)_{T(q,m)} = \sigma_{T(q,m)},\cr
\sigma(q,r)_{T(q,m)-} \to\sigma_{T(q,m)-},}
\end{equation}
(we use here the c\`adl\`ag property of $\sigma$).
Next, we introduce the following family of processes:
%
%
\begin{equation}\label{VQ8}
\left.
\begin{array}{l}
\displaystyle X(q,r)_t = X_0+\int_0^tb(q,r)_s\,ds+\int_0^t\sigma(q,r)_s\,dW_s
+\bigl(\delta1_{\{\gamma>1/q\}}\bigr)*\underline{\mu},\\
\displaystyle X'(q,r)_t=X_t-X(q,r)_t\\
\hspace*{42pt}=
\displaystyle\int_0^t\bigl(b^q_s-b(q,r)_s\bigr)\,ds+\int_0^t\bigl(\sigma_s-\sigma(q,r)_s\bigr)\,dW_s+M^q_t
\end{array}
\right\}\hspace*{-42pt}
\end{equation}
[here $M^q$ is given by (\ref{CD307})]. Finally, another notation will be
\begin{eqnarray*}
\varepsilon(q,r)^n_i &=& \mathbb{E} \biggl(\int_{i\Delta_n}^{i\Delta
_n+u_n} \bigl(|b(q,r)_s-b^q_s|^2
+|\sigma(q,r)_s-\sigma_s|^2\bigr)\,ds \biggr),\\
\varepsilon_q &=& \int_{\{z\dvtx\gamma(z)\leq1/q\}}\gamma(z)^2 \lambda(dz).
\end{eqnarray*}

By construction $X(q,r)$ satisfies (\ref{VQ11}) and (\ref{VQ10}), so
Step 1 gives
\[
\overline{V}(X(q,r))^n_t \stackrel{\mathcal L-(s)}{\longrightarrow
} \overline{U}(\sigma(q,r),\delta(q))_t
\]
for any $t$ and $q,r\geq1$, and where
$\delta(q)(\omega,t,z)=\delta(\omega,t,z) 1_{\{\gamma(z)>1/q\}}$,
and the convergence even holds in the functional sense when $X$ is
continuous.

Note that, since $\sigma$ and $\sigma(q,r)$ and $\alpha$ are
uniformly bounded
and the function $\overline{\mu}_4$ in (\ref{VQ3}) is locally
Lipschitz in
$(\eta,\zeta)$, we have
\[
\mathbb{E} \Bigl(\sup_{s\leq t}|U'(\sigma)_s-U'(\sigma
(q,r))_s|^2 \Bigr)
\leq K \mathbb{E} \biggl(\int_0^t|\sigma_s-\sigma(q,r)_s|^2 \,ds \biggr)
\leq K\varepsilon(q,r)_t.
\]
On the other hand, since $\delta(q)$ is bounded, it follows from
(\ref{CD6}) that
\begin{eqnarray*}
&&\mathbb{E} \Bigl({\sup_{s\leq t}}|U(2,\sigma,\delta(q))_s-U(2,\sigma
(q,r),\delta(q))_s|^2 \Bigr)
\\
&&\qquad\leq K\mathbb{E} \biggl(\sum_{m\geq1}\bigl|\sigma_{T(q,m)-}-\sigma
(q,r)_{T(q,m)-}\bigr|^2
1_{\{T(q,m)\leq t\}} \biggr),
\end{eqnarray*}
which goes to $0$ as $r\to\infty$ by (\ref{VQ15}). Furthermore,
\[
\mathbb{E} \Bigl({\sup_{s\leq t}}|U(2,\sigma,\delta(q))_s-U(2,\sigma
,\delta
)_s|^2 \Bigr)
\leq K\mathbb{E} \biggl(\sum_{s\leq t}|\Delta X_s|^2 1_{\{|\Delta
X_s|\leq1/q\}
} \biggr),
\]
which goes to $0$ as $q\to\infty$. Summarizing those results,
we end up with
\[
\lim_{q\to\infty} \limsup_{r\to\infty}
\mathbb{E} \Bigl(\sup_{s\leq t}|\overline{U}(\sigma,\delta
)_s-\overline{U}(\sigma
(q,r),\delta(q))_s|^2 \Bigr)
= 0.
\]
Therefore, in order to get our theorem it remains to prove that
for all $t,\eta>0$ we have, where C refers to the case, $X$ is continuous
and D to the general (discontinuous) case,
%
%
\begin{equation}\label{VQ16}\qquad\quad
\left.
\begin{array}{l}
\mbox{C: }
\displaystyle\lim_{q\to\infty}\limsup_{r\to\infty}\limsup_{n\to\infty}
\mathbb{P} \Bigl({\sup_{s\leq t}}|\overline{V}(X(q,r))^n_s-\overline
{V}(X)^n_s|>\eta\Bigr)=0,\\[8pt]
\mbox{D: }
\displaystyle\lim_{q\to\infty}\limsup_{r\to\infty}\limsup_{n\to\infty}
\mathbb{P} \bigl(|\overline{V}(X(q,r))^n_t-\overline{V}(X)^n_t|>\eta
\bigr)=0.
\end{array}
\right\}
\end{equation}

\textit{Step} 3. If $Z(q,r)=X(q,r)+\chi$, we have
\begin{eqnarray*}
&&\phi(Z,g,2)^n_i-\phi(Z(q,r),g,2)^n_i
\\
&&\qquad=(\overline{X}{}^n_i)^2-(\overline{X(q,r)}{}^n_i)^2+2\overline{\chi
}^n_i\bigl(\overline{X}{}^n_i-
\overline{X(q,r)}{}^n_i\bigr)-
\frac12 v^n_i,\\
&&\hspace*{8.89pt}v^n_i=
\sum_{j=1}^{k_n}(g'^n_j)^2 \bigl((\Delta^n_{i+j}X)^2-(\Delta^n_{i+j}X(q,r))^2
\\
&&\hspace*{76.99pt}{}
+2\Delta^n_{i+j}\chi
\bigl(\Delta^n_{i+j}X-\Delta^n_{i+j}X(q,r)\bigr) \bigr).
\end{eqnarray*}
Therefore
\[
\overline{V}(X)^n_t-\overline
{V}(X(q,r))^n_t = G^1(q,r)^n_t+G^2(q,r)^n_t-\tfrac12 V^n_t,
\]
where
\begin{eqnarray*}
V^n_t &=& \frac1{k_n\Delta_n^{1/4}}\sum
_{i=0}^{[t/\Delta_n]-k_n}v^n_i,\\
G^1(q,r)^n_t&=&\frac1{\Delta_n^{1/4}} \Biggl(
\frac1{k_n}\sum_{i=0}^{[t/\Delta_n]-k_n}
\bigl((\overline{X}{}^n_i)^2-(\overline{X(q,r)}{}^n_i)^2 \bigr)\\
&&\hspace*{29.5pt}{}-\overline{g}(2) \bigl([X,X]_t-[X(q,r),X(q,r)]_t \bigr) \Biggr),\\
G^2(q,r)^n_t &=& \frac2{k_n\Delta_n^{1/4}}
\sum_{i=0}^{[t/\Delta_n]-k_n}\overline{\chi}^n_i
\bigl(\overline{X}{}^n_i-\overline{X(q,r)}{}^n_i \bigr).
\end{eqnarray*}
We obviously have $\mathbb{E}(|v^n_i|)\leq K\Delta_n$, so
$V^n\stackrel{\mathrm{u.c.p.}}{\longrightarrow}0$.
Therefore, instead of (\ref{VQ16}), we are left to prove for $l=1,2$,
%
%
\begin{equation}\label{VQ17}
\left.
\begin{array}{l}
\mbox{C: }
\displaystyle\lim_{q\to\infty}\limsup_{r\to\infty}\limsup_{n\to\infty}
\mathbb{P} \Bigl({\sup_{s\leq t}}|G^l(q,r)^n_s|>\eta\Bigr)=0,\\[8pt]
\mbox{D: }
\displaystyle\lim_{q\to\infty}\limsup_{r\to\infty}\limsup_{n\to\infty}
\mathbb{P} \bigl(|G^l(q,r)^n_t|>\eta\bigr)=0.
\end{array}
\right\}
\end{equation}

\textit{Step} 4. We begin by proving (\ref{VQ17}) for $l=2$.
We split the sum in the definition of $G^2(q,r)^n_t$ into two parts:
$G^3(q,r)^n_t$ is the sum over those $i$'s such that the fractional part of
$i/2k_n$ is in $[0,1/2)$, and $G^4(q,r)^n_t$ which is the sum
when the fractional part is in $[1/2,1)$, so it enough to show
(\ref{VQ17}) for $l=3$ and $l=4$, and we will do it for
$l=3$ only. We have
\[
\left.
\begin{array}{l}
\displaystyle G^3(q,r)^n_t=\sum_{j=0}^{J_n+1}\zeta(q,r)^n_i,\\
\displaystyle\zeta(q,r)^n_j = \frac2{k_n\Delta_n^{1/4}}
\sum_{i=2jk_n}^{(2jk_n+k_n-1)\wedge([t/\Delta_n]-k_n)}
\overline{\chi}^n_i \bigl(\overline{X}{}^n_i-\overline{X(q,r)}{}^n_i \bigr),
\end{array}
\right\}
\]
where $J_n$ is the integer part of $([t/\Delta_n]+1-2k_n)/2k_n$
[$J_n$ depends on $t$, and all $\zeta(q,r)^n_j$ have $k_n$ summands,
except the $J_n$th one which may have less].
Note that $\zeta(q,r)^n_j$ is $\mathcal F^n_{2(j+1)k_n}$-measurable,
and by successive conditioning we have $\mathbb{E}(\zeta(q,r)^n_j\mid
\mathcal F^n_{2jk_n})=0$. Therefore by a martingale argument
(\ref{VQ17}) will follow if we prove
%
%
\begin{equation}\label{VQ21}
\lim_{q\to\infty} \limsup_{r\to\infty} \limsup_{n\to\infty}
\mathbb{E} \Biggl({\sum_{j=0}^{J(n,t)}}|\zeta(q,r)^n_j|^2 \Biggr) = 0.
\end{equation}

Now, recall (\ref{VQ8}) and (\ref{VQ11}).
Then, by (\ref{LLN14}) and standard estimates, plus (\ref{VQ15}) and
the Cauchy--Schwarz inequality, plus (\ref{LLN141}) and successive
conditioning, we get
\[
\mathbb{E} \bigl((\overline{\chi}^n_i)^2 \bigl(\overline{X}{}^n_i-\overline
{X(q,r)}{}^n_i\bigr)^2 \bigr) \leq
K\Delta_n^{1/2}\bigl(\varepsilon(q,r)^n_i+u_n\varepsilon_q\bigr)
\]
and so the expectation in (\ref{VQ21}) is smaller than
$K(\varepsilon(q,r)_t+\varepsilon_q)$. Hence (\ref{VQ21}) holds.

\textit{Step} 5. Now we turn to $l=1$ in (\ref{VQ17}). We
write $G^1(q,r)^n_t=G^5(q,r)^n_t+G^6(q,r)^n_t$
where, with the notation $A(q,r)=[X,X]-[X(q,r),X(q,r)]$,
\begin{eqnarray*}
&&\left.
\begin{array}{l}
\displaystyle G^5(q,r)^n_t =
\sum_{i=0}^{[t/\Delta_n]-k_n}\vartheta(q,r)^n_i,\\[12pt]
\displaystyle\vartheta(q,r)^n_i = \frac1{k_n\Delta_n^{1/4}} \biggl(
(\overline{X}{}^n_i)^2-(\overline{X(q,r)}{}^n_i)^2\\[12pt]
\hspace*{90pt}{}-\displaystyle\int_{i\Delta_n}^{i\Delta_n+u_n} g_n(s-i\Delta_n)^2\, dA(q,r)_s\biggr),
\end{array}
\right\}
\\
&&G^6(q,r)^n_t = \frac1{\Delta_n^{1/4}} \Biggl(\frac1{k_n}\sum
_{i=0}^{[t/\Delta_n]-k_n}
\int_{i\Delta_n}^{i\Delta_n+u_n}g_n(s-i\Delta
_n)^2\,dA(q,r)_s\\
&&\hspace*{216.8pt}{}-\overline{g}(2)
A(q,r)_t \Biggr).
\end{eqnarray*}

In this step we prove that $G^6(q,r)^n$ satisfies (\ref{VQ17}).
A simple calculation shows that [recall the notation $\overline{g}(2)_n$
of (\ref{3})]
\[
G^6(q,r)^n_t = \frac1{\Delta_n^{1/4}}\int_0^t \biggl(\frac{\overline
{g}(2)_n}{k_n}-
\overline{g}(2) \biggr)\,dA(q,r)_s+v(q,r)^n_t,
\]
where because of (\ref{SET4}) the remainder term $v(q,r)^n_t$ satisfies
with $A'(q,r)$ being the variation process of $A(q,r)$,
\[
|v(q,r)^n_t| \leq\frac K{\Delta_n^{1/4}} \bigl(
A'(q,r)_{u_n}+\bigl(A'(q,r)_t-A'(q,r)_{t-2u_n}\bigr) \bigr).
\]
In the continuous Case C, we have $A'(q,r)_{s+u_n}-A'(q,r)_s\leq Ku_n$,
hence $\sup_{s\leq t}|v(q,r)^n_s|\leq K\Delta_n^{1/4}$. In the
discontinuous\vspace*{-2pt} Case D we only have
$\mathbb{E}(A'(q,\break r)_{s+u_n}-A'(q,r)_s)\leq Ku_n$ so that $v(q,r)^n_t
\stackrel{\mathbb{P}}{\longrightarrow}0$ as $n\to\infty$. Then if
we apply (\ref{SET4}) we
obtain (\ref{VQ17}) for $l=6$.

\textit{Step} 6. It remains to prove (\ref{VQ17}) for
$l=5$. For this we use (\ref{LLN14}) again and It\^o's formula
to get, with $Y^{n,i}_t=\int_{i\Delta_n}^tg_n(s-i\Delta_n)\,dY_s$ for
any semimartingale $Y$ and for $t\geq i\Delta_n$,
\begin{eqnarray*}
&&
(\overline{X}{}^n_i)^2- \int_{i\Delta_n}^{i\Delta
_n+u_n}g_n(s-i\Delta_n)^2 \,d[X,X]_s\\
&&\qquad=2 \int_{i\Delta_n}^{i\Delta_n+u_n}X^{n,i}_s g_n(s-i\Delta_n)(
b^q_s \,ds+\sigma_s \,dW_s)\\
&&\qquad\quad{} + 2\int_{i\Delta_n}^{i\Delta_n+u_n}
X^{n,i}_{s-} \,dM^q_s\\
&&\qquad\quad{}+2\int_{i\Delta_n}^{i\Delta_n+u_n}\int_{\{\gamma(z)>1/q\}}
X^{n,i}_{s-} g_n(s-i\Delta_n) \delta(s,z) \underline{\mu}(ds,dz)
\end{eqnarray*}
and a similar expression with $(X,b^q,\sigma,\delta)$ substituted with
$(X(q,r),b(q,r)$,\break $\sigma(q,r),\delta(q))$, so the second term on the
right-hand side above vanishes in this case [remember the last part
of (\ref{VQ8})]. Therefore,
\[
\vartheta(q,r)^n_i =
\frac2{k_n\Delta_n^{1/4}}\sum_{j=1}^6\eta(q,r,j)^n_i,
\]
where, using (\ref{CD3079}) and with the notation
$I(n,i)=(i\Delta_n,i\Delta_n+u_n]$, we have
\begin{eqnarray*}
\eta(q,r,1)^n_i&=&\int_{I(n,i)}X'(q,r)^{n,i}_s g_n(s-i\Delta_n)
\,ds\, \biggl(b_s+\int_{\{|\delta(s,z)|>1\}}\delta(s,z)\lambda(dz)
\biggr),\\
\eta(q,r,2)^n_i &=& \int_{I(n,i)}X(q,r)^{n,i}_s g_n(s-i\Delta_n)
\bigl(b^q_s-b(q,r)_s\bigr)\,ds,\\
\eta(q,r,3)^n_i &=& \int_{I(n,i)}X'(q,r)^{n,i}_s g_n(s-i\Delta_n)
\sigma_s\,dW_s,\\
\eta(q,r,4)^n_i &=& \int_{I(n,i)}X(q,r)^{n,i}_s g_n(s-i\Delta_n)
\bigl(\sigma_s-\sigma(q,r)_s\bigr)\,dW_s,\\
\eta(q,r,5)^n_i &=& \int_{I(n,i)}X^{n,i}_{s-} g_n(s-i\Delta_n) \, dM^q_s,
\\
\eta(q,r,6)^n_i &=& \int_{I(n,i)}\int_{\gamma(z)>1/q}
X'(q,r)^{n,i}_{s-} g_n(s-i\Delta_n) \delta(s,z) (\underline{\mu
}-\underline{\nu})(ds,dz).
\end{eqnarray*}
Therefore, since $\eta(q,r,j)^n_i$ for $j=3,4,5,6$ are
martingale increments, (\ref{VQ17}) for $l=5$ will follow if we prove
that for all $t>0$, and as $m\to\infty$,
%
%
\begin{eqnarray}
\label{VQ22}\qquad
&&j=1,2 \nonumber\\[-8pt]\\[-8pt]
&&\quad\Rightarrow\quad
\lim_{q\to\infty} \limsup_{r\to\infty}
\limsup_{n\to\infty} \Delta_n^{1/4} \mathbb{E} \Biggl({\sum
_{i=0}^{[t/\Delta_n]-k_n}}|\eta
(q,r,j)^n_i| \Biggr) \to0,\nonumber
\\
\label{VQ23}
&&j=3,4,5,6\nonumber\\[-8pt]\\[-8pt]
&&\quad\Rightarrow\quad
\lim_{q\to\infty} \limsup_{r\to\infty}
\limsup_{n\to\infty} \Delta_n^{1/2} \mathbb{E} \Biggl({\sum
_{i=0}^{[t/\Delta_n]-k_n}}|
\eta(q,r,j)^n_i|^2 \Biggr) \to0.\nonumber
\end{eqnarray}

Then, standard estimates yield for $s\in I(n,i)$ and $p\geq2$
[recall $|b^q_t|+|b(q,r)_t|\leq Kq$],
\begin{eqnarray*}
\mathbb{E} \Bigl({\sup_{t\leq s}}|X'(q,r)^{n,i}_t|^2 \Bigr) &\leq& K\bigl(
\varepsilon(q,r)^n_i+\Delta_n^{1/2} \varepsilon_q\bigr),\\
\mathbb{E} \Bigl({\sup_{t\leq s}}|X(q,r)^{n,i}_t|^p \Bigr) &\leq& K_p
(q^p\Delta_n^{p/2}
+\Delta_n^{1/2} ),\\
\mathbb{E} \Bigl({\sup_{t\leq s}}|X^{n,i}_t|^2 \Bigr) &\leq& K\Delta_n^{1/2}
\end{eqnarray*}
and it follows that, since $|g_n|\leq K$ and $\varepsilon(q,r)^n_i\leq
K$ and
$\varepsilon_q\leq K$ and
\[
\int_{\{|\delta(s,z)|>1\}}|\delta
(s,z)|\lambda(dz)
\leq\int\gamma(z)^2\lambda(dz)<\infty,
\]
\begin{eqnarray*}
&&j=1,2 \quad\Rightarrow\quad
\mathbb{E}(|\eta(q,r,j)^n_i|) \leq
K\Delta_n^{1/2} \bigl(q\sqrt{\varepsilon(q,r)^n_i}
+\Delta_n^{1/4}\sqrt{\varepsilon_q} \bigr)\\
&&j=3,4,5,6 \quad\Rightarrow\quad\mathbb{E}(|\eta(q,r,j)^n_i|^2) \leq
K\Delta_n^{1/2} \bigl(q^2\Delta_n^{3/4}+\varepsilon(q,r)^n_i+
\Delta_n^{1/2}\varepsilon_q \bigr).
\end{eqnarray*}

By H\"older's inequality,
\[
\Biggl(\Delta_n^{3/4}\sum_{i=0}^{[t/\Delta_n]-k_n}\sqrt{\varepsilon
(q,r)^n_i} \Biggr)^2 \leq
\Delta_n^{1/2}\sum_{i=0}^{[t/\Delta_n]-k_n}\varepsilon
(q,r)^n_i \leq K\varepsilon(q,r)_t.
\]
Since $\varepsilon(q,r)\to0$ as $r\to\infty$, for each $q$, whereas
$\varepsilon_q\to0$ as $q\to\infty$. Therefore we readily
obtain (\ref{VQ22}) and (\ref{VQ23}), and the proof is finished.

\printaddresses

\end{document}